\documentclass[11pt,a4paper, twoside]{article}
\usepackage[notref,notcite]{showkeys}

\usepackage{amsmath,amssymb,amsthm,amsfonts,mathrsfs,amscd,environ}
\usepackage{dsfont}
\usepackage{latexsym,enumerate,color,geometry,extarrows}
\usepackage{xcolor}
\usepackage{verbatim,fancyhdr,enumitem}
\usepackage{hyperref}
\usepackage{url}
\hypersetup{
hidelinks,
hypertex=true,
colorlinks=true,
linkcolor=red,
citecolor=blue,
urlcolor=blue
}
 \NewEnviron{ews}{%
\begin{equation}\begin{split}
  \BODY
\end{split}\end{equation}
}

\NewEnviron{ews*}{%
\begin{equation*}\begin{split}
  \BODY
\end{split}\end{equation*}
}

\def\beg{\begin}
\def\bequ{\begin{equation}}
\def\enqu{\end{equation}}
\def\bes{\begin{split}}
\def\ens{\end{split}}
\def\bews{\begin{ews}}
\def\beqn{\begin{eqnarray}}
\def\enqn{\end{eqnarray}}
\def\beq*{\begin{equation*}}
\def\enq*{\end{equation*}}
\def\bqn*{\begin{eqnarray*}}
\def\eqn*{\end{eqnarray*}}
\def\bary{\begin{array}}
\def\eary{\end{array}}
\def\bpma{\begin{pmatrix}}
\def\epma{\end{pmatrix}}
\def\bvma{\begin{Vmatrix}}
\def\evma{\end{Vmatrix}}

 \numberwithin{equation}{section}

\def\al{\alpha}
\def\be{\beta}
\def\ga{\gamma}
\def\de{\delta}
\def\ep{\epsilon}

\def\ze{\zeta}
\def\et{\eta}
\def\th{\theta}

\def\la{\lambda}
\def\rh{\rho}

\def\si{\sigma}
\def\ta{\tau}
\def\ph{\phi}

\def\ps{\psi}
\def\om{\omega}

\def\Ga{\Gamma}

\def\Th{\Theta}
\def\Si{\Sigma}

\def\Om{\Omega}

\def\R{\mathbb R}
\def\P{\mathbb P}
\def\E{\mathbb E}

\def\N{\mathbb N}

\def\WW{\mathbb W}

\def\sE{\mathscr E}
\def\sF{\mathscr F}
\def\sD{\mathscr D}
\def\sC{\mathscr C}
\def\sB{\mathscr B}

\def\sL{\mathscr L}
\def\sG{\mathscr G}
\def\sP{\mathscr P}

\def\cH{\mathcal H}

\def\e{\operatorname{e}}

\def\d{\mathrm{d}}

\def\ff{\frac}
\def\ra{\rightarrow}
\def\nn{\nabla}
\def\pp{\partial}
\def\<{\langle}
\def\>{\rangle}
\def\sq{\sqrt}
\def\tld{\tilde}
\def\we{\wedge}
\def\1{\mathds{1}}

\def\supp{\displaystyle{\mathrm{supp}}}

\allowdisplaybreaks

\title{{\bf Local Convergence near Equilibria for Distribution-Dependent SDEs in a Generalized Kantorovich--Rubinstein Metric}
}

\author{
{\bf Shao-Qin Zhang }\\
\footnotesize{School of Statistics and Mathematics, Central University of Finance and Economics, Beijing 100081, China}\\
\footnotesize{Email: zhangsq@cufe.edu.cn}\\
}

\begin{document}

\maketitle

\begin{abstract}
Owing to exhibiting phase transitions, we investigate  the local convergence near a stationary distribution for distribution dependent stochastic differential equations.  By linearizing the nonlinear Markov semigroup associated with the distribution dependent equation around the stationary distribution, the local exponential convergence of the solution is related to the exponential convergence of a semigroup of linear operators. The generation and regularity of the linearized semigroup are investigated, and the Poincar\'e inequality for the stationary distribution is adapted to derive the exponential convergence of the linearized semigroup. Our results can be used as a criterion for the locally exponential stability of stationary distributions. Concrete examples, including the granular media equation with double-wells landscapes and quadratic interaction, are given to illustrate our main results.

\end{abstract}\noindent

AMS Subject Classification (2020): primary 60H10; secondary 35B35

\noindent

Keywords:  distribution dependent stochastic differential equations;  stationary distributions; phase transition; long-time behavior; stability

\vskip 2cm

\section{Introduction}

Stochastic differential equations (SDEs) whose coefficients depend on the own law of the solution were introduced by McKean in \cite{McK}. Such equations are named distribution dependent SDEs (DDSDEs), McKean-Vlasov SDEs, or mean-field SDEs in the literature, see e.g. \cite{BSY,BLPR,RZ,Wan18,Wan23}. Markov semigroups associated with DDSDEs are nonlinear.  \cite{Daw} established for the first time the existence of several stationary distributions, which is referred as phase transition in the literature,  for a DDSDE with a double-well confinement and a Curie-Weiss interaction on the line. Beside \cite{Daw},  phase transitions for DDSDEs were studied by many works, e.g. an equation with infinite many stationary distributions was given in \cite{AD};  local bifurcations were investigated in \cite{Tam,ZSQ23};  in a comprehensive series of papers of Tugaut et  al.,  equations with multi-wells confinement have been studied systemically, see  e.g. \cite{AA,DuTu,HT10a,Tug10,Tug14a};  for general DDSDEs, non-uniqueness of stationary distributions was discussed in \cite{ZSQ}; for  phase transitions  of McKean-Vlasov diffusions on the torus or nonlinear Markov jump processes, one can consult \cite{CGPS,FenZ}.   

For DDSDEs, one way to study the long-time behavior of  solutions is developing techniques in the ergodicity theory for distribution free SDEs, see e.g. \cite{BSWX,EGZ,LWZ,Wan23}.  However,  systems investigated there are without phase transitions.  When phase transitions exist, one stationary distribution of the DDSDE may not attract all solutions, unlike the ergodicity for distribution free SDEs.  Another way is using the gradient flow method or the free energy method, see e.g.  \cite{Bas,CGPS,MR,Tam,Tug10}. By using the free energy, local convergence of the density of solutions to some McKean-Vlasov equation under  $L^2$ norm was established in \cite{Tam}. It was proved in \cite{Tug10} that for the three distinct stationary distributions of  the classical granular media equation with double-well landscape and quadratic interaction potential, the symmetric one is unstable and the other two are stable. This method is powerful for concrete equations with phase transitions, while it relies on that coefficients of SDEs are of gradient forms.   In this paper, we introduce the third way. For a linear Markov semigroup,  it is well known that the spectrum of the generator is closely related to the ergodicity of the associated Markov process, see e.g. \cite{Chen,WBook,Wu}. Moreover,  the behavior of a dynamical system around an equilibrium can be characterized by the linearized system at the equilibrium point, see e.g. \cite{St,Tes}.  Following this idea, we want to linearize the nonlinear Markov process associated with the DDSDE around a stationary distribution in a suitable way, and find out the linearized semigroup  (or the associated infinitesimal generator) to characterize the local exponential convergence near the stationary distribution. Roughly speaking, we prove that if the linearized semigroup converges exponentially, then the stationary distribution is exponential stable, i.e. solutions of the DDSDE which start close to the stationary distribution converge exponentially to it. The convergence rate is also related to that of the linearized semigroup.

Consider the following equation
\begin{equation}\label{main-equ} 
\d X_t=b(X_t,\sL_{X_t})\d t+\si(X_t)\d B_t.
\end{equation}
where  coefficients $b:\R^d\times\sP\ra \R^d$ and $\si:\R^d\ra\R^d\otimes\R^d$ are measurable, $\sP$ is the space of probability measures on $\R^d$, and $B_t$ is a $d$-dimensional Brownian motion on a complete filtration probability space $(\Om,\sF,\{\sF_t\}_{t\ge 0},\mathbb{P})$, $\sL_{X_t}$ is the law of $X_t$ in the probability space $(\Om,\sF,\P)$.  We briefly explain the idea of linearizing the nonlinear Markov semigroup associated with $X_t$ for readers’ convenience. Let $\mu_0\in\sP$. For a solution to \eqref{main-equ} with $X_0\overset{d}{=}\mu_0$, we denote  $\mu_t=\sL_{X_t}$ and denote by $X^\mu_t$ the solution of the following equation
\beg{equation}\label{freezing}
\d X_t^{\mu}=b(X_t^{\mu},\mu_t)\d t+\si(X_t^{\mu})\d B_t.
\end{equation}
Let $P^\mu_t$ be the Markov semigroup associated with $X_t^{\mu}$, and $L_{\mu_t}$ be the infinitesimal generator:
\begin{align*}
L_{\mu_t}f(x)=\ff 1 2{\rm Tr}(\si \si^* \nn^2f)(x)+b(x,\mu_t)\cdot \nn f(x), f\in C^2,
\end{align*}
where $C^2$ is the set of  twice continuously differentiable functions on $\R^d$ and ${\rm Tr}$ is the trace of a matrix. Then $\mu_t$ is a solution of the following nonlinear Fokker-Planck equation on $\sP$:
\[\pp_t\mu_t=L_{\mu_t}^*\mu_t,\]
in the sense that 
\[\mu_t(f):=\int_{\R^d}f\d \mu_t=\mu_0(f)+\int_0^t\mu_s(L_{\mu_s}f)\d s,~t\geq 0,f\in C^2_0,\]
where $C^2_0$ is the set of twice continuously differentiable functions with compact support on $\R^d$.  
If there is a stationary distribution for \eqref{main-equ}, saying $\mu_\infty$, then we have a SDE as follows
\beg{equation}\label{equ-mu-inf}
\d X_t^{\mu_\infty}=b(X_t^{\mu_\infty},\mu_\infty)\d t+\si(X_t^{\mu_\infty})\d B_t.
\end{equation}
Let $P_t^{\mu_\infty}$ be the Markov semigroup associated with $X_t^{\mu_\infty}$, and $L_{\mu_\infty}$ be the infinitesimal generator:
\[L_{\mu_\infty}f(x)=\ff 1 2{\rm Tr}(\si \si^* \nn^2f)(x)+b(x,\mu_\infty)\cdot \nn f(x),~f\in C^2.\]
Then $\mu_{\infty}$ satisfies the stationary Fokker-Planck equation $L_{\mu_\infty}^*\mu_\infty=0$. By using equations of $\mu_t$ and $\mu_\infty$, we have that 
\beg{align*}
\pp_t(\mu_t-\mu_\infty)&=L_{\mu_t}^*\mu_t-L_{\mu_\infty}^*\mu_{\infty}\\
&=(L_{\mu_t}-L_{\mu_\infty})^*(\mu_t-\mu_\infty)+(L_{\mu_t}-L_{\mu_\infty})^*\mu_\infty+L_{\mu_\infty}^*(\mu_t-\mu_\infty).
\end{align*}
For simplicity, we take  that $b(x,\mu)=\mu(B(x,\cdot))$ for some function $B(\cdot,\cdot)$ on $\R^d\times\R^d$ as an example. Then, formally, for $f\in C^2_0$
\beg{align*}
\pp_t\left(\mu_t(f)-\mu_{\infty}(f)\right)&=(\mu_t-\mu_\infty)(L_{\mu_\infty}f)+\mu_\infty\left((L_{\mu_t}-L_{\mu_\infty})f\right)+(\mu_t-\mu_\infty)\left((L_{\mu_t}-L_{\mu_\infty})f\right)\\
&=\int_{\R^d}L_{\mu_\infty}f(z)(\mu_t-\mu_\infty)(\d z)+\int_{\R^d}\mu_\infty\left(B(\cdot,z)\cdot \nn f(\cdot)\right)(\mu_t-\mu_\infty)(\d z) \\
&\quad\,+(\mu_t-\mu_\infty)\left((b(\cdot,\mu_t)-b(\cdot,\mu_\infty))\cdot \nn f(\cdot)\right)\\
&=\int_{\R^d}\left(L_{\mu_\infty}f(z)+\mu_\infty\left(B(\cdot,z)\cdot \nn f(\cdot)\right)\right)(\mu_t-\mu_\infty)(\d z)\\
&\quad\,+(\mu_t-\mu_\infty)\left((b(\cdot,\mu_t)-b(\cdot,\mu_\infty))\cdot \nn f(\cdot)\right),~f\in C^2_0.
\end{align*}
The first term of the last equality above is a linear term of  $\mu_t-\mu_\infty$, and the second term  is a ``second order infinitesimal of $\mu_t-\mu_\infty$". Let 
\[\bar A f(z)=\int_{\R^d}\left(B(x,z)\cdot \nn  f(x)\right)\mu_\infty(\d x).\]
Then one can see that the linearization of $\pp_t\mu_t=L_{\mu_t}^*\mu_t$ at $\mu_\infty$ is  the linear equation $\pp_t\mu_t=(L_{\mu_\infty}+\bar A)^*\mu_t$.
Note that $B(x,\cdot)$ is an external derivative or linear functional derivative of $b(x,\cdot)$ indeed. When $b(x,\mu)$ does not depends linearly on  $\mu$, we can use the linear functional derivative to linearize $b(x,\cdot)$. For derivatives of functionals on $\sP$, one can consult \cite{CarDe,RW} for detailed discussions.  

We study the semigroup generated by $L_{\mu_\infty}+\bar A$ in Section 4. Then we prove that the exponential convergence of the semigroup generated by $L_{\mu_\infty}+\bar A$ leads to the local convergence of solutions for \eqref{main-equ}. That is our main result of this paper, see  Theorem \ref{thm0} for precise statement in Section 2 and its proof in Section 5.  The Poincar\'e inequality for $L_{\mu_\infty}$ can be used to derive the exponential convergence of the semigroup generated by $L_{\mu_\infty}+\bar A$, see Proposition \ref{prop-22} and Proposition \ref{lem-2.3} below.   Moreover, the convergence rate is also given by using the convergence rates of $P_t^{\mu_\infty}$ and the semigroup generated by $L_{\mu_\infty}+\bar A$.

Recently, the long-time behavior of DDSDEs is also investigated by  linearizing  nonlinear Markov processes in \cite{Co}. Comparing with \cite{Co}, where the Lions derivative is used, we adopt the linear functional derivative,  and a clearer formulation of the linearization semigroup is provided in this paper. Moreover, their criteria cannot deal with models in \cite{Daw,Tug14a}.  In Section 6, we apply our result (Theorem \ref{thm0}) to concrete models studied in \cite{Co} and  \cite{Daw}, and establish the local convergence for them. The local convergence rate for McKean-Vlasov equations is also investigated in \cite{MR}, where their argument is based on the Wasserstein gradient flow and nonlinear functional inequality, which is quite different from us.  Moreover, log-Sobolev inequality is used in \cite{MR}, which excludes models whose stationary distributions do not satisfy the Gaussian concentration. 

This paper is structured as follows. In Section 2, we present our main results and concrete examples.  In Section 3, we study the regularity of decoupled SDEs associated with \eqref{main-equ}. In Section 4, we establish a generation theorem  for $C_0$-semigroups of $L_{\mu_\infty}+\bar A$ and investigate also the regularity of the semigroup.  The proof of our main result is given in Section 5. Concrete examples are discussed in Section 6.     

{\bf Notation:} The following notations are used in the sequel.\\
$\bullet$ We denote by $L^p$ (resp. $L^p(\mu)$) the space of functions for which the $p$-th power of the absolute value is Lebesgue integrable (resp. integrable w.r.t. the measure $\mu$); $C_0$ (resp.  $C^1_0$) the space of all the continuous (resp. continuously differentiable) functions with compact support on $\R^d$; $C^k$ the $k$ continuously differentiable functions on $\R^d$; $C^k_b$ the bounded continuously differentiable functions up to  $k$-th order on $\R^d$. For $\R^{d_1}$-valued functions, similar notations ($C^k(\R^d,\R^{d_1}),C_b^k(\R^d,\R^{d_1})$, etc) are used; for $\R^{d_1}\otimes\R^{d_2}$-valued functions,  similar notations ($C^k(\R^d,\R^{d_1}\otimes \R^{d_2}),C_b^k(\R^d,\R^{d_1}\otimes \R^{d_2})$, etc) are used.\\
$\bullet$ For a probability measure $\mu$, let 
\[W^{1,2}_{\mu}=\{f\in W^{1,2}_{loc}~|~f,\nn f\in L^2(\mu)\},~~\|f\|_{W^{1,2}_{\mu}}=\|f\|_{L^2(\mu)}+\|\nn f\|_{L^2(\mu)}.\]
We denote by $\sD(L)$ the domain of the linear operator $L$ and by $\Si(L)$ the spectrum of $L$.\\
$\bullet$ To emphasise the initial distribution, we use $X^{\mu_\infty}_t(\nu)$ and $X^{\mu}_t(\nu)$, i.e. $X^{\mu_\infty}_0(\nu)\overset{d}{=}\nu$ and $X^{\mu}_0(\nu)\overset{d}{=}\nu$, respectively. When $\nu=\de_x$ with $x\in\R^d$, we use $X^{\mu_\infty}_t(x)$ and $X^{\mu}_t(x)$.\\
$\bullet$ For any measurable function $V\geq 1$ on $\R^d$, let
\[\sP_{V}:=\{\mu\in\sP~|~\|\mu\|_V:=\mu(V)<+\infty\},\]
Define the weight total variance distance on $\sP_V$:
\[\|\mu-\nu\|_{V}=\sup_{|f|\leq V} \left|\mu(f)-\nu(f)\right|,~\mu,\nu\in\sP_V.\]

\section{Main Result}
We first introduce the distance used to  investigate the local convergence near a stationary distribution of \eqref{main-equ}. Let
\[\sG=\left\{\ph:[0,+\infty)\mapsto [0,+\infty)~\text{increasing,}~\ph(0+)=0,~\ph^2(\sq{\cdot})~\text{is concave,}~\int_0^1\ff {\ph(s)} {s}\d s<+\infty\right\},\]
and for $\ph\in\sG$ and a measurable function $V\geq 1$, let 
\beg{align*}
\mathscr{G}_{V,\ph}&=\left\{g\in\sB(\R^d)\left|~\|g\|_{V,\ph}:=\sup_{x\neq y}\ff {|g(x)-g(y)|} {\ph(|x-y|)(V(x)+V(y))}<+\infty \right.\right\},\\
\sP_{V,\ph}&=\left\{\mu\in\sP~|~\mu\left((1+\ph(|\cdot|))V(\cdot)\right)<+\infty\right\}.
\end{align*}
For $\mu,\nu\in\sP_{V,\ph}$, the generalized Kantorovich--Rubinstein distance between $\mu$ and $\nu$ is defined by
\beg{align*}
\|\mu-\nu\|_{V,\ph}=\sup_{g\in\mathscr{G}_{V,\ph},~\|g\|_{V,\ph}\leq 1}\int_{\R^d}g(x)(\mu-\nu)(\d x).
\end{align*}
When $V(x)=(1+|x|^2)^{\ff p 2}$ for some $p\geq 0$,  $\sP_{V}$ ($\|\cdot\|_V,\sG_{V,\ph},\sP_{V,\ph},\|\cdot\|_{V,\ph}$) will be denote by $\sP_p$ (resp. $\|\cdot\|_p,\sG_{p,\ph},\sP_{p,\ph},\|\cdot\|_{p,\ph}$). 
\beg{rem}\label{ad-rem000}
We remark that  $\|\mu-\nu\|_{V,\ph}$ gives a distance between $\mu$ and $\nu$. Since $\ph^2(\sq{\cdot})$ is concave and $\ph(0)=0$, we have that
\[\ph(r)=\sq{\ph^2(\sq{r^2})}=\sq{\ph^2(\sq{r^2\cdot 1+(1-r^2)\cdot 0})}\geq \sq{r^2\ph^2(1)}=r\ph(1),~r\in [0,1].\]
Combining this with that $\ph$ is increasing, there is a constant $C>0$ such that 
\beg{equation}\label{Cr<ph}
C(r\we 1)\leq \ph(r),~r\geq 0.
\end{equation} 
Then, for any  bounded and Lipschitz  function $g$, there is 
\beg{align*}
\sup_{x\neq y}\ff {|g(x)-g(y)|} {\ph(|x-y|)(V(x)+V(y))}\leq \sup_{x\neq y}\ff {|x-y|\we (2\|g\|_\infty)} {C(|x-y|\we 1)}<+\infty.
\end{align*}
Hence, $\mathscr{G}_{V,\ph}$ is separating for $\sP_{V,\ph}$, and $\|\mu-\nu\|_{V,\ph}$ is a distance between $\mu$ and $\nu$ in $\sP_{V,\ph}$. 
\end{rem}
 
To linearize \eqref{main-equ}, we introduce the linear functional derivative of a function on $\sP_{V,\ph}$ for some $\ph\in \sG$ and $V\geq 1$.
\beg{defn}
A function $u:\sP_{V,\ph}\mapsto \R$ is called linear functional differentiable on $\sP_{V,\ph}$ if there is a measurable function
\[\R^d\times\sP_{V,\ph}\ni (x,\mu)\mapsto D^F_{\mu}u(x)\]
such that
\beg{equation}\label{DFuu}
\sup_{x\in\R^d}\ff {\sup_{\|\mu-\de_0\|_{V,\ph}\leq L} |D^F_{\mu} u(x)|} {(1+\ph(|x|))V(x)}<+\infty,~L\geq 0,
\end{equation}
and 
\[u(\mu)-u(\nu)=\int_0^1\d r\int_{\R^d}D^F_{r\mu+(1-r)\nu} u(x)(\mu-\nu)(\d x),~\mu,\nu\in\sP_{V,\ph}.\]
$D_{\mu}^Fu$ is called the linear functional derivative of $u$ at $\mu$ on $\sP_{V,\ph}$.
\end{defn}
\beg{rem}\label{re-DF}
The linear functional derivative is unique up to a constant. We always choose the linear functional derivative $D^F_{\mu}u$ such that  $\mu(D^F_{\mu}u)=0$. 
\end{rem}

We aim to investigate the asymptotic behaviour of the solution to \eqref{main-equ} near a stationary distribution, so we also make the following assumption
\beg{description}
\item [(H0)] \eqref{main-equ} has a stationary distribution $\mu_\infty$  such that for any $p>0$, $\mu_\infty(|\cdot|^p)<+\infty$.  

\end{description}
For the existence of stationary distributions for DDSDEs, one can consult \cite{ZSQ}. We assume that coefficients $b,\si$ satisfy following assumptions.
\beg{description}[align=left, noitemsep]
\item [(H1)] There is $1\leq U_0\in C^2$ with $\lim_{|x|\ra +\infty}U_0(x)=+\infty$, and 
\beg{align}
&\mu_\infty(U_0^p)<+\infty,~p\geq 1,\label{mu-inf-U0}\\
&\sup_{x\in\R^d,\mu\in\sP_{U_0}}\ff {\<b(x,\mu),\nn U_0(x)\>} {U_0(x)+\|\mu\|_{U_0}}<+\infty,\label{nnb1-Lypu0}\\
&\sup_{x\in\R^d}\ff {|\si^*\nn U_0(x)|+\|\si\si^*\nn^2 U_0(x)\|}  {U_0(x)}<+\infty,\label{nnb1-Lypu1}
\end{align}
such that for any $\mu\in \sP_{U_0}$, the drift $b(\cdot,\mu)\in C^2(\R^d,\R^d)$, and there exist a locally bounded  function $K:\sP_{U_0}\ra [0,+\infty)$ and a nonnegative constant $\be_1$  such that
\beg{align}\label{Inequ-nnb1}
\<\nn b(\cdot,\mu)(x) v,v\>&\leq K(\mu)|v|^2,\\
|\nn^2 b(\cdot,\mu)(x)|&\leq K(\mu)(1+|x|)^{\be_1},~x,v\in\R^d,~\mu\in\sP_{U_0}.\label{Inequ-nnb2}
\end{align}

\item [(H2)] The diffusion term $\si\in C^2_b(\R^d,\R^d\otimes\R^d)$, and  there is a positive constant $\si_0$ such that 
\beg{equation}\label{nondege0}
\si(x)\si^*(x)\geq \si_0^2,~x\in\R^d.
\end{equation}

\item [(H3)] There exist  $\ph_0\in\sG$, $1\leq V_0\in C^2$ satisfying \eqref{nnb1-Lypu1} with $U_0$ replaced by $V_0$, and  $K_V>0$ such that   
\beg{align}\label{ph-p-q0}
(1+\ph_0(|x|))V_0(x)&\leq K_VU_0(x),\\
\sup_{|v|\leq 1} V_0(x+v)&\leq K_V V_0(x),\label{vV0V0}\\
 \<b(x,\mu_\infty),\nn V_0(x)\> &\leq K_V V_0(x),\label{LypuV0}\\
|b(x,\mu)-b(x,\nu)|&\leq K_V\|\mu-\nu\|_{V_0,\ph_0},~x\in\R^d,~\mu,\nu\in\sP_{V_0,\ph_0}.\label{b1-mu-nu}
\end{align}

\end{description}

\beg{rem}
Since \eqref{ph-p-q0}, for any $f\in \sG_{V_0,\ph_0}$ with $\|f\|_{V_0,\ph_0}\leq 1$, there is 
\beg{align*}
|f(x)-f(0)|&\leq \ph_0(|x|)(V_0(x)+V_0(0))\\
&\leq (1+V_0(0))\ph_0(|x|)V_0(x)\leq (1+V_0(0)) K_V U_0(x).
\end{align*}
Thus
\beg{align}
\|\mu-\nu\|_{V_0,\ph_0}&=\sup_{f\in\sG_{V_0,\ph_0},\|f\|_{V_0,\ph_0}\leq 1}|\mu(f)-\nu(f)|\nonumber\\
&\leq \sup_{|f-f(0)|\leq (1+V_0(0)) K_V U_0}|\mu(f-f(0))-\nu(f-f(0))|\nonumber\\
&\leq (1+V_0(0)) K_V \sup_{|g|\leq U_0}|\mu(g)-\nu(g)|\nonumber\\
&=(1+V_0(0)) K_V\|\mu-\nu\|_{U_0}.\label{norm-p0q0}
\end{align}
This implies that there is $C>0$ such that
\beg{equation}\label{ad-0q0}
|b(x,\mu)-b(x,\nu)|\leq C\|\mu-\nu\|_{U_0},~x\in\R^d,\mu,\nu\in\sP_{U_0}.
\end{equation} 
By using \eqref{nnb1-Lypu0}-\eqref{Inequ-nnb1}, {\bf (H2)}, and \eqref{ad-0q0}, we can derive as \cite[Theorem 1.1]{Ren} that,  for any $T>0$, $p\geq 1$ and $X_0$ with $\sL_{X_0}\in \sP_{U_0^p}$, \eqref{main-equ} has a unique solution $X_t$ with $\sL_{X_\cdot}\in C([0,T],\sP_{U_0^p})$. 

\end{rem}

\beg{rem}
It is clear that  \eqref{mu-inf-U0}, \eqref{Inequ-nnb1} and {\bf (H2)} imply the wellposedness of \eqref{equ-mu-inf}. Moreover, it follows from \cite[(1) of Theorem 1.1]{Wang-11} that the log-Harnack inequality holds for $P_t^{\mu_\infty}$. Thus $\mu_\infty$ is the only invariant probability measure of \eqref{equ-mu-inf}.

\end{rem}

We also assume that $b$ satisfies following conditions.
\beg{description}
\item [(H4)] There exists  $C>0$  such that for any $x_1,x_2\in\R^d$ and $\mu,\nu\in\sG_{V_0,\ph_0}$, 
\beg{align}
|b(x_1,\mu)-b(x_1,\nu)-(b(x_2,\mu)-b(x_2,\nu))|&\leq C\ph_0(|x_1-x_2|)\|\mu-\nu\|_{V_0,\ph_0}.\label{hhh0}
\end{align}
For any $x\in\R^d$, $b(x,\cdot)$ is linear functional differentiable on $\sP_{V_0,\ph_0}$. There is $F\in L^2(\mu_\infty)$ such that
\beg{align}
\sup_{z\in\R^d}\ff {|D^F_{\mu_\infty}b(x,\cdot)(z)|} {(1+\ph_0(|z|))V_0(z)}+\|D^F_{\mu_\infty}b(x,\cdot)\|_{V_0,\ph_0}\leq F(x),~x\in\R^d,\label{DF-DF}
\end{align}
and there is $C>0$ so that
\beg{equation}\label{DFb-pph}
\left\| D^F_{\mu}b(x,\cdot)-D^F_{\nu}b(x,\cdot)\right\|_{V_0,\ph_0}\leq C\|\mu-\nu\|_{V_0,\ph_0},~\mu,\nu\in\sP_{V_0,\ph_0},~x\in\R^d.
\end{equation}

\end{description}
In the following discussion, we denote 
\[D^F_{\mu}b(x, z)=D^F_{\mu}b(x, \cdot)(z)\] 
for simplicity. Due to {\bf (H4)} and $\mu_\infty\in\sP_{V_0,\ph_0}$, we let 
\beg{equation}\label{def-barA}
\bar Af(z)=\int_{\R^d} D^F_{\mu_\infty}b(x, z)\cdot\nn f(x) \mu_\infty(\d x),~f\in W^{1,2}_{\mu_\infty}.
\end{equation}
\beg{rem}\label{Rem2.5}
Combining this with Remark \ref{re-DF}, \eqref{DF-DF}, $(1+\ph_0(|\cdot|)V_0(\cdot)\in L^2(\mu_\infty)$ (due to \eqref{ph-p-q0} and \eqref{mu-inf-U0}) and the Fubini theorem, we have that $\mu_\infty(\bar Af)=0$. 
\end{rem}

Our main result indicates that the exponential convergence of the semigroup generated by $L_{\mu_\infty}+\bar A$ can lead to the exponential convergence of $\sL_{X_t}$ to $\mu_\infty$ in the metric $\|\cdot\|_{V_0,\ph_0}$. 
\beg{thm}\label{thm0}
(1) Assume $\mu_\infty$ is an invariant  probability measure of \eqref{equ-mu-inf}, $(b(\cdot,\mu_\infty),\si(\cdot))$ satisfies \eqref{Inequ-nnb1}, \eqref{Inequ-nnb2}, {\bf (H2)},  and the integral kernel of $\bar A$ satisfies $D_{\mu_\infty}^F b\in L^2(\mu_\infty\times\mu_\infty)$. Then $L_{\mu_\infty}+\bar A$ generates a $C_0$-semigroup, saying $Q_t$, on $L^2(\mu_\infty)$. In particular, if {\bf (H0)}-{\bf (H3)} hold except \eqref{vV0V0} and  \eqref{b1-mu-nu},  and \eqref{DF-DF} holds, then $L_{\mu_\infty}+\bar A$ generates a $C_0$-semigroup, and 
\beg{equation}\label{invQtmu}
\mu_\infty(Q_tf)=\mu_\infty(f),~f\in L^2(\mu_\infty).
\end{equation}
(2) Suppose  that {\bf (H0)}-{\bf (H4)} hold, there are $C_W\geq 1$ and $\la_P>0$ such that
\beg{equation}\label{con-W1}
\|P_t^{\mu_\infty}f\|_{V_0,\ph_0}\leq C_We^{-\la_Pt}\|f\|_{V_0,\ph_0},~f\in \mathscr{G}_{V_0,\ph},
\end{equation}
and there are $C_Q\geq 1$ and $\la_Q>0$ such that
\beg{equation}\label{Q-exp-dec}
\|Q_tf\|_{L^2(\mu_\infty)}\leq C_Q e^{-\la_Q t}\|f\|_{L^2(\mu_\infty)},~f\in L^{2}(\mu_\infty),~\mu_\infty(f)=0.
\end{equation}
Then there exists $C\geq 1$ depending on $C_Q,\la_Q$, $C_W,\la_P$, $V_0$, $\|F\|_{L^2(\mu_\infty)}$, $\mu_\infty$  such that for any $\ep\in (0,1)$ and any $\mu_0\in\sP_{U_0}$ with 
\beg{equation}\label{add-con-init}
\|\mu_0-\mu_\infty\|_{V_0,\ph_0}<\ff {[(\ep(\la_P\we\la_Q))^2\we 1]} {16C^2}\left(\left(\int_0^1\ff{\ph_0(s)} s\d s\right)^{-1}\we\left(\ff { 2(1-\ep)(\la_p\we\la_Q)} {\ph_0(1)}\right)\right),
\end{equation}
the law of the solution for \eqref{main-equ} with  $X_0\overset{d}{=}\mu_0$ converges to $\mu_\infty$ exponentially: 
\beg{equation}\label{th-loc-cov}
\|\sL_{X_{t}}-\mu_\infty\|_{V_0,\ph_0}<  \ff {2C e^{-(1-\ep)(\la_P\we\la_Q) t}} {(\ep(\la_P\we\la_Q))\we 1}\|\sL_{X_0}-\mu_\infty\|_{V_0,\ph_0},~t\geq 0.
\end{equation}
If $\la_P\neq \la_Q$, then there is $C\geq 1$ depending on $C_Q,\la_Q$, $C_W,\la_P$, $V_0$, $\|F\|_{L^2(\mu_\infty)}$, $\mu_\infty$  such that
\beg{equation*}
\|\sL_{X_{t}}-\mu_\infty\|_{V_0,\ph_0}<  \ff {2C e^{-(\la_P\we\la_Q) t}} {|\la_P-\la_Q|\we 1}\|\sL_{X_0}-\mu_\infty\|_{V_0,\ph_0},~t\geq 0,
\end{equation*}
when 
\beg{equation*}
\|\mu_0-\mu_\infty\|_{V_0,\ph_0}<\ff {(\la_P-\la_Q)^2} {16C^2}\left(\left(\int_0^1\ff{\ph_0(s)} s\d s\right)^{-1}\we\left(\ff { 2(\la_p\we\la_Q)} {\ph_0(1)}\right)\right).
\end{equation*}
\end{thm}
\beg{rem}
This theorem is a criteria for the exponential stability of $\mu_\infty$ : The stationary distribution $\mu_\infty$ of \eqref{main-equ} is exponentially stable under the metric $\|\cdot\|_{V,\ph}$ in $\sP_{V,\ph}$ for some $V\geq 1, \ph\in\sG$, if there are positive constants $\de,C,\la$ such that
\[\|\sL_{X_t}-\mu_\infty\|_{V,\ph}\leq Ce^{-\la t}\|\sL_{X_0}-\mu_\infty\|_{V,\ph},~\|\sL_{X_0}-\mu_\infty\|_{V,\ph}<\de.\]
For the instability of  stationary distributions, one can see our another paper \cite{ZSQ25-ins}. 
\end{rem}

\beg{rem}
According to \eqref{add-con-init} and  \eqref{th-loc-cov}, we obtain a convergence rate $\la:=(1-\ep)\la_P\we\la_Q<\la_P\we\la_Q$.  The closer to $\mu_\infty$ the initial value $\mu_0$ is, the faster the convergence rate  for $\sL_{X_t}$ approaches $\la_P\we\la_Q$. 

In the case that $\ph_0(r)=r,~r\leq 1$, there is 
\[\left(\left(\int_0^1\ff{\ph_0(s)} s\d s\right)^{-1}\we\left(\ff {1\we (2\la)} {\ph_0(1)}\right)\right)=1\we (2\la).\]
Then
\[\max_{\la\in (0,\la_P\we\la_Q)}\left[((\la_P\we\la_Q-\la)^2\we 1)\left(\ff 1 2\we \la\right)\right]=\ff {4(\la_P\we\la_Q)^3} {27}\we \ff 1 2,\]
and the maximal value is achieved at $\la=\ff {\la_P\we\la_Q} 3\we\ff 1 2$. Then, for $\mu_0$ with
\[\|\mu_0-\mu_\infty\|_{V_0,\ph_0}<\ff {1} {16C^2}\left(\ff {4(\la_P\we\la_Q)^3} {27}\we \ff 1 2\right),\]
the convergence rate $\la\geq \ff {\la_P\we\la_Q} 3\we\ff 1 2$.
\end{rem}

Due to  the (weak) spectral mapping theorem, the exponential decay \eqref{Q-exp-dec} can be established by analysing the spectrum of the generator $L_{\mu_\infty}+\bar A$  for $Q_t$ possessing some  regularity properties, see \cite[Corollary V.2.10 and Corollary V.2.11]{EN-short}. Let $L^2_{\mathbb{C}}(\mu_\infty)$ be the complexification of $L^2(\mu_\infty)$ with the inner product defined as follows  
\[\<f,g\>_{L^2_{\mathbb{C}}}=\mu_\infty(f\bar g),~f,g\in L^2_{\mathbb{C}}(\mu_\infty),\]
where $\bar g$ is the complex conjugate of $g$. A real operator $L$ in $L^2(\mu_\infty)$ is extended to $L^2_{\mathbb{C}}(\mu_\infty)$ in the following way
\[L(f_1+\sq{-1}f_2)=Lf_1+\sq{-1}Lf_2,~f_1,f_2\in\sD(L).\]
We denote by $\Sigma(L)$ the spectrum of the  complexification of $L$ in $L^2_{\mathbb{C}}(\mu_\infty)$. To prove \eqref{Q-exp-dec}, we give following two propositions.  The first one is an abstract result, and the second one concerns the case where $L_{\mu_\infty}$ is self-adjoint and $D^F_{\mu_\infty}b$ is of the gradient form.  The Poincar\'e inequality is used in both propositions. 
\beg{prp}\label{prop-22}
Assume that $P^{\mu_\infty}_t$ is a $C_0$-contraction semigroup on $L^2(\mu_\infty)$, there are positive constants $C,\tld\la$ such that 
\beg{equation}\label{exp-L2-mu}
\|P^{\mu_\infty}_tf\|_{L^2(\mu_\infty)}\leq Ce^{-\tld\la  t}\|f\|_{L^2(\mu_\infty)},~\mu_\infty(f)=0,~t\geq 0,
\end{equation}
and $\bar A$ can be extended to a compact operator on $L^2(\mu_\infty)$. If $I+\bar A(L_{\mu_\infty}-z)^{-1}$ has a bounded  inverse on $\cH:=\{f\in L^2_{\mathbb{C}}(\mu_\infty)~|~\mu_\infty(f)=0\}$  for any $z$ with ${\rm Re}z\geq 0$, where $L_{\mu_\infty}^{-1}$ is the inverse on $\cH$, then \eqref{Q-exp-dec} holds.   
\end{prp}
\beg{proof}
We first prove that 
\beg{equation}\label{Si_L+AC+}
\Sigma\left(L_{\mu_\infty}+\bar A\right)\cap \{z|~{\rm Re}z\geq 0,~z\neq 0\}=\emptyset.
\end{equation}
It follows from \eqref{exp-L2-mu} that $\Sigma(L_{\mu_\infty})\subset \{z|~{\rm Re}z\leq -\tld\la\}\cup\{0\}$.
For any $z\in \{z|~{\rm Re}z\geq 0,~z\neq 0\}$,  it holds that $(L_{\mu_\infty}-z)^{-1}$ is a bounded operator on $L^2_{\mathbb{C}}(\mu_\infty)$ since $z\not\in \Sigma(L_{\mu_\infty})$. For any $g\in L^2_{\mathbb{C}}(\mu_\infty)$, $\bar A  (L_{\mu_\infty}-z)^{-1}g\in \cH$ due to that $\bar A$ is  bounded and Remark \ref{Rem2.5}. We set 
\beg{align*}
h&=(I+\bar A(L_{\mu_\infty}-z)^{-1})^{-1}\bar A (L_{\mu_\infty}-z)^{-1}g,\\
f&=(L_{\mu_\infty}-z)^{-1}(-h+g).
\end{align*}
Then 
\beg{align*}
\bar Af 
&=-\bar A(L_{\mu_\infty}-z)^{-1}(I+\bar A(L_{\infty}-z)^{-1})^{-1}\bar A (L_{\mu_\infty}-z)^{-1} g+ \bar A(L_{\mu_\infty}-z)^{-1}g\\
& = -\bar A (L_{\mu_\infty}-z)^{-1} g +(I+\bar A(L_{\infty}-z)^{-1})^{-1}\bar A (L_{\mu_\infty}-z)^{-1} g\\
&\quad\, +\bar A(L_{\mu_\infty}-z)^{-1}g\\
&= h.
\end{align*}
Hence 
\beg{align*}
(L_{\mu_\infty}+\bar A-z)f=(L_{\mu_\infty}-z)f+\bar Af=-h+g+h=g.
\end{align*}
Moreover, if the equation $(L_{\mu_\infty}+\bar A-z)f=0$ has a solution, then
\beg{align}\label{IA-1A}
(I+\bar A(L_{\mu_\infty}-z)^{-1})\bar Af&=\bar A \left( f+(L_{\mu_\infty}-z)^{-1}\bar A f\right)\nonumber\\
&=\bar A (L_{\mu_\infty}-z)^{-1}\left((L_{\mu_\infty}-z)f+\bar A f\right)\nonumber\\
&=0.
\end{align}
This, together with $\mu_\infty(\bar Af)=0$ and $I+\bar A(L_{\mu_\infty}-z)^{-1}$ has a bounded  inverse on $\cH$, implies that $\bar Af=0$. Consequently, $(L_{\mu_\infty}-z)f=0$, and $f=0$. Hence, $(L_{\mu_\infty}+\bar A-z)^{-1}$ is well-defined. 
Moreover,  
\beg{align*}
&\|(L_{\mu_\infty}+\bar A-z)^{-1}g\|_{L^2(\mu_\infty)}=\|f\|_{L^2(\mu_\infty)}\\
&\leq \|(L_{\mu_\infty}-z)^{-1}h\|_{L^2(\mu_\infty)}+\|(L_{\mu_\infty}-z)^{-1}g\|_{L^2(\mu_\infty)}\\
& \leq \|(L_{\mu_\infty}-z)^{-1}\|_{L^2(\mu_\infty)}\|(I+\bar A(L_{\infty}-z)^{-1})^{-1}\|_{L^2(\mu_\infty)}\|\bar A(L_{\mu_\infty}-z)^{-1}\|_{L^2(\mu_\infty)}\|g\|_{L^2(\mu_\infty)}\\
&\quad\, +\|(L_{\mu_\infty}-z)^{-1}\|_{L^2(\mu_\infty)}\|g\|_{L^2(\mu_\infty)}.
\end{align*}
Therefore, $(L_{\mu_\infty}+\bar A-z)^{-1}$ is bounded operator on $L^2(\mu_\infty)$.

Next, we prove that $0$ is a simple eigenvalue of $L_{\mu_\infty}+\bar A$, i.e. the algebraic multiplicity equals to $1$. It is clear that $\1$ is an eigenfunction of $L_{\mu_\infty}+\bar A$.  On $\cH$,  $L_{\mu_\infty}^{-1}$ exists and is a bounded operator  since \eqref{exp-L2-mu}.   If there is a $f\in\sD(L_{\mu_\infty})\cap\cH$ with $f\neq 0$ and $(L_{\mu_\infty}+\bar A)f=0$, then,  as \eqref{IA-1A}, there is 
\[(I+\bar AL_{\mu_\infty}^{-1})\bar A f=\bar A(f+L_{\mu_\infty}^{-1}\bar Af)=\bar A L_{\mu_\infty}^{-1}(L_{\mu_\infty}f+\bar Af)=0.\]
Thus, $\bar A f=0$, and $L_{\mu_\infty}f=0$. Consequently, $f$ is a constant, which leads to a contradiction. Hence, the geometric multiplicity of the eigenvalue $0$ for for $L_{\mu_\infty}+\bar A$ is $1$.  If $0$ is not a simple eigenvalue, then there is a nontrivial $f\in \sD(L_{\mu_\infty})\cap\cH$ such that $(L_{\mu_\infty}+\bar A)f=\1$. This is a contradiction since $\mu_\infty((L_{\mu_\infty}+\bar A)f)=0$. 

Finally, we prove \eqref{Q-exp-dec}.  It follows from \eqref{exp-L2-mu} and \cite[Proposition V.4.5]{EN-short} that $P^{\mu_\infty}_t$ is a quasi-compact strongly continuous semigroup. Combining this with that $\bar A$ is a compact operator and \cite[Proposition V.4.9]{EN-short}, we have that $Q_t$ is also a quasi-compact strongly continuous semigroup. Therefore, \eqref{Q-exp-dec} follows from \cite[Theorem V.4.6]{EN-short}, \eqref{Si_L+AC+} and that the algebraic multiplicity of the eigenvalue $0$ for $L_{\mu_\infty}+\bar A$ equals $1$.

\end{proof}

\beg{rem}\label{rem2.8}
It is well-known that \eqref{exp-L2-mu} is equivalent to the following Poincar\'e inequality: 
\beg{align}\label{Poincare}
\mu_\infty(f^2)\leq \ff 1 {\tld\la}\sE(f,f)+\mu_\infty(f)^2,~f\in\sD(\sE),
\end{align}
where $\sE(f,g):=-\mu_\infty(fL_{\mu_\infty}g)$ for $f,g\in\sD(L_{\mu_\infty})$, and $\sD(\sE)$ is the closure of $\sD(L_{\mu_\infty})$ under the norm $\|f\|_{\sE}:=\|f\|_{L^2(\mu_\infty)}+\sq{\sE(f,f)}$.
\end{rem}
To verify the condition that $I+\bar A(L_{\mu_\infty}-z)^{-1}$ has a bounded  inverse on $\cH$  for any $z$ with ${\rm Re}z\geq 0$ in Proposition \ref{prop-22}, we give the following proposition.
\beg{prp}\label{lem-2.3}
Assume that $L_{\mu_\infty}$ is self-adjoint satisfying \eqref{Poincare} and
\[\sE(f,g)=\ff 1 2\mu_\infty(\<\si^*\nn f,\si^*\nn g\>),~f,g\in W_{\mu_\infty}^{1,2}.\]
Suppose that there is $W\in L^2(\R^d\times\R^d,\mu_\infty\times\mu_\infty)$ satisfying 
\beg{align}
&W(\cdot,y)\in W_{\mu_\infty}^{1,2}\cap \sD(L_{\mu_\infty}),\quad \mu_\infty(W(x,\cdot))=0,~x,y\in\R^d,\label{WLW=0}\\
\int_{\R^d\times\R^d}&\left(|\nn W(\cdot,y)(x)|^2+ |L_{\mu_\infty}W(\cdot,y)(x)|^2\right)\mu_\infty(\d x)\mu_\infty(\d y)<+\infty,\label{LW-L2}\\
&\qquad\int_{\R^d\times\R^d} |W(x,y)|^2\mu_\infty(\d x)\mu_\infty(\d y)<1,\label{WL2}
\end{align}
such that
\beg{equation}\label{DF=si2W}
D^F_{\mu_\infty}b(x,y)=\ff 1 2\si(x)\si(x)^*\nn W(\cdot,y)(x).
\end{equation}
Then $\bar A$ can be extended to a compact operator on $L^2(\mu_\infty)$, and $I+\bar A(L_{\mu_\infty}-z)^{-1}$ has a bounded  inverse on $\cH$  for any $z$ with ${\rm Re}z\geq 0$. As a consequence, \eqref{Q-exp-dec} holds. 
\end{prp}
\beg{proof}
Because of \eqref{WLW=0} and \eqref{LW-L2}, we have that $\mu_\infty(\nn W(x,\cdot))=0$. Thus, $D^F_{\mu_\infty}b(x,y)$ satisfies $\mu_\infty(D^F_{\mu_\infty}b(x,\cdot))=0$, which meets the requirement in Remark \ref{re-DF}. For any $f\in W^{1,2}_{\mu_\infty}$, since $W(\cdot,y)\in W_{\mu_\infty}^{1,2}\cap \sD(L_{\mu_\infty})$, there is 
\beg{align*}
\bar Af(y)&=\ff 1 2 \int_{\R^d}\left\<\si^*(x)\nn W(\cdot,y)(x),\si^*(x)\nn f(x)\right\>\mu_\infty(\d x)\\
&=-\int_{\R^d}f(x)L_{\mu_\infty}W(\cdot,y)(x)\mu_\infty(\d x).
\end{align*} 
Then $\bar A$ is an integral operator with kernel $L_{\mu_\infty}W(\cdot,y)(x)$. Hence,  \eqref{LW-L2} yields that $\bar A$ is a compact operator. 

Next, we prove $\sup_{{\rm Re}z\geq 0}\|\bar A(L_{\mu_\infty}-z)^{-1}\|_{\cH}<1$, then the proof is complete. Since $L_{\mu_\infty}$ is self-adjoint and satisfies \eqref{Poincare}, for any $f\in\cH$ and $z$ with ${\rm Re}z\geq 0$, there is 
\beg{align*}
\|L_{\mu_\infty}(L_{\mu_\infty}-z)^{-1}f\|_{L^2_{\mathbb{C}}}&= \int_{\tld\la}^{+\infty} \left|\ff {\la} {\la+z}\right|^2\d\<E_{\la}f,f\>_{L^2_{\mathbb{C}}} \\
&\leq \int_{\tld\la}^{+\infty} \d\<E_{\la}f,f\>_{L^2_{\mathbb{C}}} \\
&= \|f\|_{L^2_{\mathbb{C}}},
\end{align*}
where $\{E_{\la}\}_{\la\geq 0}$ is the resolution of the identity for $-L_{\mu_\infty}$. Then
\beg{align*}
\bar A(L_{\mu_\infty}-z)^{-1}f(y)&=\ff 1 2\int_{\R^d}\left\<\si^*(x)\nn W(\cdot,y)(x),\si^*(x)(\nn (L_{\mu_\infty}-z)^{-1}f)(x) \right\>\mu_\infty(\d x)\\
&= - \int_{\R^d} W(x,y) L_{\mu_\infty}(L_{\mu_\infty}-z)^{-1}f(x)\mu_\infty(\d x).
\end{align*}
Consequently, 
\beg{align*}
\sup_{{\rm Re}z\geq 0}\|\bar A(L_{\mu_\infty}-z)^{-1}f\|_{L^2_{\mathbb{C}}}&\leq \left(\int_{\R^d\times\R^d}W(x,y)^2\mu_\infty(\d x)\mu_\infty(\d y)\right)^{\ff 1 2}\|L_{\mu_\infty}(L_{\mu_\infty}-z)^{-1}f\|_{L^2_{\mathbb{C}}}\\
&< \|W\|_{L^2(\mu_\infty\times\mu_\infty)}\|f\|_{L^2_{\mathbb{C}}}.
\end{align*}
This, together with \eqref{WL2}, implies $\sup_{{\rm Re}z\geq 0}\|\bar A(L_{\mu_\infty}-z)^{-1}\|_{\cH}<1$.

\end{proof}

To illustrate our main result, we give following examples, where  Poincar\'e inequalities are used.  The first example comes from \cite{Daw}, the second one comes from \cite{Co}. The third example is a two dimensional version of the second example, and the final one can not be covered by \cite{MR} since the log-Sobolev fails in this case. Proposition \ref{prop-22} and Proposition \ref{lem-2.3}  are used in the proofs of these examples. Detailed discussions are given in Section 6.
\beg{exa}\label{exa2}
Consider  Dawson's model in \cite{Daw}: 
\beg{equation}\label{eq-Daw}
\d X_t=-(X_t^3-X_t)\d t-\be\int_{\R}(X_t-y)\sL_{X_t}(\d y)\d t+\si\d B_t.
\end{equation}
Due to \cite[Theorem 3.3.1 and Theorem 3.3.2]{Daw}, there is $\si_c>0$ such that for $0<\si<\si_c$, \eqref{eq-Daw} has three stationary distributions, saying $\mu_+,\mu_S,\mu_-$, which satisfy
\[\int_{\R^d}x\mu_+(\d x)>0,\qquad \int_{\R^d}x\mu_{S}(\d x)=0,\qquad\int_{\R^d} x\mu_-(\d x)<0.\]
By using our theorem, there exist $C\geq 1$ and $\la>0$ such that for any $\sL_{X_0}\in\sP_1$ with 
\[\WW_1(\sL_{X_0},\mu_{\pm })<\ff 1 {4C^2}\left(\la\we\ff {1} {2}\right),\]
there is 
\[\WW_1(\sL_{X_t},\mu_{\pm})\leq 2Ce^{-\la t}\WW_1(\sL_{X_0},\mu_\pm).\] 
Moreover, for any $p_0\geq 1$ and $\ph_0(r):=r\we 1$, there exist  $C\geq 1$ and $\la>0$ such that 
\beg{equation}\label{th-loc-cov00} 
\|\sL_{X_{t}}-\mu_{\pm}\|_{p_0,\ph_0}<  2C e^{-\la t}\|\sL_{X_0}-\mu_{\pm}\|_{p_0,\ph_0},~t\geq 0,
\end{equation}
for $\sL_{X_0}\in\sP_{p_0}$ satisfying 
\beg{equation}\label{add-con-init00}
\|\sL_{X_0}-\mu_{\pm}\|_{p_0,\ph_0}<\ff {1} {4C^2}\left(\ff 1 2\we \la\right).
\end{equation}
\end{exa}
\beg{rem}
In \cite[Theorem 3.4]{Tug10}, only the weak convergence of  $\mu_t$ to $\mu_\pm$ was proved, and smooth density w.r.t.  the Lebesgue measure and finite moment of any order were required for $\mu_0$ in addition.
\end{rem}

\beg{exa}\label{exa1}
Let $\be\neq 0$. Consider the following equation
\beg{equation}\label{exa-Gaus}
\d X_t=-X_t\d t+\be \int_{\R}\cos(y)\sL_{X_t}(\d y)\d t+\sq 2\d B_t.
\end{equation}
Let $m$ satisfy 
\beg{equation}\label{eq-mm}
\left\{\beg{array}{c}
\cos(\be m)  =\sq e m,\\
\be\sin(\be m)> -\sq e.
\end{array}
\right.
\end{equation}
Then $\mu_m$, which is a Gaussian measure with mean $\be m$  and covariance $1$,  is a stationary distribution of  \eqref{exa-Gaus} and there exist $C\geq 1$ and $\la>0$ such that for any $\sL_{X_0}\in\sP_1$ with 
\[\WW_1(\sL_{X_0},\mu_{m})<\ff {1} {4C^2e^{\la}}\left(\la\we\ff {1} {2}\right),\]
there is 
\[\WW_1(\sL_{X_t},\mu_m)\leq 2Ce^{-\la t}\WW_1(\sL_{X_0},\mu_m).\]
Moreover, for any $p_0\geq 1$ and $\ph_0(r):=r\we 1$, there exist  $C\geq 1$ and $\la>0$ such that \eqref{th-loc-cov00} holds for $\sL_{X_0}\in\sP_{p_0}$ satisfying \eqref{add-con-init00} with $\mu_{\pm}$ replaced by $\mu_{m}$ in both inequalities.

\end{exa}

The following example is a two dimensional version of Example \ref{exa1}. 
\beg{exa}\label{exa3}
Let $\be\neq 0$. Consider the following system
\beg{equation}\label{exa-XY}
\left\{
\beg{array}{c}
\d X_t=-X_t\d t+\be \int_{\R}\cos(y)\sL_{Y_t}(\d y)\d t+\sq 2\d B_t^{(1)},\\
\d Y_t=-Y_t\d t+\be \int_{\R}\cos(y)\sL_{X_t}(\d y)\d t+\sq 2 \d B_t^{(2)}.
\end{array}
\right.
\end{equation}
Let $(m_1,m_2)$ satisfy 
\beg{equation}\label{XY-m1m2}
\left\{
\beg{array}{c}
\cos(\be m_1)=\sq e m_2,\\
\cos(\be m_2)=\sq e m_1.
\end{array}
\right.
\end{equation}
Then $\mu_m$, which is a two dimensional Gaussian measure with mean $m=(m_2\be,m_1\be)$ and covariance matrix equal to the identity, is a stationary distribution of \eqref{exa-XY}. Furthermore, if 
\beg{align}\label{bm1m2}
\be^2\sin(\be m_1)\sin(\be m_2)< e,
\end{align}
then conclusions in Example \ref{exa1} hold for $\mu_m$. 
\end{exa}

According to the concentration of measures for the log-Sobolev inequality, see e.g. \cite[Corollary 5.3.2]{WBook}, the final example can not be treated by \cite{MR}.
\beg{exa}\label{exa-2.5}
Let $\be> 0$. Consider 
\beg{align}\label{equ-2.5}
\d X_t & =\ff {1+\ff 1 3X_t^2} {(1+X_t^2)^{\ff 4 3}}\left(-\ff {X_t^3 } {1+X_t^2 } +\ff {(1-\be)X_t } {(1+X_t^2)^{\ff 1 3}}+\be\int_{\R} \ff {x} {(1+x^2)^{\ff 1 3}} \mu(\d x)\right)\d t\nonumber\\
&\quad\, -\ff {\si^2 X_t(1+\ff 1 9 X_t^2)} {(1+\ff 1 3 X_t^2)(1+X_t^2)}\d t+\si\d B_t.
\end{align}
Then there is $\si_c>0$ such that for $\si\in (0,\si_c)$, \eqref{equ-2.5} has three stationary distributions, saying $\mu_+,\mu_s,\mu_-$, which satisfy
\[\int_{\R^d}\ff x {(1+x^2)^{\ff 1 3}}\mu_+(\d x)>0,\qquad \int_{\R^d}x\mu_{S}(\d x)=0,\qquad\int_{\R^d}\ff x {(1+x^2)^{\ff 1 3}}\mu_-(\d x)<0.\]
Let  $\ph_0(r)=r\we 1$, 
\[V_0(x)=e^{(1+x^2)^{\ff 1 3}},~x\in\R.\] 
Then there exist $C\geq 1$ and $\la>0$ such that for any $\sL_{X_0}\in\sP_{V_0}$ with 
\[\|\sL_{X_0}-\mu_{\pm }\|_{V_0,\ph_0}<\ff 1 {4C^2}\left(\la\we\ff {1} {2}\right),\]
there is 
\[\|\sL_{X_t}-\mu_{\pm}\|_{V_0,\ph_0}\leq 2Ce^{-\la t}\|\sL_{X_0}-\mu_\pm\|_{V_0,\ph_0}.\]
\end{exa}

We give some remarks on the conditions of this theorem. The first remark is on $\|\cdot\|_{V_0,\ph_0}$ and  \eqref{con-W1}.
\beg{rem}\label{rem-WW-pp}
If $V_0\equiv 1$ and $\ph_0(r)=r$, then $\|\mu-\nu\|_{V_0,\ph_0}=\WW_1(\mu,\nu)$, where $\WW_1$ is the $L^1$-Wasserstein metric:
\beg{equation}\label{WW1}
\WW_1(\mu,\nu)=\inf_{\pi\in\sC(\mu,\nu)}\int_{\R^d\times\R^d}|x-y|\pi(\d x,\d y),~\mu,\nu\in\sP_1,
\end{equation}
and $\sC(\mu,\nu)$ consists of all couplings of $(\mu,\nu)$. 

We remark that \eqref{con-W1} holds if $X_t^{\mu_\infty}$ satisfies 
\beg{equation}\label{WW-p-ph}
\WW_{V_0,\ph_0}\left(\sL_{X_t^{\mu_\infty}(x)},\sL_{X_t^{\mu_\infty}(y)}\right)\leq C_We^{-\la_P t}\ph_0(|x-y|)\left(V_0(x)+V_0(y)\right),~x,y\in\R^d,
\end{equation}
where
\[\WW_{V_0,\ph_0}(\mu,\nu)=\inf_{\pi\in \sC(\mu,\nu)}\int_{\R^d\times\R^d}\ph_0(|x-y|)\left(V_0(x)+V_0(y)\right)\pi(\d x,\d y).\]
Indeed, let $\pi_{x,y}$ be the optimal coupling of $\left(\sL_{X_t^{\mu_\infty}(x)},\sL_{X_t^{\mu_\infty}(y)}\right)$. Then  
\beg{align*}
\left|P_t^{\mu_\infty}f(x)-P_t^{\mu_\infty}f(y)\right|&=\left|\int_{\R^d\times\R^d}(f(u)-f(v))\pi_{x,y}(\d u,\d v)\right|\\
&\leq  \int_{\R^d\times\R^d}\ph_0(|u-v|)\left(V_0(v)+V_0(u)\right)\pi_{x,y}(\d u,\d v)\\
&= \WW_{V_0,\ph_0}\left(\sL_{X_t^{\mu_\infty}(x)},\sL_{X_t^{\mu_\infty}(y)}\right)\\
&\leq C_{W}e^{-\la_P t}\ph_0(|x-y|)\left(V_0(x)+V_0(y)\right),~x,y\in\R^d,~f\in \sG_{V_0,\ph_0},
\end{align*}
which implies \eqref{con-W1}. 

The coupling method can be adopted to establish \eqref{WW-p-ph} for $X_t^{\mu_\infty}$.  We refer to \cite{EGZ,Wan23} and references therein for the exponential convergence of $X_t^{\mu_\infty}$ in $\WW_{V_0,\ph_0}$.

\end{rem}
The second remark is a sufficient condition for \eqref{hhh0}. 
\beg{rem}\label{rem-hhh0}
If there is a constant $C>0$ such that
\beg{align*}
|\nn b(\cdot,\mu)(x)-\nn b(\cdot,\nu)(x)|+|b(x,\mu)-b(x,\nu)|\leq C\|\mu-\nu\|_{V_0,\ph_0},
\end{align*}
then \eqref{hhh0} holds. Indeed, denoting 
\[h(x)=b(x,\mu)-b(x,\nu),\]
we have that 
\beg{align*}
|h(x_1)-h(x_2)|&\leq |b(x_1,\mu)-b(x_1,\nu)|+|b(x_2,\mu)-b(x_2,\nu)|\\
&\leq  2C\|\mu-\nu\|_{V_0,\ph_0},\\
|h(x_1)-h(x_2)|&\leq \left|\int_0^1\nn_{x_1-x_2}h(x_2+\th(x_1-x_2))\d\th \right|\\
&\leq C|x_1-x_2|\cdot\|\mu-\nu\|_{V_0,\ph_0},
\end{align*}
which imply that
\beg{equation}\label{hhh}
|h(x_1)-h(x_2)|\leq 2C(|x_1-x_2|\we 1)\|\mu-\nu\|_{V_0,\ph_0}.
\end{equation}
This, together with  \eqref{Cr<ph}, yields \eqref{hhh0}. 
\end{rem}
The final remark is a sufficient condition for \eqref{DFb-pph}.
\beg{rem}\label{rem-DFb-pph}
If $D^F_{\mu}b$ satisfies 
\beg{align*}
\left| D^F_{\mu}b(x,z)- D^F_{\nu}b(x,z)\right|&\leq C\|\mu-\nu\|_{V_0,\ph_0}V_0(z),~x,z\in\R^d,\\
\left|\nn^{(2)} D^F_{\mu}b(x,\cdot)(z)-\nn^{(2)} D^F_{\nu}b(x,\cdot)(z)\right|&\leq C\|\mu-\nu\|_{V_0,\ph_0}V_1(z),~x,z_1,z_2\in\R^d,
\end{align*}
where $\nn^{(2)}$ is the gradient w.r.t. the second variable $z$ of $D^F_{\mu}b(x,z)$ for any $\mu$, and $V_1$ satisfies
\[\sup_{\th\in [0,1]}V_1(\th x+(1-\th)y)\leq C(V_0(x)+V_0(y))\]
then $D^F_{\mu}b(x,\cdot)-D^F_{\nu}b(x,\cdot)$ satisfies \eqref{DFb-pph}.\\
Indeed, 
\beg{align*}
& \left|\left(D^F_{\mu}b(x,z_1)-D^F_{\nu}b(x,z_1)\right)-\left(D^F_{\mu}b(x,z_2)-D^F_{\nu}b(x,z_2)\right)\right| \\
&\quad\, = \left|\int_0^1\nn_{z_1-z_2}^{(2)}\left(D^F_{\mu}b(x,\cdot)-D^F_{\nu}b(x,\cdot)\right)(z_2+\th(z_1-z_2))\d\th\right|\\
&\quad\, \leq C|z_1-z_2| \int_0^1\|\mu-\nu\|_{V_0,\ph_0}V_1(z_1+\th(z_1-z_2))\d\th\\
&\quad\,\leq C\|\mu-\nu\|_{V_0,\ph_0}|z_1-z_2|\left(V_0(z_1)+V_0(z_2)\right),
\end{align*}
and 
\beg{align*}
&\left|\left(D^F_{ \mu }b(x,z_1)-D^F_{\nu}b(x,z_1)\right)-\left(D^F_{\mu}b(x,z_2)-D^F_{\nu}b(x,z_2)\right)\right|\\
&\quad\, \leq 2C\|\mu-\nu\|_{V_0,\ph_0}\left(V_0(z_1)+V_0(z_2)\right).
\end{align*}
Thus, together with \eqref{Cr<ph}, we have that
\beg{align*}
&\left|\left(D^F_{\mu}b(x,z_1)-D^F_{\nu}b(x,z_1)\right)-\left(D^F_{\mu}b(x,z_2)-D^F_{\nu}b(x,z_2)\right)\right|\\
&\quad\, \leq C\|\mu -\nu\|_{V_0,\ph_0}(|z_1-z_2|\we 1)\left(V_0(z_1)+V_0(z_2)\right)\\
&\quad\, \leq C\|\mu -\nu\|_{V_0,\ph_0}\ph_0(|z_1-z_2|)\left(V_0(z_1)+V_0(z_2)\right).
\end{align*}
\end{rem}

\section{Regularity of diffusions with polynomial growth coefficients}
The equation \eqref{equ-mu-inf} is a distribution free SDE with a polynomial growth drift term. We study the regularity of this kind of SDEs and associated Markov semigroups, which will be used in following sections. To this end, we consider 
\beg{align}\label{equ-Y}
\d Y_t=Z(Y_t)\d t+\si(Y_t)\d B_t,
\end{align}
where  $\si$ and $Z$ satisfy following assumptions.  
\beg{description}[align=left, noitemsep]
\item [(A1)] The drift term $Z\in C^1(\R^d,\R^d)$ and $|\nn Z|$ has polynomial growth. There exist a nonnegative constant $K_1$ and a nonnegative function $K_2\in C^{0}$  such that for all $x,v\in\R^d$
\beg{align}\label{Inequ-nnZ1}
\<\nn_v Z(x),v\>\leq (K_1-K_2(x))|v|^2.
\end{align}
\end{description}
Let 
\beg{equation}\label{lim-K2}
K_2^*(r)=\left\{\beg{array}{cc}
r^{-1}\displaystyle\inf_{|v|=1,x\in\R^d}\left(\int_0^r K_2(x+\th v)\d \th\right),& r>0,\\[15pt]
\displaystyle\inf_{x\in\R^d}K_2(x),& r=0.
\end{array}\right.
\end{equation}
Then $K_2^*$ is a nonnegative and locally bounded function on $[0,+\infty)$. We also assume that 
\beg{description}
\item [(A2)] The diffusion term $\si\in C^1(\R^d,\R^d\otimes\R^d)$ has bounded and continuous first derivative, and  there is a positive constant $\si_0$ such that \eqref{nondege0} holds. There exist $\al_1\in [0,1]$, $\al_2>1$ and $K_\si>0$ such that
\beg{align}
\|\si(x)\|_{HS}&\leq K_{\si}(1+|x|)^{\al_1},~x\in\R^d,\label{sub-lin-si}\\
\varlimsup_{r\ra +\infty}&\ff  {(1+r)^{(\al_1+1)\al_2}}  {K_2^*(r)r^2}<+\infty.\label{lim-al1al2}
\end{align}
\end{description}

\beg{rem}
We can assume that $K_2\geq  c$ for some positive constant. In fact, if $K_2(\cdot)\equiv 0$, we can replace $K_1$ by $K_1+c$ and set $K_2(x)\equiv c$, then \eqref{Inequ-nnb1} implies that {\bf (A1)} holds with $K_1=K_b+c$ and $K_2(x)\equiv c$.

\eqref{Inequ-nnZ1} and $K_2^*$ can be used to model some drift which is dissipative in long distance. For instance, $Z=-|x|^{\be}x$ for some $\be>0$ and large $x$, then $K_2(x)=K_2|x|^{\be}$ for some positive constant $K_2$. Conversely, if $Z\in C^1(\R^d,\R^d)$ satisfies \eqref{Inequ-nnZ1} with $K_2(x)=K_2|x|^{\be}$, then $Z$ is dissipative in long distance: there is $\tld K_2>0$ such that
\[\<Z(y)-Z(x),y-x\>\leq \left(K_1-\tld K_2|y-x|^{\be }\right)|y-x|^2.\]
This proof is elementary, so we omit it. 
\end{rem}

We denote by $\nn Y_{t}$ and $\nn^2Y_t$ the first and second order gradient of $Y_t$  w.r.t. initial value, respectively. We denote by $D Y_{t}$ the Malliavin derivative. To emphasise the initial value, we use $Y_t^y,\nn Y_t^y$, $\nn^2 Y_t^y$ and $D Y_t^y$ when $Y_0=y$.    

\beg{lem}\label{lem-Malli}
Assume  that {\bf (A1)} and {\bf (A2)} hold except \eqref{nondege0}. Then \\
(1) Equation \eqref{equ-Y} has a unique strong solution, and for any $p\geq 1$, there is $C_p>0$ such that
\beg{equation}\label{Eest-Y-2p}
\E\left(\sup_{s\in [0,t]}|Y_t|^{2p}+ \int_0^t(K^*_2(|Y_s|)+1)|Y_s|^{2p}\d s\Big|\sF_0\right)\leq C_pe^{C_{p} t}\left(\sq t(\sq t+1)+ |Y_0|^{2p}\right).
\end{equation}

(2) Let $v\in\R^d$ and $Y_0=y\in\R^d$. $Y_{t}$ is differentiable w.r.t. initial value along $v$, and the derivative process  $\nn_v Y_t$ satisfies
\[\d \nn_vY_t=\nn_{\nn_vY_t }Z(Y_{t})\d t+\nn_{\nn_vY_t}\si(Y_{t})\d B_t,~\nn_vY_0=v\in\R^d, t\geq 0.\]
Let $\{h_{t}\}_{t\geq 0}$ be an adapted process with $h'_{t}=w(t)\nn_v Y_{t}$, where $w(\cdot)$ is a $\R^d\otimes\R^d$-valued, bounded and adapted process on $[0,\infty)$. Then $Y_{t}$ is Malliavin differentiable along $h$,  and the derivative process $D_h Y_{t}$  satisfies 
\bequ\label{equ-DhY}
\d D_h Y_{t}=\nn_{D_h Y_{t}}Z(Y_{t})\d t+\nn_{D_h Y_{t}}\si(Y_{t})\d B_t+\si(Y_{t})h'_{t}\d t,~D_h Y_{0}=0.
\enqu
For $Y_t,\nn_vY_t,D_h Y_t$, we have following estimates for any $p\geq 1$
\beg{align}
\E\left(\sup_{t\in[0,s]}|\nn_vY_{t}|^{2p}+\int_0^s(K_2(Y_t)+1)|\nn_vY_t|^{2p}\d t\right)&\leq C_pe^{C_ps}|v|^{2p},\label{Eest-nnY-2p}\\
\E\left(\sup_{t\in[0,s]}|D_hY_{t}|^{2p}+\int_0^s(K_2(Y_t)+1)|D_hY_t|^{2p}\d t\right)&\leq C_pe^{C_ps}\E \int_0^s |\si(Y_t)h'_t|^{2p}\d t,~s\geq 0,\label{Eest-DhY-2p}
\end{align}
where $C_p$ is a positive constant which is different in different formula.

(3) If, furthermore,  $Z\in C^2(\R^d,\R^d)$, $\si$ has bounded and continuous second order derivatives, and there are nonnegative constants $K_3$ and $\be_1$ such that 
\beg{equation}\label{Inequ-nnZ2}
|\nn^2 Z(x)|\leq K_3(1+|x|)^{\be_1},~x\in\R^d,
\end{equation}  
then, the second order derivative of $Y_t$ w.r.t. the initial value exists and satisfies
\bequ\label{Inequ-nn2Y}
\sup_{|u|,|v|\leq 1}\E\sup_{t\in [0,s]}|\nn_u\nn_v Y_t^y|^{2p}\leq C_p e^{C_p s}(1+|y|)^{2p\be_1},~s\geq 0
\enqu 
for some $C_p$ depending on  $p,K_1,K_3,\be_1,\al_1,\al_2,\|\nn\si\|_\infty,\|\nn^2\si\|_\infty$.\\

\end{lem}
The proof of this lemma is fundamental, and we give it in the appendix.  When $\nn Z,\nn\si$ or $\nn^2 Z$ and $\nn^2 \si$ are continuous,   we have the following corollary.  
\beg{cor}\label{cor-nnY-con}
Let $v,u\in\R^d$. Assume {\bf (A1)} and {\bf (A2)} except \eqref{nondege0}. Then for any $p\geq 1$ 
\beg{align}\label{cont-nn1}
\lim_{x\ra y}\E\sup_{t\in [0,s]}| \nn_v Y_t^y- \nn_v Y_t^x|^{2p}&=0,~s>0,\\
\lim_{x\ra y,s\ra t}\E|\nn_vY_s^x-\nn_v Y_t^y|^{2p}&=0.\label{cont-nn1-ad}
\end{align}
If, furthermore,  the assumption of (3) in Lemma \ref{lem-Malli} holds, then
\beg{align}\label{cont-nn2}
\lim_{x\ra y}\E\sup_{t\in [0,s]}|\nn_u\nn_v Y_t^y-\nn_u\nn_v Y_t^x|^{2p}&=0,~s>0\\
\lim_{x\ra y,s\ra t}\E|\nn_u\nn_v Y_s^x-\nn_u\nn_v Y_t^y|^{2p}&=0.\label{cont-nn2-ad}
\end{align}

\end{cor}
This corollary can be proved by using estimates in Lemma \ref{lem-Malli} and the dominated convergence theorem, and we omit it.

Let $P_t^Z$ be the Markov semigroup associated with $Y_t$. Next, we introduce a lemma on the gradient estimate of $P_{t}^Z$ for $t>0$. The method to prove this lemma is due to \cite[Lemma 2.1]{Wang16} essentially. 
\beg{lem}\label{lem-grad-est}
Assume that {\bf (A1)}, {\bf (A2)} and the assumption of (3) of Lemma \ref{lem-Malli} hold. For any $t>0$ and any $f\in\sB(\R^d)$ with $|f|\leq V$, where  $1\leq V\in C^{0}$ is an arbitrary function satisfying 
\beg{align}\label{EVp}
\sup_{s\in [0,t],y\in\R^d}\ff { P^Z_t V^p(y)} {V^p(y)}<+\infty,~t>0,~p\geq1,
\end{align}
then there is $P_t^Zf\in C^2$ and
\beg{equation}\label{Bismut} 
\nn_vP_t^Zf(y)=\ff 1 t\E\left(f(Y_t^{y})\int_0^{t}\<\si^{-1}(Y_{ r}^{y})\nn_v Y_{ r}^{y},\d B_r\>\right).
\end{equation}
Moreover,  for any $q\in(1,+\infty)$,  there is $C>0$ such that
\beg{align}
|\nn P_t^Zf|(y)&\leq \ff {C} {\si_0\sq {t\we 1}} \left(P_t^Z|f|^q\right)^{\ff 1 q}(y),\label{nnpzf}\\
|\nn^2 P_{t}^Z f|(y)&\leq \ff {C} {\si_0^2 (t\we 1)}\left(P^Z_{t}|f|^q\right)^{\ff 1 q}(y),~t>0,y\in\R^d.\label{grad-2-est}
\end{align}
Furthermore, if $f\in \sG_{V_1,\ph}$ for some  $\ph\in\sG$ and $V_1\geq 1$ satisfying \eqref{EVp} and 
\beg{align}\label{ZnnV}
\< Z(x),\nn V_1(x)\>\leq C V_1(x),~x\in\R^d,
\end{align}
then $P_t^Zf\in C^2(\R^d)$ satisfies  \eqref{Bismut}-\eqref{grad-2-est} and 
\beg{align}\label{nn1Pph}
|\nn P_{t}^Z f|(y)&\leq\ff {Ce^{Ct}\|f\|_{V_1,\ph}} {\si_0}V_1(y)\left(1+|y|\right)^{\al_1\we \ff {2\al_1} {(1+\al_1)\al_2}}t^{-\ff 1 2}\ph\left(t^{\ff 1 2\we \ff {(1+\al_1)\al_2-2\al_1} {2(1+\al_1)\al_2}}\right),\\
|\nn^2 P_{t}^Z f|(y)&\leq \ff {Ce^{Ct}\|f\|_{V_1,\ph}} {\si_0^2}V_1(y)\left(1+|y|\right)^{\al_1\we \ff {2\al_1} {(1+\al_1)\al_2}}t^{-1}\ph\left(t^{\ff 1 2\we \ff {(1+\al_1)\al_2-2\al_1} {2(1+\al_1)\al_2}}\right).\label{nn2Pph}
\end{align}
\end{lem}

\beg{proof}[{Proof of Lemma \ref{lem-grad-est}}]
For any $t> 0$ and $y,v\in \R^d$. It follows from Lemma \ref{lem-Malli} that $Y_t^y$ is differentiable w.r.t. initial value. We set $h_{r}=\ff 1 {t} \int_0^r\si^{-1}(Y_{s}^y)\nn_v Y_{s}^y\d s$.  Due to Lemma \ref{lem-Malli} again, $Y_t$ is Malliavin differentiable along  $h$ on $[0,t]$.  Moreover,  $\{D_h Y_{r}^y\}_{r\in [0,t]}$ and $\{\ff {r} {t}\nn_v Y_{r}^y\}_{r\in [0,t]}$ satisfy  the same equation
$$\d \xi_{r}=\nn_{\xi_{ r}}Z(Y_{ r}^y)\d r+\nn_{\xi_{ r}}\si(Y_{ r}^y)\d B_r+\ff 1 {t} \nn_vY_{ r}^y\d r,~\xi_0=0,~r\in [0,t].$$
Thus $D_hY_t^y=\nn_vY_t^y$, and for $f\in C_b^1$ 
\beg{align*}
\nn_v P^Z_{ t}f(y)&=\E \<\nn f(Y_{ t}^y), \nn_v Y_{ t}^y\>=\E \<\nn f(Y_{ t}^y), D_h Y_{ t}^y\>\\
&=\E D_h( f(Y_{ t}^y))=\ff 1 {t}\E f(Y_{t}^y)\int_0^t\<\si^{-1}(Y_{ r}^y)\nn_v Y_{ r}^y,\d B_r\>.
\end{align*}
This, together with $f\in C^1_b$, \eqref{nondege0}, \eqref{Eest-nnY-2p} and \eqref{cont-nn1}, yields by the dominated convergence theorem that $\nn_v P_{t}^Zf\in C^0$. Moreover,  
\beg{align}\label{Bis-int}
&P^Z_{t}f(y+v)-P^Z_{t}f(y) =\int_0^1 \nn_v P_{t}^Z f(y+\th v)\d \th\nonumber\\
&\quad\,= \ff 1 {t }\int_0^1\left(\E f(Y_{ t}^{y+\th v})\int_0^t\<\si^{-1}(Y_{ r}^{y+\th v})\nn_v Y_{ r}^{y+\th v},\d B_r\>\right)\d \th.
\end{align}

Next, we prove \eqref{Bismut} and \eqref{nnpzf} for $f$ with $|f|\leq V$ and $V$ satisfying \eqref{EVp}. For any $p\geq 1$, following from \eqref{nondege0}, \eqref{Eest-nnY-2p}, the B-D-G inequality and the Minkowski inequality, there are positive constants $C_{1}$ which depends only on $p$, and $C_{2}$ which depends on $p,K_1,\|\nn\si\|_\infty$, such that
\beg{equation}\label{in-MM}
\beg{split}
\E\sup_{s\in[0,t]}\left|\int_0^s\<\si^{-1}(Y_{ r}^{y+\th v})\nn_v Y_{ r}^{y+\th v},\d B_r\>\right|^{p}&\leq C_{1} \E\left(\int_0^t |\si^{-1}(Y_{ r}^{y+\th v})\nn_v Y_{ r}^{y+\th v}|^2\d r\right)^{\ff p 2}\\
&\leq \ff {C_{1 }} {\si_0^p} \left(\int_0^t \left(\E|\nn_v Y_{r}^{y+\th v}|^{p\vee 2}\right)^{\ff 2 {p\vee 2}}\d r \right)^{\ff {p} 2}\\
&\leq \ff {C_{1 }} {\si_0^p}\left(\int_0^t e^{2C_2r}|v|^{2}\d r \right)^{\ff {p} 2}\\
&\leq  C_{1 }\si_0^{-p}t^{\ff p 2}e^{pC_2t} |v|^p.
\end{split}
\end{equation}
Thus, letting
\[\nu_1(A)=\int_0^1 \E\left(\1_{A}(Y_t^{y+\th v})\left|\int_0^t\<\si^{-1}(Y_{ r}^{y+\th v})\nn_v Y_{ r}^{y+\th v},\d B_r\>\right|\right)\d\th,~A\in\sB(\R^d),\]
$\nu_1$ is a finite measure.  Let $\nu=\nu_1+\sL_{Y_t^{y+v}} +\sL_{Y_t^y}$. Since $C_b^1$ is dense in $L^1(\nu)$, the equality \eqref{Bis-int} can be extended to any  $f\in L^1(\nu)$. Because $V$ satisfies \eqref{EVp} and $|f|\leq V$, for any $p\geq 1$, there is 
\[\E |f|^p(Y_t^y)\leq \E V(Y_t^y)^p\leq C_t V(y)^p,~y\in\R^d,~t\geq 0\]
where $C_t>0$ is locally bounded on $t$. Then \eqref{in-MM} and  \eqref{Eest-Y-2p} imply that $f\in L^1(\nu)$.  Moreover,  there are positive constants $C_3,C_4$ such that
\beg{equation}\label{PPint}
\beg{split}
&|P^Z_{t}f(y+v)-P^Z_{t}f(y)|\\
&\quad\,\leq \ff 1 t \int_0^1\left|\E f(Y_{ t}^{y+\th v})\int_0^t\<\si^{-1}(Y_{ r}^{y+\th v})\nn_v Y_{ r}^{y+\th v},\d B_r\>\right|\d \th\\
&\quad\,\leq \ff 1 t \int_0^1\left(\E V(Y_t^{y+\th v})^{2}\right)^{\ff 1 2}\left(\E\left|\int_0^t\<\si^{-1}(Y_{ r}^{y+\th v})\nn_v Y_{ r}^{y+\th v},\d B_r\>\right|^2\right)^{\ff 1 2}\d\th\\
&\quad\, \leq \ff 1 {\sq t}\left( C_t\sup_{0\leq \th\leq 1}V(y+\th v)\right) \sq{C_{1 }} \si_0^{-1}e^{C_2t} |v|.
\end{split}
\end{equation}
This implies that $P_t^Zf\in C^0$. Replacing $v$ by $\ep v$ for $\ep>0$, for any $0<\de<t$,  there is
\beg{align*}
&\ff {P_t^Z f(y+\ep v)-P_t^Z f(y)} {\ep} = \ff 1 {\ep}\left( P_{t-\de}^Z (P_{\de}^Z f)(y+\ep v)-P_{t-\de}^Z(P_{\de}^Zf)(y)\right)\\
&\quad\,=\ff 1 {\ep (t-\de)}\int_0^{\ep} \left(\E P_{\de}^Zf(Y_{ t-\de}^{y+\th v})\int_0^{t-\de}\<\si^{-1}(Y_{ r}^{y+\th v})\nn_v Y_{ r}^{y+\th v},\d B_r\>\right)\d \th.
\end{align*}
Due to Lemma \ref{lem-Malli} and Corollary \ref{cor-nnY-con}, $Y_t^y,~\nn Y_t^y$ are continuous in $y$. Combining this with that $V$ is locally bounded,   \eqref{EVp}, $P_{\de}^Zf\in C^0$ and
\beg{align}\label{plygrow}
P_{\de}^Z|f|^p(y)\leq C_{\de}V^p(y),~y\in\R^d,~p\geq 1,
\end{align}
the dominated convergence theorem yields
\beg{equation}\label{bis1}
\beg{split}
\nn_v P_t^Z(y)&=\lim_{\ep\ra 0^+}\ff {P_t^Z f(y+\ep v)-P_t^Z f(y)} {\ep}\\
&= \ff 1 {t-\de} \E P_{\de}^Zf(Y_{ t-\de}^{y})\int_0^{t-\de}\<\si^{-1}(Y_{ r}^{y})\nn_v Y_{ r}^{y},\d B_r\>.
\end{split}
\end{equation}
By using Lemma \ref{lem-Malli}, Corollary \ref{cor-nnY-con} and $P_{\de}^Zf\in C^0$ again, we can derive from the dominated convergence that $\nn_v P_t^Zf(y)$ is continuous in $(v,y)\in \R^d\times\R^d$. Moreover, by using the Markov property
\[P_{\de}^Zf(Y_{ t-\de}^{y})=\E[f(Y_{\de}^{Y_{t-\de}^y})|\sF_{t-\de}]=\E[f(Y_t^{y})|\sF_{t-\de}],\]
we have that 
\beg{align*}
&\E P_{\de}^Zf(Y_{ t-\de}^{y})\int_0^{t-\de}\<\si^{-1}(Y_{ r}^{y})\nn_v Y_{ r}^{y},\d B_r\>\\
&\quad\,=\E \left(\E[f(Y_t^{y})|\sF_{t-\de}]\int_0^{t-\de}\<\si^{-1}(Y_{ r}^{y})\nn_v Y_{ r}^{y},\d B_r\>\right)\\
&\quad\,=\E \left(\E\left[f(Y_t^{y})\int_0^{t-\de}\<\si^{-1}(Y_{ r}^{y})\nn_v Y_{ r}^{y},\d B_r\>\Big|\sF_{t-\de}\right]\right)\\
&\quad\, =  \E\left(f(Y_t^{y})\int_0^{t-\de}\<\si^{-1}(Y_{ r}^{y})\nn_v Y_{ r}^{y},\d B_r\>\right).
\end{align*}
Putting this into \eqref{bis1} and letting $\de\ra 0^+$, we obtain \eqref{Bismut}. Combining \eqref{Bismut} with \eqref{Eest-nnY-2p}, \eqref{in-MM} and the semigroup property of $P_t^Z$, we can prove \eqref{nnpzf} (the proof is similar to that of \eqref{grad-2-est} below, and we omit it).

Next, we prove that $P^Z_{t}f\in C^2$. By \eqref{bis1} with $\de=\ff t 2$, 
\[\nn_v P^Z_{t}f(y)=\ff 2 {t}\E P^Z_{\ff {t} 2 }f(Y_{ \ff { t} 2}^y)\int_0^{\ff { t} 2}\<\si^{-1}(Y_{ r}^y)\nn_v Y_{ r}^y,\d B_r\>.\]
Noting that $P_{\ff t 2}^Zf\in C^1$ and satisfies \eqref{EVp} due to \eqref{plygrow}, and taking into account  \eqref{Eest-nnY-2p}, \eqref{Inequ-nn2Y} and
\beg{equation}\label{nnsi-1}
|\nn  \si^{-1}(x)|=|\si^{-1}(x)(\nn\si(x))\si^{-1}(x)| \leq \si_0^{-2}|\nn\si(x)|\leq\si_0^{-2}\|\nn\si\|_\infty,
\end{equation}
we can derive by using the dominated convergence theorem that 
\beg{equation}\label{nn2P}
\beg{split}
\nn_u\nn_v P_{t}f(y)= & \ff 2 {t}\E \nn_{\nn_u Y_{ \ff {t} 2}^y}P^Z_{\ff { t} 2}f(Y_{ \ff { t} 2}^y)\int_0^{\ff { t} 2}\<\si^{-1}(Y_{ r}^y)\nn_v Y_{ r}^y,\d B_r\> \\
&+\ff 2 {t }\E P^Z_{\ff { t} 2 }f(Y_{\ff {t} 2}^y)\int_0^{\ff {t} 2}\<\nn_{\nn_u Y_{r}^y}\si^{-1}(Y_{r}^y)\nn_v Y_{ r}^y,\d B_r\> \\
&+\ff 2 {t }\E P^Z_{\ff { t} 2}f(Y_{ \ff { t} 2}^y)\int_0^{\ff {t} 2}\<\si^{-1}(Y_{r}^y)\nn_u\nn_v Y_{r}^y,\d B_r\>,
\end{split}
\end{equation}
and $\nn_u\nn_v P_t^Zf\in C^0$. By using \eqref{Eest-nnY-2p}, \eqref{in-MM}, \eqref{nnpzf} and \eqref{nnsi-1},  for any $q\in (1,\infty)$,  there is $C>0$ such that 
\[|\nn_u\nn_v P_t^Zf(y)|\leq \ff {C} {\si_0^2 t}|u|\cdot |v|\left(P_t^Z|f|^q\right)^{\ff 1 q}(y),~t\in (0,1].\]
Combining this with the semigroup property of $P_t^Z$, we have for $t>1$
\beg{align*}
|\nn^2P_t^Zf(y)|&=|\nn^2P_1^ZP_{t-1}^Zf(y)|\leq \ff {Ce^{C}} {\si_0^2}\left(P_1^Z|P_{t-1}^Zf|^q\right)^{\ff 1 q}(y)\\
&\leq \ff {Ce^{C}} {\si_0^2}\left(P_1^ZP_{t-1}^Z|f|^q\right)^{\ff 1 q}(y)=\ff {Ce^{C}} {\si_0^2}\left(P_{t}^Z|f|^q\right)^{\ff 1 q}(y).
\end{align*}
Thus \eqref{grad-2-est} holds.

Let $\xi_t$  be the solution of the following equation and $\xi_0=y$
\beg{equation}\label{eq-xi}
\d \xi_t=Z(\xi_t)\d t.
\end{equation}
Then, due to \eqref{ZnnV} and {\bf (A1)},
\beg{align}\label{V1et}
&V_1(\xi_t)\leq e^{Ct}V_1(y),\\
|\xi_t|^2+2\int_0^t &K^*_2(|\xi_s|)|\xi_s|^2\d s\leq e^{2K_1t}|y|^2,~t\geq 0.\nonumber
\end{align}
Then 
\beg{equation*}
\beg{split}
\d |Y_t^y-\xi_t|^2&\leq 2\left(K_1- K^*_2(|Y_t^y -\xi_t|)\right)|Y_t^y-\xi_t|^2\d t+\|\si(Y_t^y)\|_{HS}^2\d t+2\<Y_t^y-\xi_t,\si(Y_t^y)\d B_t\>\\
&\leq 2\left(K_1-K^*_2(|Y_t^y -\xi_t|)\right)|Y_t^y-\xi_t|^2\d t+K_\si^2(1+|Y_t^y|)^{2\al_1}\d t\\
&\quad\,+2\<Y_t^y-\xi_t,\si(Y_t^y)\d B_t\>\\
&\leq 2\left(K_1-K^*_2(|Y_t^y -\xi_t|)\right)|Y_t^y-\xi_t|^2\d t+K_\si^2(1+|Y_t^y-\xi_t|+|\xi_t|)^{2\al_1}\d t\\
&\quad\,+2\<Y_t^y-\xi_t,\si(Y_t^y)\d B_t\>.\\
&\leq 2\left(K_1-K^*_2(|Y_t^y -\xi_t|)\right)|Y_t^y-\xi_t|^2\d t+2^{(2\al_1-1)^+}K_\si^2(1+|Y_t^y-\xi_t|^{2\al_1})\d t\\
&\quad\,+2^{(2\al_1-1)^+}K_\si^2(1+|\xi_t|)^{2\al_1}\d t+2\<Y_t^y-\xi_t,\si(Y_t^y)\d B_t\>.
\end{split}
\end{equation*}
Combining this with \eqref{lim-al1al2} and \eqref{rr_0}, there is $C>0$ which is independent of $y$ and whose value may vary at each line such that
\beg{equation*}
\beg{split}
\E |Y_s^y-\xi_s|^2 &\leq Ce^{Cs}\int_0^s(1+|\xi_t|)^{2\al_1}\d t\\
&\leq Ce^{Cs}s^{ \ff {(1+\al_1)\al_2-2\al_1} {(1+\al_1)\al_2}}\left(\int_0^s(1+|\xi_t|)^{(1+\al_1)\al_2}\d t\right)^{\ff {2\al_1} {(1+\al_1)\al_2}}\\
&\leq  Ce^{Cs}s^{ \ff {(1+\al_1)\al_2-2\al_1} {(1+\al_1)\al_2}}\left(s+\int_0^s K^*_2(|\xi_t|)|\xi_t|^2\d t\right)^{\ff {2\al_1} {(1+\al_1)\al_2}}\\
&\leq C e^{Cs}\left(s+ s^{ \ff {(1+\al_1)\al_2-2\al_1} {(1+\al_1)\al_2}}|y|^{\ff {4\al_1} {(1+\al_1)\al_2}}\right)
\end{split}
\end{equation*}
and
\beg{align*}
\E |Y_s^y-\xi_s|^2 &\leq Ce^{Cs}\int_0^s(1+|\xi_t|)^{2\al_1}\d t\\
&\leq 2^{(2\al_1-1)^+}Ce^{Cs}\int_0^s\left(1+|\xi_t|^{2\al_1}\right)\d t\\
&\leq 2^{(2\al_1-1)^+}Ce^{Cs}\int_0^s\left(1+e^{2\al_1 K_1t}|y|^{2\al_1}\right)\d t\\
&=2^{(2\al_1-1)^+}Ce^{Cs} \left(s+\ff {e^{2\al_1 K_1s}-1} {2\al_1K_1}|y|^{2\al_1}\right).
\end{align*}
Consequently 
\beg{align}\label{Yxit}
\E |Y_s^y-\xi_s|^2&\leq  Ce^{Cs}\left\{s+\left(s|y|^{2\al_1}\right)\we\left(s^{ \ff {(1+\al_1)\al_2-2\al_1} {(1+\al_1)\al_2}}|y|^{\ff {4\al_1} {(1+\al_1)\al_2}}\right)\right\}\nonumber\\
&\leq  C e^{ Cs}s^{1\we \ff {(1+\al_1)\al_2-2\al_1} {(1+\al_1)\al_2}}\left(1+|y|^{2\al_1}\we |y|^{\ff {4\al_1} {(1+\al_1)\al_2}}\right)\nonumber\\
&\leq  C e^{ Cs}s^{1\we \ff {(1+\al_1)\al_2-2\al_1} {(1+\al_1)\al_2}}\left(1+|y|\right)^{2\al_1\we \ff {4\al_1} {(1+\al_1)\al_2}}.
\end{align}
Since $\ph^2(\sq{\cdot})$ is concave, we have 
\[\ff {\ph^2(\sq{ r^2})} {1+r^2}\leq \ff {\ph^2(\sq{ r^2})} {1+r^2}+\ff {r^2\ph^2(\sq{ r^{-2}})} {1+r^2}\leq \ph^2(\sq{1})=\ph^2(1),~r>0.\]
Thus,  there is $C>0$ such that
\[\ph(r)=\sq{\ph^2(\sq{r^2})}\leq \sq{C(1+r^2)},~r\geq 0.\]
Then
\beg{align*}
|f(y)-f(x)|\leq \ph(|x-y|)\left(V_1(x)+V_1(y)\right)\leq \sq{C(1+|x-y|^2)}\left(V_1(x)+V_1(y)\right),~x,y\in\R^d.
\end{align*}
Then $|f|\leq V_2$  with  $V_2(x):=C(1+|x|)V_1(x)$. The H\"older inequality, \eqref{Eest-Y-2p} and that $V_1$ satisfies \eqref{EVp} imply \eqref{EVp} holds for $V_2$.
Thus, $P_t^Zf\in C^2$ and \eqref{Bismut}-\eqref{grad-2-est} hold. By using \eqref{Eest-Y-2p}, \eqref{EVp} for $V_1$,  \eqref{Yxit}, \eqref{V1et} and $\ph^2(\sq{\cdot})$ is concave, we have for $1<p_2<2$ that 
\beg{align}\label{PZf-f}
(P^Z_{t}|f-f(\xi_{t})|^{p_2})^{\ff 1 {p_2}}(y)&=(\E |f(Y_t^y)-f(\xi_{t})|^{p_2})^{\ff 1 {p_2}}\nonumber\\
&\leq  \|f\|_{V_1,\ph}\left(\E\ph(|Y_{t}^y-\xi_{t}|)^{2}\right)^{\ff 1 2}\left[\left(\E V_1(Y_t^y)^{\ff {2p_2} {2-p_2}}\right)^{\ff {2-p_2} {2p_2}}+V_1(\xi_t)\right]\nonumber\\
& \leq  Ce^{Ct}\|f\|_{V_1,\ph} \ph\left(\left(\E|Y_{t}^y-\xi_{t}|^2\right)^{\ff 1 2}\right)V_1(y) \nonumber\\
&\leq Ce^{Ct}\|f\|_{V_1,\ph} \left(1+|y|\right)^{\al_1\we \ff {2\al_1} {(1+\al_1)\al_2}} V_1(y) \ph\left(t^{\ff 1 2\we \ff {(1+\al_1)\al_2-2\al_1} {2(1+\al_1)\al_2}}\right).
\end{align}
where in the last inequality we have used 
\beg{align*}
&\ph\left(C e^{ Ct} t^{\ff 1 2\we \ff {(1+\al_1)\al_2-2\al_1} {2(1+\al_1)\al_2}}\left(1+|y|\right)^{\al_1\we \ff {2\al_1} {(1+\al_1)\al_2}}\right)\\
&\quad\,\leq   Ce^{ Ct}\left(1+|y|\right)^{\al_1\we \ff {2\al_1} {(1+\al_1)\al_2}}\ph\left(t^{\ff 1 2\we \ff {(1+\al_1)\al_2-2\al_1} {2(1+\al_1)\al_2}}\right),
\end{align*}
since $\ph^2(\sq{\cdot})$ is concave. Then, \eqref{PZf-f} and \eqref{nnpzf} with $f$ replaced by $f-f(\xi_t)$ imply that 
\beg{align*}
|\nn P_{t}^Z f|(y)&=|\nn P_{t}^Z (f-f(\xi_t))|(y)\\
&\leq \ff {C} {\si_0 \sq{t\we 1}}(P^Z_{t}|f-f(\xi_{t})|^{p_2})^{\ff 1 {p_2}}(y)\\
&\leq \ff {Ce^{Ct}\|f\|_{V_1,\ph}} {\si_0^2}V_1(y)\left(1+|y|\right)^{\al_1\we \ff {2\al_1} {(1+\al_1)\al_2}}t^{-\ff 1 2}\ph\left(t^{\ff 1 2\we \ff {(1+\al_1)\al_2-2\al_1} {2(1+\al_1)\al_2}}\right).
\end{align*}
Similarly
\beg{align*}
|\nn^2P_{t}^Z f|(y)&\leq \ff {C\|f\|_{V_1,\ph}} {\si_0^2 (t\we 1)}(P^Z_{t}|f-f(\xi_{t})|^{p_2})^{\ff 1 {p_2}}(y)\\
&\leq \ff {Ce^{Ct}\|f\|_{V_1,\ph}} {\si_0^2}V_1(y)\left(1+|y|\right)^{\al_1\we \ff {2\al_1} {(1+\al_1)\al_2}}t^{-1}\ph\left(t^{\ff 1 2\we \ff {(1+\al_1)\al_2-2\al_1} {2(1+\al_1)\al_2}}\right).
\end{align*}

\end{proof}

\section{Generation and regularity of semigroups}

This section is devoted to study the perturbation of $P_t^{\mu_\infty}$ by $\bar A$ on $L^2(\mu_\infty)$.  In our concrete examples, see Example \ref{exa1} and Example \ref{exa2}, $\bar A$ can be extended to a bounded operator on $L^2(\mu_\infty)$ by using the integration by part formula and the additional regularity of the density of $\mu_\infty$ w.r.t. the Lebesgue measure.  To investigate the perturbation of $P_t^{\mu_\infty}$ by $\bar A$ under the minimal condition on $\mu_\infty$ (i.e. {\bf (H0)}), we introduce a slightly general framework, which allows us to treat $\bar A$ as an unbounded operator on $L^2(\mu_\infty)$ directly. Let $\{P_t\}_{t\geq 0}$ be a $C_0$-semigroup of contractions  on $L^2(\mu_\infty)$ with generator $L$, and let $S$ be a  linear operator in $L^2(\mu_\infty)$. In Subsection \ref{4.1},  a generation theorem for the perturbation of $P_t$ by $S$ is established.  Denote by $Q_t$ the semigroup generated by $L+S$. We investigate the regularity of $Q_t$ in Subsection \ref{4.2}.  By using theorems in the first two subsections, we study in Subsection \ref{4.3}  the $C_0$-semigroup generated by $L_{\mu_\infty}+\bar A$.

\subsection{A generation theorem for $C_0$-semigroups}\label{4.1}

  We assume that $P_t$ and $S$ satisfy following assumptions. 
\beg{description}
\item {\bf Assumption (B)}
\item [(B1)] Restricting on $W^{1,2}_{\mu_\infty}$, $\{P_t\}_{t\geq 0}$ is a $C_0$-semigroup on $W^{1,2}_{\mu_\infty}$, i.e. for any $t>0$ and $f\in W^{1,2}_{\mu_\infty}$, $P_tf\in W^{1,2}_{\mu_\infty}$ and  
\[\lim_{t\ra 0^+}\|P_tf-f\|_{W^{1,2}_{\mu_\infty}}=0.\]
There is $C_0>0$ such that
\beg{equation}\label{nnP-L2}
\|\nn P_tf\|_{L^2(\mu_\infty)}\leq \ff {C_0} {\sq{t\we 1}}\|f\|_{L^2(\mu_\infty)},~t>0.
\end{equation}

\item [(B2)] $S$ is a bounded operator from $W^{1,2}_{\mu_\infty}$ to $L^2(\mu_\infty)$, and $[0,+\infty)\ni t\mapsto P_tS$ is continuous under the operator norm from $W^{1,2}_{\mu_\infty}$ to $L^2(\mu_\infty)$. 
\end{description}

We denote by  $\sD(L)$ the domain of $L$ and by $\|\cdot\|_{W,2}$ the operator norm of bounded operators from $W^{1,2}_{\mu_\infty}$ to $L^2(\mu_\infty)$. Then {\bf (B2)} indicates that 
\beg{equation}\label{S-W120}
\lim_{t\ra 0^+}\|P_tS-S\|_{W,2}\equiv\lim_{t\ra 0^+}\sup_{\|f\|_{W^{1,2}_{\mu_\infty}}\leq 1}\|P_tSf-Sf\|_{L^2(\mu_\infty)}=0.
\end{equation}
\beg{rem}\label{rem-comp}
A sufficient condition for \eqref{S-W120} is that $S$ is a compact operator from $W^{1,2}_{\mu_\infty}$ to $L^2(\mu_\infty)$. Indeed,  if $S$ is a compact operator from $W^{1,2}_{\mu_\infty}$ to $L^2(\mu_\infty)$, then the set $\{Sf~|~\|f\|_{W^{1,2}_{\mu_\infty}}\leq 1\}$ is a precompact set in $L^2(\mu_\infty)$. Thus, for each $\ep>0$, there exist $\{f_n\}_{n=1}^m$ with $\|f_n\|_{W^{1,2}_{\mu_\infty}}\leq 1$ such that 
\[\{Sf~|~\|f\|_{W^{1,2}_{\mu_\infty}}\leq 1\}\subset\bigcup_{n=1}^m \{Sf~|~\|Sf-Sf_n\|_{L^2(\mu_\infty)}<\ep\}\equiv \bigcup_{n=1}^m \mathscr{S}_n.\]
This indicates  that
\beg{align*}
\varlimsup_{t\ra 0^+}\sup_{\|f\|_{W^{1,2}_{\mu_\infty}}\leq 1}\|P_tSf-Sf\|_{L^2(\mu_\infty)}&\leq \varlimsup_{t\ra 0^+}\max_{1\leq n\leq m}\|P_tSf_n-Sf_n\|_{L^2(\mu_\infty)}\\
&\quad\, +\varlimsup_{t\ra 0^+}\max_{1\leq n\leq m}\sup_{Sf\in\mathscr{S}_n}\|(P_t-I)(Sf-Sf_n)\|_{L^2(\mu_\infty)}\\
&\leq 2\ep,
\end{align*}
which implies \eqref{S-W120}.
\end{rem}
Under the assumption {\bf (B)}, we have the following generation theorem for the $C_0$-semigroup by $L+S$.
\beg{thm}\label{pert}
Assume that $P_t$ and $S$ satisfy {\bf (B)}. Then there is a unique $C_0$-semigroup $Q_t$ on $L^2(\mu_\infty)$ satisfying 
\beg{equation}\label{Qt-Pt-int}
Q_tf=P_t f+\int_0^t P_{s}SQ_{t-s}f\d s,~t\geq 0,~f\in L^2(\mu_\infty)
\end{equation}
where the integral converges in $W^{1,2}_{\mu_\infty}$, and there is $C_1>0$ which depends on $C_0$ in \eqref{nnP-L2} and $\|S\|_{W,2}$ such that
\beg{equation}\label{Q-W12}
\|Q_tf\|_{W^{1,2}_{\mu_\infty}}\leq \ff {C_1} {\sq t} e^{C_1t}\|f\|_{L^2(\mu_\infty)},~t>0,~f\in L^2(\mu_\infty).
\end{equation}
Moreover, \\
(1) for any $f\in L^2(\mu_\infty)$ and $t\geq 0$, 
\beg{equation}\label{exp-Qt}
\|Q_tf\|_{L^2(\mu_\infty)}\leq  \left(1+ 2C_1\|S\|_{W,2}\sq{t} e^{C_1t}\right)\|f\|_{L^2(\mu_\infty)},
\end{equation}
(2) $Q_t$ is a $C_0$-semigroup on $W^{1,2}_{\mu_\infty}$, \\
(3) the generator of $Q_t$ is $L+S$ and $\sD(L+S)=\sD(L)$.
\end{thm} 

\beg{proof}
Since $S$ is a bounded operator from $W^{1,2}_{\mu_\infty}$ to $L^2(\mu_\infty)$,  for any $\ep>0$, we have by $\|P_{\ep}\|_{L^2(\mu_\infty)}\leq 1$ and \eqref{nnP-L2} that 
\beg{align}\label{PS-W2}
\|P_{\ep}Sf\|_{L^2(\mu_\infty)}&\leq \|Sf\|_{L^2(\mu_\infty)}\leq \|S\|_{W,2}\|f\|_{W^{1,2}_{\mu_\infty}},\\
\|\nn P_{\ep}Sf\|_{L^2(\mu_\infty)}&\leq \ff {C_0} {\sq{\ep\we 1}}\|Sf\|_{L^2(\mu_\infty)}\leq \ff {C_0} {\sq{\ep\we 1}}\|S\|_{W,2}\|f\|_{W^{1,2}_{\mu_\infty}},\nonumber
\end{align}
which yield that $P_{\ep}S$ is a bounded operator on $W^{1,2}_{\mu_\infty}$. Denote by $L^{W}$ the generator of the semigroup $P_t$ restricted on $W^{1,2}_{\mu_\infty}$. It follows from \cite[Proposition 3.1.2 and (1.10)]{Pazy} that $L^W+P_{\ep}S$ generates a $C_0$-semigroup on $W^{1,2}_{\mu_\infty}$, which is denoted by $Q_t^{\ep}$, and
\beg{equation}\label{Qep=Pt}
Q^{\ep}_tf=P_t f+\int_0^t P_s(P_{\ep}S) Q_{t-s}^{\ep}f\d s,~f\in W^{1,2}_{\mu_\infty}.
\end{equation}
Then 
\beg{equation}\label{nnQep0}
\beg{split}
\|Q^\ep_tf\|_{L^2(\mu_\infty)} &\leq \|P_tf\|_{L^2(\mu_\infty)}+\int_0^t \|P_{s}(P_{\ep}S) Q^{\ep}_{t-s}f\|_{L^2(\mu_\infty)} \d s\\
&\leq  \|f\|_{L^2(\mu_\infty)} +\int_0^t \| S\|_{W,2} \| Q^{\ep}_{t-s}f\|_{W^{1,2}_{\mu_\infty}} \d s,
\end{split}
\end{equation} 
and by using \eqref{nnP-L2},
\beg{equation}\label{nnQep}
\beg{split}
\|\nn Q^\ep_tf\|_{L^2(\mu_\infty)} &\leq \|\nn P_tf\|_{L^2(\mu_\infty)}+\int_0^t \|\nn P_{s}(P_{\ep}S) Q^{\ep}_{t-s}f\|_{L^2(\mu_\infty)} \d s\\
&\leq \ff {C_0} {\sq {t\we 1}}\|f\|_{L^2(\mu_\infty)} +\int_0^t \ff {C_0} {\sq{s\we 1}}\|P_{\ep}S\|_{W,2}\cdot\| Q^{\ep}_{t-s}f\|_{W^{1,2}_{\mu_\infty}} \d s\\
&=\ff {C_0} {\sq{t\we 1}}\|f\|_{L^2(\mu_\infty)} +\int_0^t \ff {C_0\|S\|_{W,2}} {\sq{(t-s)\we 1}}\| Q^{\ep}_{s}f\|_{W^{1,2}_{\mu_\infty}} \d s.
\end{split}
\end{equation} 
By setting 
\[\ga_1(t)= \ff {1+C_0} {\sq{t\we 1}}\|f\|_{L^2(\mu_\infty)},\qquad \ga_2(t)= (1+C_0)\|S\|_{W,2},\] 
it follows from \eqref{nnQep0} and \eqref{nnQep} that 
\[\|Q^{\ep}_tf\|_{W^{1,2}_{\mu_\infty}}\leq \ga_1(t)\|f\|_{L^2(\mu_\infty)}+\int_0^t\ff {\ga_2(s)} {\sq{(t-s)\we 1}}\| Q^{\ep}_{s}f\|_{W^{1,2}_{\mu_\infty}}\d s.\]
Then, by using  Lemma \ref{lem-Gron},  we have that 
\beg{align*}
\Ga_1(t)&\leq\int_0^t\ff {(1+C_0)\|f\|_{L^2(\mu_\infty)}} {\sq{((t-s)\we 1)(s\we 1)}}\d s= (1+C_0)\varphi_0(t)\|f\|_{L^2(\mu_\infty)},\\
\Ga_2(t)&\leq (1+ C_0)\varphi_1(t)\|f\|_{L^2(\mu_\infty)}+(1+C_0)^2\|S\|_{W,2}\|f\|_{L^2(\mu_\infty)}\varphi_2(t),
\end{align*}
and
\beg{align*}
&\| Q^\ep_tf\|_{W^{1,2}_{\mu_\infty}} \leq (1+C_0)\left(\ff {1} {\sq{t\we 1}}+(1+C_0)\|S\|_{W,2} \varphi_0(t)\right) \|f\|_{L^2(\mu_\infty)}\\
&\quad\,+(1+C_0)^3\|S\|_{W,2}^2\varphi_0(t)\left[\varphi_1(t)+(1+C_0)\|S\|_{W,2}\varphi_2(t) \right]\e^{(1+C_0)^2\|S\|_{W,2}^2\varphi_2(t)}\|f\|_{L^2(\mu_\infty)},
\end{align*}
where $\varphi_0(t)$ is defined by \eqref{app-phi0} and
\beg{align*}
\varphi_1(t)&=2\sq{t}\1_{[0\leq t\leq 1]}+(t+1)\1_{[t>1]},\\
\varphi_2(t)&=\int_0^t\varphi_0(s)\d s=B\left(\ff 1 2,\ff 1 2\right)t+\ff 4 3\sq{t-1}\1_{[1\leq t\leq 2]}+\left(\ff {t^2} 2-\ff 2 3\right)\1_{[t>2]},~t\geq 0.
\end{align*}
By using the following fundamental inequality
\beg{equation}\label{in-funda}
2\sq r\leq r+1\leq e^r,~r\geq 0,
\end{equation}
there is a  constant $C_{1}>0$ which is independent of $\ep$ and depends only on $C_0$, $\|S\|_{W,2}$ such that 
\beg{equation}\label{Qep-e}
\|Q^\ep_tf\|_{W^{1,2}_{\mu_\infty}}\leq \ff {C_1e^{C_1t}} {\sq t}\|f\|_{L^2(\mu_\infty)},~f\in W^{1,2}_{\mu_\infty},~t>0.
\end{equation}
By the approximation argument, we see that for all $t>0$, $Q^{\ep}_t$ is a bounded operator from $L^2(\mu_\infty)$ to $W^{1,2}_{\mu_\infty}$ and satisfies \eqref{Qep-e} for any $f\in L^2(\mu_\infty)$. Consequently,  
\beg{align}\label{int-Qep-e}
\int_0^t\| Q^\ep_sf\|_{W^{1,2}_{\mu_\infty}}\d s&\leq 2C_1\sq{t}e^{C_1t}\|f\|_{L^2(\mu_\infty)},\\
\int_0^t\ff{\|Q^\ep_s f\|_{W^{1,2}_{\mu_\infty}}} {\sq{(t-s)\we 1}}\d s&\leq C_1e^{C_1t}\|f\|_{L^2(\mu_\infty)}\int_0^t \ff {\d s} {\sq{s[(t-s)\we 1]}}\nonumber\\
&= C_1e^{C_1t} \left(B\left(\ff 1 2,\ff 1 2\right)+2\sq{(t-1)^+}\right)\|f\|_{L^2(\mu_\infty)},\label{int-Qe-tse}
\end{align} 
and \eqref{Qep=Pt} holds for any $f\in L^2(\mu_\infty)$. We derive from \eqref{int-Qep-e} and \eqref{nnQep0} that 
\beg{equation}\label{exp-Qt-ep0}
\beg{split}
\|Q^\ep_tf\|_{L^2(\mu_\infty)}&\leq \|f\|_{L^2(\mu_\infty)}+2C_1\|S\|_{W,2}\sq t e^{C_1t} \|f\|_{L^2(\mu_\infty)}\\
&= \left(1+ 2C_1\|S\|_{W,2}\sq{t} e^{C_1t}\right)\|f\|_{L^2(\mu_\infty)},~t\geq 0.
\end{split}
\end{equation}

Next, we take $\ep\downarrow 0$.  For any $\ep_1,\ep_2>0$ and $f\in L^2(\mu_\infty)$, we have that
\beg{equation}\label{Qep-ep}
\left(Q^{\ep_1}_t-Q^{\ep_2}_t\right)f=\int_0^t P_s\left(\left(P_{\ep_1}-P_{\ep_2}\right)S\right)Q^{\ep_1}_{t-s}f\d s+\int_0^t P_s\left(P_{\ep_2} S\right)\left(Q^{\ep_1}_{t-s}-Q^{\ep_2}_{t-s}\right)f\d s.
\end{equation} 
This, together with \eqref{nnP-L2} and \eqref{int-Qe-tse}, implies that 
\beg{align*}
\|\left(Q^{\ep_1}_t-Q^{\ep_2}_t\right)f\|_{W^{1,2}_{\mu_\infty}}&\leq \int_0^t \| P_s\left(\left(P_{\ep_1}-P_{\ep_2}\right)S\right)Q^{\ep_1}_{t-s}f\|_{W^{1,2}_{\mu_\infty}}\d s\\
&\quad\, +\int_0^t\| P_s\left(P_{\ep_2}S\right)\left(Q^{\ep_1}_{t-s}-Q^{\ep_2}_{t-s}\right)f\|_{W^{1,2}_{\mu_\infty}}\d s\\
&\leq \|(P_{\ep_1}-P_{\ep_2})S\|_{W,2}\int_0^t \ff {1+C_0} {\sq{s\we 1}}\| Q^{\ep_1}_{t-s}f\|_{W^{1,2}_{\mu_\infty}}\d s\\
&\quad\, +\|S\|_{W,2}\int_0^t\ff {1+C_0} {\sq{s\we 1}}\| \left(Q^{\ep_1}_{t-s}-Q^{\ep_2}_{t-s}\right)f\|_{W^{1,2}_{\mu_\infty}}\d s\\
&\leq (1+C_0)C_2e^{C_2t}\|(P_{\ep_1}-P_{\ep_2})S\|_{W,2}\|f\|_{L^2(\mu_\infty)}\\
&\quad\, +(1+C_0)\|S\|_{W,2}\int_0^t\ff {\| \left(Q^{\ep_1}_{s}-Q^{\ep_2}_{s}\right)f\|_{W^{1,2}_{\mu_\infty}}} {\sq{(t-s)\we 1}}\d s,
\end{align*}
where $C_2>0$ is a constant independent of $\ep_1,\ep_2$. By setting 
\[\ga_1(t)=(1+C_0)C_2e^{C_2t}\|(P_{\ep_1}-P_{\ep_2})S\|_{W,2}\|f\|_{L^2(\mu_\infty)},\qquad \ga_2(t)=(1+C_0)\|S\|_{W,2}\] 
in Lemma \ref{lem-Gron}, we have that 
\beg{align*}
\Ga_1(t)&\leq (1+C_0)C_2e^{C_2t}\varphi_1(t)\|(P_{\ep_1}-P_{\ep_2})S\|_{W,2}\|f\|_{L^2(\mu_\infty)},\\
\Ga_2(t)&\leq (1+C_0)\left(e^{C_2t}-1\right) \|(P_{\ep_1}-P_{\ep_2})S\|_{W,2}\|f\|_{L^2(\mu_\infty)}\\
&\quad\,+(1+C_0)^2C_2\|S\|_{W,2}\|(P_{\ep_1}-P_{\ep_2})S\|_{W,2}\|f\|_{L^2(\mu_\infty)}\int_0^te^{C_2s}\varphi_1(s)\d s\\
&\leq (1+C_0)^2\left(e^{C_2t}-1\right)\|(P_{\ep_1}-P_{\ep_2})S\|_{W,2}\|f\|_{L^2(\mu_\infty)}\left(1+\|S\|_{W,2}\varphi_1(t)\right).
\end{align*}
Then there is a constant $C_3>0$ which is also independent of $\ep_1,\ep_2$ such that
\beg{align*}
\|\left(Q^{\ep_1}_t-Q^{\ep_2}_t\right)f\|_{W^{1,2}_{\mu_\infty}}&\leq C_3e^{C_3t}\|(P_{\ep_1}-P_{\ep_2})S\|_{W,2}\|f\|_{L^2(\mu_\infty)}.
\end{align*}
Combining this with \eqref{S-W120}, we have proven that $Q^{\ep}_tf$ is a Cauchy sequence in $W^{1,2}_{\mu_\infty}$. Thus, there is $Q_tf\in W^{1,2}_{\mu_\infty}$ such that 
\beg{equation}\label{limQQ}
\lim_{\ep\ra 0^+}\|Q^{\ep}_tf-Q_tf\|_{W^{1,2}_{\mu_\infty}}=0,~f\in L^2(\mu_\infty),~t>0.
\end{equation}
Moreover, it is clear that $\{Q_t\}_{t>0}$ is a semigroup from $L^2(\mu_\infty)$ to $W^{1,2}_{\mu_\infty}$ since so is $\{Q^{\ep}_t\}_{t>0}$.  Following from \eqref{limQQ}, \eqref{Qep-e} and \eqref{exp-Qt-ep0}, $Q_tf$ satisfies  \eqref{Q-W12} and \eqref{exp-Qt}.  Taking $\ep\ra 0^+$ on both side of \eqref{Qep=Pt} in $W^{1,2}_{\mu_\infty}$, we find that $Q_tf$ satisfies \eqref{Qt-Pt-int} and the integral converges in $W^{1,2}_{\mu_\infty}$.  Repeating the argument of estimating $\left(Q^{\ep_1}_t-Q^{\ep_2}_t\right)f$, we can prove that there is a unique $Q_tf$ satisfying \eqref{Qt-Pt-int} and \eqref{Q-W12}. 

By \eqref{Qt-Pt-int} and \eqref{Q-W12}, we have that
\beg{align*}
\|Q_tf-f\|_{L^2(\mu_\infty)}&\leq \|P_tf-f\|_{L^2(\mu_\infty)}+\|S\|_{W,2}\int_0^t \|  Q_{t-s}f\|_{W^{1,2}_{\mu_\infty}}\d s\\
&\leq \|P_tf-f\|_{L^2(\mu_\infty)}+C_1\|S\|_{W,2}\int_0^t \ff {e^{C_1(t-s)}} {\sq{(t-s)}}\d s\|f\|_{L^2(\mu_\infty)}\\
&\leq \|P_tf-f\|_{L^2(\mu_\infty)}+2C_1\|S\|_{W,2}\sq t e^{C_1t}\|f\|_{L^2(\mu_\infty)}.
\end{align*}
Thus
\[\lim_{t\ra 0^+}\|Q_tf-f\|_{L^2(\mu_\infty)}=\lim_{t\ra 0^+}\|P_tf-f\|_{L^2(\mu_\infty)}=0.\]
Hence, $\{Q_t\}_{t\geq 0}$ is a $C_0$-semigroup on $L^2(\mu_\infty)$.

Next, we prove (2). We first investigate $Q_t$ on $W^{1,2}_{\mu_\infty}$. Since $P_t$ is a $C_0$-semigroup on $W^{1,2}_{\mu_\infty}$, according to \cite[Theorem 1.2.2]{Pazy}, there are $C\geq 1$ and $w\geq 0$ such that
\[\| P_tf\|_{W^{1,2}_{\mu_\infty}}\leq Ce^{wt}\|f\|_{W^{1,2}_{\mu_\infty}}.\]
Combining this with \eqref{nnQep0} and the first inequality in \eqref{nnQep},   we have for any $f\in W^{1,2}_{\mu_\infty}$ that
\beg{align*}
\| Q^{\ep} _tf\|_{W^{1,2}_{\mu_\infty}}&\leq \| P_tf\|_{W^{1,2}_{\mu_\infty}}+\int_0^t \ff {(1+C_0)\|S\|_{W,2}} {\sq{s\we 1}}\| Q^{\ep}_{t-s}f\|_{W^{1,2}_{\mu_\infty}} \d s\\
&\leq Ce^{wt}\|f\|_{W^{1,2}_{\mu_\infty}}+(1+C_0)\|S\|_{W,2}\int_0^t \ff {\| Q^{\ep}_{t-s}f\|_{W^{1,2}_{\mu_\infty}} } {\sq{s\we 1}}\d s.
\end{align*}
This, together with Lemma \ref{lem-Gron}, implies that there is $C_4>0$ which is independent of $\ep$ so that
\[\|Q^{\ep} _tf\|_{W^{1,2}_{\mu_\infty}}\leq C_4e^{C_4t}\| f\|_{W^{1,2}_{\mu_\infty}}.\]
Hence, due to  \eqref{limQQ}, we arrive at
\[\|Q_tf\|_{W^{1,2}_{\mu_\infty}}\leq C_4e^{C_4t}\| f\|_{W^{1,2}_{\mu_\infty}}.\]
Moreover, 
\beg{align*}
\|Q_tf-f\|_{W^{1,2}_{\mu_\infty}}&\leq \|P_tf-f\|_{W^{1,2}_{\mu_\infty}}+\int_0^t\|P_sS Q_{t-s}f\|_{W^{1,2}_{\mu_\infty}}\d s\\
&\leq \|P_tf-f\|_{W^{1,2}_{\mu_\infty}}+\int_0^t\left(1+\ff {C_0} {\sq{s\we 1}}\right)\|S\|_{W,2}\| Q_{t-s}f\|_{W^{1,2}_{\mu_\infty}}\d s\\
&\leq \|P_tf-f\|_{W^{1,2}_{\mu_\infty}}+(1+C_0)C_4\|S\|_{W,2}\int_0^t \ff {e^{C_4(t-s)}} {\sq{s\we 1}}\d s\|f\|_{W^{1,2}_{\mu_\infty}}\\
&\leq \|P_tf-f\|_{W^{1,2}_{\mu_\infty}}+2(1+C_0)C_4\|S\|_{W,2}e^{C_4t}\varphi_1(t)\|f\|_{W^{1,2}_{\mu_\infty}}.
\end{align*}
This implies that $Q_t$ is also a $C_0$-semigroup on $W^{1,2}_{\mu_\infty}$.

Finally, we investigate the generator of $Q_t$. We first prove that  $\sD(L)\subset W^{1,2}_{\mu_\infty}$. Since $P_t$ is contractive, for any $\la>0$, there is 
\beg{equation}\label{resol}
(\la-L)^{-1}f=\int_0^{+\infty} e^{-\la t} P_tf\d t,~f\in L^2(\mu_\infty).
\end{equation}
Since $\{P_t\}_{t\geq 0}$ is a $C_0$-semigroup on $W^{1,2}_{\mu_\infty}$, there is $C,w>0$ such that
\[\|P_tf\|_{W^{1,2}_{\mu_\infty}}\leq Ce^{wt}\|f\|_{W^{1,2}_{\mu_\infty}},~f\in W^{1,2}_{\mu_\infty}.\]
Then the right hand side of \eqref{resol} converges in $W^{1,2}_{\mu_\infty}$ for any $f\in W^{1,2}_{\mu_\infty}$ and $\la>w$. This, together with \eqref{nnP-L2}, implies that
\beg{align*}
\|(\la-L)^{-1}f\|_{W^{1,2}_{\mu_\infty}}&\leq \int_0^{+\infty}e^{-\la t} \|P_tf\|_{W^{1,2}_{\mu_\infty}}\d t\\
& \leq \int_0^{+\infty} \left(\ff {C_0} {\sq{t\we 1}} +1\right)e^{-\la t}\|f\|_{L^2(\mu_\infty)}\d t\\
&\leq (C_0+1)\int_0^{+\infty} \ff {e^{-\la t}} {\sq{t\we 1}}\d t \|f\|_{L^2(\mu_\infty)},~f\in W^{1,2}_{\mu_\infty}. 
\end{align*}
By using the approximation argument, we see that there is $C>0$ such that
\[\|(\la-L)^{-1}f\|_{W^{1,2}_{\mu_\infty}}\leq C\|f\|_{L^2(\mu_\infty)},~f\in L^2(\mu_\infty).\]
Thus $\sD(L)\subset W^{1,2}_{\mu_\infty}$.

For $f\in \sD(L)$, there is  
\beg{align*}
\ff {Q_tf-f} t &= \ff {P_tf-f} {t}+\ff 1 t\int_0^t P_{t-s}SQ_sf\d s\\
&=\ff {P_tf-f} {t}+\ff 1 t\int_0^t \left(P_{t-s}-I\right)S Q_sf\d s+\ff 1 t\int_0^t SQ_sf\d s\\
&=: I_1+I_2+I_3.
\end{align*}
For $I_2$, we have that
\beg{align*}
\left\|\ff 1 t\int_0^t \left(P_{t-s}-I\right)S Q_sf\d s\right\|_{L^2(\mu_\infty)}&\leq \ff 1 t \int_0^t \|\left(P_{t-s}-I\right)S\|_{W,2}\|Q_sf\|_{W^{1,2}_{\mu_\infty}}\d s\\
&\leq \ C_4e^{C_4t}\|f\|_{W^{1,2}_{\mu_\infty}}\ff 1 t\int_0^t \|\left(P_{s}-I\right)S\|_{W,2}\d s,
\end{align*}
where $f\in \sD(L)\subset W^{1,2}_{\mu_\infty}$ has been used in the last inequality. Combining this with \eqref{S-W120}, there is $\displaystyle\lim_{t\ra 0^+}\| I_2\|_{L^2(\mu_\infty)}=0$.\\
Since $\sD(L)\subset W^{1,2}_{\mu_\infty}$, $S$ is bounded from $W^{1,2}_{\mu_\infty}$ to $L^2(\mu_\infty)$ and $Q_\cdot f$ is continuous in $W^{1,2}_{\mu_\infty}$, we find that
\[\lim_{t\ra 0^+}\ff 1 t\int_0^t S Q_sf\d s=S f,~\text{in}~L^2(\mu_\infty).\]
Since $f\in\sD(L)$, we have that $\lim_{t\ra 0^+}I_1=Lf$ in $L^2(\mu_\infty)$.\\
Hence,  there is 
\[\lim_{t\ra 0^+}\left\|\ff {Q_tf-f} t -Lf-Sf\right\|_{L^2(\mu_\infty)}=0.\]
This yields that the generator of $Q_t$ on $L^2(\mu_\infty)$ is $L+S$ and $\sD(L)\subset \sD(L+S)$. We can prove similarly that $\sD(L+S)\subset \sD(L)$ since 
\[\ff {P_tf-f} t=\ff {Q_tf-f} t-\ff 1 t\int_0^t \left(P_{t-s}-I\right)S Q_sf\d s-\ff 1 t\int_0^t SQ_sf\d s,~f\in\sD(L+S).\]

\end{proof}

\subsection{Regularity of $Q_t$}\label{4.2}

In Theorem \ref{pert}, we have proven that $Q_tf(\cdot)$ is  differentiable in the sense of Sobolev for any $f\in L^2(\mu_\infty)$ and $t>0$. Next, we investigate the  regularity of $Q_tf$ in the strong sense when $f$ is regular. Stronger assumptions are imposed on $P_t$ and $S$.  For simplicity, we denote by $G_f(t,z)$ the second term in \eqref{Qt-Pt-int}  and
\[G_{f}^{(1)}(t,z)=\int_0^t\nn P_sSQ_{t-s}f(z)\d s,\qquad G_{f}^{(2)}(t,z)=\int_0^t\nn^2 P_sSQ_{t-s}f(z)\d s.\]
\beg{thm}\label{thm-regu}
Let $P_t$ and $S$ satisfy {\bf (B)}. Suppose that for any $t>0$ and $f\in W^{1,2}_{\mu_\infty}$,  $P_tSf\in C^1$ and there is a locally bounded measurable function $1\leq V_1\in L^2(\mu_\infty)$ such that
\beg{equation}\label{nnPS-p}
|\nn P_tSf|(z)\leq \ff {Ce^{Ct}} {\sq{t}}V_1(z)\|f\|_{W^{1,2}_{\mu_\infty}}.
\end{equation}
Then for any $f\in L^2(\mu_\infty)$, $G_f(\cdot,\cdot)\in C^{0,1}((0,+\infty)\times\R^d)$ and $\nn G_f=G_f^{(1)}$ with
\beg{equation}\label{nnG1}
\left|G_{f}^{(1)}(t,z)\right|\leq Ce^{Ct}V_1(z)\|f\|_{L^2(\mu_\infty)},~t> 0,~z\in\R^d,~f\in L^2(\mu_\infty).
\end{equation}
Moreover,   if $f\in W^{1,2}_{\mu_\infty}$, then 
\[
\left|G_{f}^{(1)}(t,z)\right|\leq C\sq t e^{Ct}V_1(z)\|f\|_{W^{1,2}_{\mu_\infty}},~t\geq 0,~z\in\R^d,~f\in W^{1,2}_{\mu_\infty},
\]
and there is $G_f(\cdot,\cdot)\in C^{0,1}([0,+\infty)\times\R^d)$. Consequently, for  $f\in C^1_b$, if $P_\cdot f(\cdot)\in C^{0,1}([0,+\infty)\times\R^d)$, then $Q_{\cdot}f(\cdot)\in C^{0,1}([0,+\infty)\times\R^d)$. 

Assume in addition that for any $t>0$ and $f\in W^{1,2}_{\mu_\infty}$, $P_tSf\in C^2$ and there exist a locally bounded function $1\leq V_2\in L^2(\mu_\infty)$ and a positive function $\th\in L^1_{loc}([0,+\infty))$ such that  
\beg{align}\label{int-ths}
\Th(t)&:=\int_0^t \ff {\th(s)} {\sq{t-s}}\d s<+\infty,~t>0,\\
|\nn^2 P_tSf|(z)&\leq  Ce^{Ct}\th(t) V_2(z)\|f\|_{W^{1,2}_{\mu_\infty}},~z\in\R^d,t> 0.\label{nn2PS-p}
\end{align}
Then  for any $f\in L^2(\mu_\infty)$, $G_f(\cdot,\cdot)\in C^{0,2}((0,+\infty)\times\R^d)$ and $\nn^2 G_f=G_f^{(2)}$ with
\beg{equation}\label{nnG2}
\left|G_{f}^{(2)}(t,z)\right|\leq C\Th(t)e^{Ct}V_2(z)\|f\|_{L^2(\mu_\infty)},~t> 0,~z\in\R^d,~f\in L^2(\mu_\infty).
\end{equation}
Moreover, if $f\in W^{1,2}_{\mu_\infty}$, then there is $G_f(\cdot,\cdot)\in C^{0,2}([0,+\infty)\times\R^d)$ and
\beg{align}\label{Gf2W12}
\left|G_{f}^{(2)}(t,z)\right|\leq C\|\th\|_{L^1[0,t]}e^{Ct}V_2(z)\|f\|_{W^{1,2}_{\mu_\infty}},~t\geq 0,~z\in\R^d,~f\in L^2(\mu_\infty).
\end{align}
Consequently, for $f\in C^2_b$, if $P_\cdot f(\cdot)\in C^{0,2}([0,+\infty)\times\R^d)$, then $Q_{\cdot}f(\cdot)\in C^{0,2}([0,+\infty)\times\R^d)$.

\end{thm}

\beg{proof}
(1) For each $f\in L^2(\mu_\infty)$ and $s\in (0,t)$,  it follows from \eqref{Q-W12} and \eqref{nnPS-p} that   there is $C>0$ such that
\beg{equation}\label{PSQf-grow}
\beg{split}
|\nn P_sSQ_{t-s}f(z)|&\leq \ff {Ce^{Ct}} {\sq s}\|Q_{t-s}f\|_{W^{1,2}_{\mu_\infty}}V_1(z)\\
&\leq \ff {CC_1e^{Cs}e^{C_1(t-s)}} {\sq {s(t-s)}}\|f\|_{L^2(\mu_\infty)}V_1(z).
\end{split}
\end{equation}
Then $\nn P_sSQ_{t-s}f\in L^2(\mu_\infty)$, there is $C>0$ such that \eqref{nnG1} holds, and for any  $m\geq 1$
\beg{equation}\label{nnintP}
\beg{split}
\int_0^t \sup_{|z|\leq m}|\nn  P_{s}S Q_{t-s}f(z)|\d s&\leq \int_0^t \ff {Ce^{Ct}} {\sq {s(t-s)}}\|f\|_{L^2(\mu_\infty)}\left(\sup_{|z|\leq m}V_1(z)\right)\d s\\
&\leq C\left(\sup_{|z|\leq m}V_1(z)\right)\|f\|_{L^2(\mu_\infty)}B\left(\ff 1 2,\ff 1 2\right)e^{Ct}.
\end{split}
\end{equation}
Combining this with $P_{s}SQ_{t-s}f(\cdot)\in C^1$, the dominated convergence theorem yields that
\beg{equation}\label{nnPt-Ptnn}
\beg{split}
\nn_v G_f(t,z)&=\lim_{\ep\ra 0^+}\int_0^t  \ff 1 {\ep} \left(P_{t-s}S Q_sf(z+\ep v)-P_{t-s}S Q_sf(z)\right)\d s\\
&=\int_0^t \nn_v P_{t-s}S Q_sf(z)\d s.
\end{split}
\end{equation}
Using  $P_{s}S Q_{t-s}f(\cdot)\in C^1$ again, we can derive from \eqref{nnPt-Ptnn}, \eqref{nnintP} and the dominated convergence theorem  that  $G_f^{(1)}(t,\cdot)\in C^0$ and $\nn G_f=G^{(1)}_f$. \\
For $f\in W^{1,2}_{\mu_\infty}$, since (2) of Theorem \ref{pert}, there is $C>0$ such that
\beg{align*}
\int_0^t \sup_{|z|\leq m}|\nn  P_{s}S Q_{t-s}f(z)|\d s&\leq \int_0^t \ff {Ce^{Ct}} {\sq {s}}\|Q_{t-s}f\|_{W^{1,2}_{\mu_\infty}}\left(\sup_{|z|\leq m}V_1(z)\right)\d s\\
&\leq C\left(\sup_{|z|\leq m}V_1(z)\right)\|f\|_{L^2(\mu_\infty)}\sq t e^{Ct}.
\end{align*}
Hence, we set 
\[G_f^{(1)}(0,z)=0,~f\in W^{1,2}_{\mu_\infty},~z\in\R^d.\]
For any $m\geq 1$ and $t_1>t_2\geq 0$, it follows from \eqref{nnPS-p} and \eqref{nnintP} that there are positive constants $C_m$ which depending on $m$ and $C$ independent of $m$   such that
\beg{align*}
&\sup_{|z|\leq m}\left|G_f^{(1)}(t_1,z)-G_f^{(1)}(t_2,z)\right|\\
&\quad\,\leq\int_0^{t_2} \sup_{|z|\leq m}|\nn P_sS Q_{t_2-s}(Q_{t_1-t_2}-I)f(z)|\d s+\int_{t_2}^{t_1} \sup_{|z|\leq m}|\nn P_{s}S Q_{t_1-s}f(z)|\d s\\
&\quad\, \leq C_m\left(\int_0^{t_2} \ff {e^{Ct_2}\|(Q_{t_1-t_2}-I)f\|_{L^2(\mu_\infty)}} {\sq{s(t_2-s)}}\d s+ \int_{t_2}^{t_1}\ff {e^{Ct_1}\|Q_{t_1-s}f\|_{W^{1,2}_{\mu_\infty}}} {\sq{s}}\d s\right)\\
&\quad\, =: I_1+I_2.
\end{align*}
Since $Q_t$ is a $C_0$-semigroup on $L^2(\mu_\infty)$, there is 
\[\lim_{|t_1-t_2|\ra 0^+}\|(Q_{t_1-t_2}-I)f\|_{L^2(\mu_\infty)}=0,\]
which yields $\displaystyle\lim_{|t_1-t_2|\ra 0^+}I_1=0$. For $I_2$, if there is $c>0$ such that $t_2\geq c$, then
\beg{align*}
\varlimsup_{|t_1-t_2|\ra 0^+\atop t_2\geq c>0}\int_{t_2}^{t_1}\ff {\|Q_{t_1-s}f\|_{W^{1,2}_{\mu_\infty}}} {\sq{s\we 1}}\d s&\leq \varlimsup_{|t_1-t_2|\ra 0^+\atop t_2\geq c>0}\left(\int_{t_2}^{t_1}\ff {Ce^{C(t_1-s)}} {\sq{(s\we 1) (t_1-s)}}\d s\right)\|f\|_{L^2(\mu_\infty)}\\
&\leq  \varlimsup_{|t_1-t_2|\ra 0^+ \atop t_2\geq c>0}\left(\int_{t_2}^{t_1}\ff {Ce^{C(t_1-s)}} {\sq{t_1-s}}\d s\right)\ff {\|f\|_{L^2(\mu_\infty)}} {\sq{c\we 1}}\\
&=0;
\end{align*}
if $t_2=0$, then for $f\in W^{1,2}_{\mu_\infty}$, we find that
\[
\varlimsup_{|t_1-t_2|\ra 0^+}\int_{t_2}^{t_1}\ff {\|Q_{t_1-s}f\|_{W^{1,2}_{\mu_\infty}}} {\sq{s\we 1}}\d s\leq \varlimsup_{t_1\ra 0^+}\left(\int_{0}^{t_1}\ff {Ce^{C(t_1-s)}} {\sq{s\we 1}}\d s\right)\|f\|_{W^{1,2}_{\mu_\infty}}=0.
\]
Hence,  $G_f^{(1)}(\cdot,z)$ is continuous  locally uniformly w.r.t. $z$  on $(0,+\infty)$ when $f\in L^2(\mu_\infty)$ and on $[0,+\infty)$ when $f\in W^{1,2}_{\mu_\infty}$.  This implies that  
\beg{equation*}
\lim_{(t_1,z_1)\ra (t_2,z_2)}\left|G_f^{(1)}(t_1,z_1)-G_f^{(1)}(t_2,z_2)\right|=0,~\left\{\beg{array}{l}
t_2>0,~f\in L^2(\mu_\infty),\\
t_2=0,~f\in W^{1,2}_{\mu_\infty}.
\end{array}\right.
\end{equation*}
For $f\in C^1_b\subset W^{1,2}_{\mu_\infty}$, since \eqref{Qt-Pt-int} and $P_\cdot f(\cdot) \in C^{0,1}([0,+\infty)\times\R^d)$, we have that $Q_{\cdot}f(\cdot)\in C^{0,1}([0,+\infty)\times\R^d)$. 

(2)  It follows from \eqref{Q-W12}, \eqref{int-ths} and \eqref{nn2PS-p} that $\{\nn^2P_{s}S Q_{t-s}f\}_{0<s<t}\subset L^2(\mu_\infty)$, and for any $m\geq 1$
\beg{equation}\label{int-nn2PAQ}
\beg{split}
\int_0^t \sup_{|z|\leq m}|\nn^2P_{s}S Q_{t-s}f(z)|\d s&\leq \int_0^t C  e^{C s} \th(s)\| Q_{t-s}f\|_{W^{1,2}_{\mu_\infty}}\left(\sup_{|z|\leq m} V_2(z)\right)\d s\\
&\leq CC_1\left(\sup_{|z|\leq m} V_2(z)\right)\|f\|_{L^2(\mu_\infty)}\int_0^t  \ff {e^{Cs}e^{C_1(t-s)} \th(s)} {\sq{t-s}}\d s\\
&\leq CC_1\left(\sup_{|z|\leq m} V_2(z)\right)e^{(C\vee C_1)t}\Th(t)\|f\|_{L^2(\mu_\infty)}.
\end{split}
\end{equation}
We can prove \eqref{nnG2} similarly. Since for $t>s>0$, $\nn^2P_{s}SQ_{t-s}f(\cdot)\in C^0$. This, together with \eqref{int-nn2PAQ}, \eqref{int-ths} and  the dominated convergence theorem implies $G_f^{(2)}(t,\cdot)\in C^0$  and $\nn^2G_f=G_f^{(2)}$.  \\
For $f\in W^{1,2}_{\mu_\infty}$, since (2) of Theorem \ref{pert}, there is $C>0$ such that
\beg{align*}
\int_0^t \sup_{|z|\leq m}|\nn^2P_{s}S Q_{t-s}f(z)|\d s&\leq \int_0^t C  e^{C s} \th(s)\| Q_{t-s}f\|_{W^{1,2}_{\mu_\infty}}\left(\sup_{|z|\leq m} V_2(z)\right)\d s\\
&\leq CC_1\left(\sup_{|z|\leq m} V_2(z)\right)\|f\|_{W^{1,2}_{\mu_\infty}}\int_0^t   e^{Cs}e^{C(t-s)} \th(s) \d s\\
&\leq CC_1\left(\sup_{|z|\leq m} V_2(z)\right)e^{Ct}\|\th\|_{L^1[0,t]}\|f\|_{W^{1,2}_{\mu_\infty}}.
\end{align*}
We can prove \eqref{Gf2W12} similarly. Hence, we can set 
\[G_f^{(2)}(0,z)=0,~f\in W^{1,2}_{\mu_\infty},~z\in\R^d.\]
For $m\geq 1$ and $t_1>t_2\geq 0$, it follows from \eqref{nn2PS-p} and \eqref{int-nn2PAQ} that there exist positive constants $C$ and $C_m$ such that
\beg{align*}
&\sup_{|z|\leq m}\left|G_f^{(2)}(t_1,z)-G_f^{(2)}(t_2,z)\right|\\
&\quad\,\leq\int_0^{t_2} \sup_{|z|\leq m}|\nn^2 P_sS Q_{t_2-s}(Q_{t_1-t_2}-I)f(z)|\d s+\int_{t_2}^{t_1} \sup_{|z|\leq m}|\nn^2 P_{s}S Q_{t_1-s}f(z)|\d s\\
&\quad\,\leq C_m\left(e^{Ct_2}\|(Q_{t_1-t_2}-I)f\|_{L^2(\mu_\infty)}\int_0^{t_2} \ff {\th(s)} {\sq {t_2-s}}\d s+e^{Ct_1}\int_{t_2}^{t_1}\|Q_{t_1-s}f\|_{W^{1,2}_{\mu_\infty}}\th(s) \d s\right).
\end{align*}
Thus, following \eqref{int-ths} and the same argument as proving assertions in (1), we have that
\beg{equation*}
\lim_{(t_1,z_1)\ra (t_2,z_2)}\left|G_f^{(2)}(t_1,z_1)-G_f^{(2)}(t_2,z_2)\right|=0,~\left\{\beg{array}{l}
t_2>0,~f\in L^2(\mu_\infty),\\
t_2=0,~f\in W^{1,2}_{\mu_\infty}.
\end{array}\right.
\end{equation*}
As a consequence, for $f\in C^2_b$, if $P_\cdot f(\cdot)\in C^{0,2}([0,+\infty)\times\R^d)$,  then there is $Q_\cdot f(\cdot)\in C^{0,2}([0,+\infty)\times\R^d)$.

\end{proof}

\subsection{Semigroup generated by $L_{\mu_\infty}+\bar A$}\label{4.3}

We come back to investigate the perturbation of $P^{\mu_\infty}_t$ by $\bar A$, and also denote by $Q_t$ the semigroup generated by $L_{\mu_\infty}+\bar A$.
\beg{cor}\label{cor-L+A}
Assume that $\mu_\infty$ is an invariant probability measure of \eqref{equ-mu-inf}. \\
(1) If the coefficient $(b(\cdot,\mu_\infty),\si(\cdot))$ of \eqref{equ-mu-inf} satisfies  {\bf (A1)}, {\bf (A2)}, and the assumption of (3) of Lemma \ref{lem-Malli}, and the integral kernel $D^F_{\mu_\infty}b\in L^2(\mu_\infty\times\mu_\infty)$, then the assumption in Theorem \ref{pert} holds for $P^{\mu_\infty}_t$ and $\bar A$.\\ 
(2) If, in addition,  \eqref{LypuV0} holds  and $D^F_{\mu_\infty}b$ satisfies \eqref{DF-DF} for some $\ph_0\in\sG$  and $1\leq V_0\in C^2$ satisfying \eqref{nnb1-Lypu1} with $U_0$ replaced by $V_0$ and 
\beg{align}\label{xphV0-p}
(1+|\cdot|+\ph_0(|\cdot|))V_0\in L^{p}(\mu_\infty),~p\geq 2,
\end{align} 
then Theorem \ref{thm-regu} holds for $P^{\mu_\infty}_t$ and $\bar A$. Moreover, following properties hold for $Q_tf$.\\
(a) $C_b^2\subset \sD(L_{\mu_\infty}+\bar A)$, $Q_t\1=\1$ and for any $f\in C_b^2$, there hold that $Q_\cdot f(\cdot)\in C^{1,2}([0,+\infty)\times \R^d)$ and $Q_tf$ is a classical solution of the following equation 
\beg{equation}\label{equ-clas}
\pp_t v_t(x)= \ff 1 2{\rm Tr}(\si \si^* \nn^2v_t)(x)+b(x,\mu_\infty)\cdot \nn v_t(x)+(\bar A v_t)(x),~v_0(x)=f(x).
\end{equation}
(b) There is $C>0$  such that for any $t\geq 0$ and $f\in L^2(\mu_\infty)$, 
\beg{align}\label{0Q1-pp}
|Q_tf(x)|\leq  P_t^{\mu_\infty}|f|(x)+Ce^{Ct}(1+\ph_0(|x|))V_0(x)\|f-\mu_\infty(f)\|_{L^2(\mu_\infty)},~x\in\R^d.
\end{align}
For any $t>0$ and any $f\in \sB(\R^d)$ with $|f|\leq V$ for some $1\leq V\in C^0\cap L^2(\mu_\infty)$ satisfying \eqref{EVp} with $P_t^Z$ replaced by $P_t^{\mu_\infty}$, there is $Q_tf\in C^2$, and for any $q\in (1,+\infty)$, there exists $C>0$, which is independent of $t$, $x$ and $f$,  such that 
\beg{align}\label{nnQ1-pp}
|\nn Q_tf(x)|&\leq \ff {C} {\sq{t\we 1}}\inf_{c\in\R}\left(P_t^{\mu_\infty}|f+c|^q(x)\right)^{\ff 1 q}+Ce^{Ct}V_0(x)(1+|x|)^{ \al_1\we \ff {2\al_1} {(1+\al_1)\al_2}}\|\bar f\|_{L^2(\mu_\infty)},\\
|\nn^2 Q_tf(x)|&\leq \ff {C} {t\we 1}\inf_{c\in\R}\left(P_t^{\mu_\infty}| f+c|^q(x)\right)^{\ff 1 q}+Ce^{Ct}\Th(t)V_0(x)(1+|x|)^{\al_1\we \ff {2\al_1} {(1+\al_1)\al_2}}\|\bar f\|_{L^2(\mu_\infty)},\label{nnQ2-pp}
\end{align}
where $\bar f=f-\mu_\infty(f)$ and 
\[\Th(t)=\int_0^t \ff {\ph_0\left(s^{\ff 1 2\we\ff {(1+\al_1)\al_2-2\al_1} {2(1+\al_1)\al_2}}\right)} {s\sq{t-s}}\d s,~t>0. \]
\end{cor}
Since $\mu_\infty$ is an invariant probability measure of $P^{\mu_\infty}_t$, $P_t^{\mu_\infty}$ can be extended to be a $C_0$-semigroup of contractions on $L^2(\mu_\infty)$. The proof of Corollary \ref{cor-L+A} is divided into following lemmas. 
\beg{lem}\label{den-C1-W12}
$C_0^1$ is dense in $W^{1,2}_{\mu_\infty}$.
\end{lem}
\beg{proof}
Let $N\geq 1$, $\ze_N\in C^2$ such that $\1_{[|x|\leq N]}\leq \ze_N(x)\leq \1_{[|x|\leq N+1]}$ and $|\nn \ze_N|\leq 2$. For any $f\in W^{1,2}_{\mu_\infty}$,  we find that $f\ze_N\in W^{1,2}_0(B_{N+1})$, where $B_{N+1}=\{x~|~|x|\leq N+1\}$, and
\beg{align*}
&\varlimsup_{n\ra +\infty}\left(\|f-f\ze_N\|_{L^2(\mu_\infty)}+\|\nn(f-f\ze_N)\|_{L^2(\mu_\infty)}\right)\\
&\quad\,\leq \lim_{N\ra +\infty}\left(\|f\1_{[|\cdot|\geq N]}\|_{L^2(\mu_\infty)}+ \|(\nn f)\1_{[|\cdot|\geq N]}\|_{L^2(\mu_\infty)}+2\|f\1_{[N\leq |\cdot|\leq N+1]}\|_{L^2(\mu_\infty)}\right)=0.
\end{align*}
There are $g_n\in C^1_0$ with $\supp(g_n)\subset B_{N+1}$ such that  $g_n \xlongrightarrow[~]{W^{1,2}} f\ze_N$. Note that the coefficients of $L_{\mu_\infty}$ are locally Lipschtiz and $\si$ is invertible,  and $\mu_\infty$ satisfies 
\[\mu_\infty(L_{\mu_\infty}f)=0,~f\in C_0^2.\]
It follows from \cite[Theorem 1.2.2]{BKR} or \cite[Theorem 2]{Sj} that $\mu_\infty(\d x)\ll \d x$ and the density $\ff {\d \mu_\infty} {\d x}$ is continuous.  Thus $\ff {\d\mu_\infty} {\d x}$ is locally bounded.  Then, we have that
\[\varlimsup_{n\ra +\infty}\|g_n-f\ze_N\|_{W^{1,2}_{\mu_\infty}}\leq C_{N+1}\lim_{n\ra +\infty}\|g_n-f\ze_N\|_{W^{1,2}}=0.\]
Hence, $C_0^1$ is dense in $W^{1,2}_{\mu_\infty}$.

\end{proof}

Next lemma shows that $P_t^{\mu_\infty}$ satisfies {\bf (B1)}.
\beg{lem}\label{lem-PWnn}
Assume $(b(\cdot,\mu_\infty),\si(\cdot))$ satisfy {\bf (A1)}, {\bf (A2)}, and the assumption of (3) of Lemma \ref{lem-Malli}. Then $P_t^{\mu_\infty}$ is strongly continuous on $W^{1,2}_{\mu_\infty}$ and satisfies \eqref{nnP-L2}.
\end{lem}

\beg{proof}
For all $f\in C^1_b$ and $|v|\leq 1$, due to \eqref{Eest-nnY-2p}, there exist $C\geq 1,w\geq 0$ such that
\beg{equation}\label{nnp-pnn}
|\nn_vP_t^{\mu_\infty}f(x)|=|\E\<\nn f(X_t^{\mu_\infty}),\nn_vX_t^{\mu_\infty}\>|\leq Ce^{wt}\sq{P_t^{\mu_\infty}|\nn f|^2(x)},~a.e.~x\in\R^d.
\end{equation}
This implies that
\beg{equation}\label{Pt-lip-grow}
\left\| \nn P_t^{\mu_\infty} f \right\|_{L^2(\mu_\infty)}^2\leq C^2 e^{2wt}\mu_\infty(P_t^{\mu_\infty}|\nn f|^2)\leq C^2 e^{2wt} \|\nn f\|_{L^2(\mu_\infty)}^2.
\end{equation}
This, together with $\mu_\infty(|P^{\mu_\infty}_tf|^2)\leq \mu_\infty(f^2)$, yields that 
\beg{equation}\label{Pt-W12}
\|P_t^{\mu_\infty}f\|_{W^{1,2}_{\mu_\infty}}\leq (Ce^{wt}\vee 1)\|f\|_{W^{1,2}_{\mu_\infty}}\leq Ce^{wt}\|f\|_{W^{1,2}_{\mu_\infty}},~f\in C_b^1.
\end{equation}
For $f\in C^1_b$ and $|v|\leq 1$, we have that 
\beg{align*}
\left|\nn_v P_t^{\mu_\infty} f(x)-\nn_v f(x)\right|&\leq \left|\E\<\nn f(X_t^{\mu_\infty}(x))-\nn f(x),\nn_v X_t^{\mu_\infty}(x)\>\right|\\
&\quad\,+ \left|(\nn f(x))\left(\E\nn_v X_t^{\mu_\infty}(x)-v\right)\right|\\
&\leq C e^{wt}\left(\E|\nn f(X_t^{\mu_\infty}(x))-\nn f(x)|^2\right)^{\ff 1 2}\\
&\quad\, +\|\nn f\|_\infty\left(\E|\nn X_t^{\mu_\infty}(x)-I|^2\right)^{\ff 1 2}.
\end{align*}
This, together with \eqref{cont-nn1-ad} and the dominated convergence theorem, yields that 
\[\lim_{t\ra 0^+}\mu_{\infty}\left(\left|\nn P_t^{\mu_\infty} f -\nn  f\right|^2\right)=0,~f\in C^1_b.\] 
Combining this with \eqref{Pt-W12} and the approximation argument, we prove that $P_t^{\mu_\infty}$ is strong continuous on $W^{1,2}_{\mu_\infty}$. Since $f$ is bounded,  \eqref{EVp} holds with $V\equiv \|f\|_\infty$. Then \eqref{nnpzf} holds for $f$,  and \eqref{nnP-L2} can be derived from \eqref{nnpzf} and the density of $C^1_b$ in $L^{2}(\mu_\infty)$.

\end{proof}
 
The following lemma shows that $P_t^{\mu_\infty}$ and $\bar A$ satisfy {\bf (B2)}. 
\beg{lem}\label{lem-PAB2}
The same assumption as Lemma \ref{lem-PWnn} holds. If the kernel $D^F_{\mu_\infty}b\in L^2(\mu_\infty\times\mu_\infty)$,  then $P_t^{\mu_\infty}$ and $\bar A$ satisfy {\bf (B2)}.  Particularly, if \eqref{DF-DF} holds with $F\in L^2(\mu_\infty)$ and $(1+\ph_0(|\cdot|))V_0(\cdot)\in L^2(\mu_\infty)$, then $D^F_{\mu_\infty}b\in L^2(\mu_\infty\times\mu_\infty)$.
\end{lem}
\beg{proof}
If $D^F_{\mu_\infty}b\in L^2(\mu_\infty\times\mu_\infty)$, then $\bar A$ is a compact operator from $W^{1,2}_{\mu_\infty}$ to $L^2(\mu_\infty)$. Thus Remark \ref{rem-comp} yields that $P_t^{\mu_\infty}$ and $\bar A$ satisfy {\bf (B2)}.  In particular, it is clear that \eqref{DF-DF} and $F,(1+\ph_0(|\cdot|))V_0\in L^2(\mu_\infty)$  imply $D^F_{\mu_\infty}b\in L^2(\mu_\infty\times\mu_\infty)$. 

\end{proof}

Up to now, we have shown that assumptions in Theorem \ref{pert} hold for $P_t^{\mu_\infty}$ and $\bar A$. Next, we check assumptions in Theorem \ref{thm-regu}.
\beg{lem}\label{nnPS-nn2PS}
All assumptions of  Corollary \ref{cor-L+A} hold. Then $P_t^{\mu_\infty}\bar A$ satisfies \eqref{nnPS-p} and \eqref{nn2PS-p} with
\beg{align*}
V_1(z)&=V_2(z)=V_0(z)(1+|z|)^{\al_1\we \ff {2\al_1} {(1+\al_1)\al_2}},\\
\th(t)&=t^{-1}\ph\left(t^{\ff 1 2\we\ff {(1+\al_1)\al_2-2\al_1} {2(1+\al_1)\al_2}}\right).
\end{align*}
Moreover, $V_1,V_2\in L^2(\mu_\infty)$ and $\th\in L^1_{loc}([0,+\infty))$ satisfies \eqref{int-ths}.
\end{lem}

\beg{proof}
For all $f\in W^{1,2}_{\mu_\infty}$, it follows from \eqref{DF-DF} that
\beg{equation}\label{Af-Af}
\beg{split}
\left|\bar Af(z_1)-\bar Af(z_2)\right|&\leq \int_{\R^d} \left| D^F_{\mu_\infty}b(x,z_1)- D^F_{\mu_\infty}b(x,z_2)\right|\cdot|\nn f(x)|\mu_\infty(\d x)\\
&\leq \mu_\infty\left(F|\nn f|\right)\ph_0(|z_1-z_2|)(V_0(z_1)+V_0(z_2))\\
&\leq \|F\|_{L^2(\mu_\infty)}\|\nn f\|_{L^2(\mu_\infty)}\ph_0(|z_1-z_2|)(V_0(z_1)+V_0(z_2)).
\end{split}
\end{equation}
It follows from the It\^o formula, \eqref{nnb1-Lypu1} with $U_0$ replaced by $V_0$ and  \eqref{LypuV0} that for any $p\geq 1$, there is $C>0$ such that 
\[\E V_0(X^{\mu_\infty}_t(x))^p\leq e^{Ct}V_0(x)^p,~t\geq 0,~x\in\R^d.\]
Combining this with \eqref{Af-Af}, {\bf (A1)}, {\bf (A2)}, and the assumption of (3) of Lemma \ref{lem-Malli}, we can derive from Lemma \ref{lem-grad-est}  that $P_t^{\mu_\infty}\bar A f\in C^2(\R^d)$ and 
\beg{align*}
|\nn P_t^{\mu_\infty}\bar Af(z) |&\leq \ff {Ce^{Ct}} {\sq t}\|f\|_{W^{1,2}_{\mu_\infty}}V_0(z)(1+|z|)^{\al_1\we \ff {2\al_1} {(1+\al_1)\al_2}}\ph_0\left(t^{\ff 1 2\we \ff {(1+\al_1)\al_2-2\al_1} {2(1+\al_1)\al_2}}\right)\\
&\leq  \ff {Ce^{Ct}} {\sq t}(1+t)^{\ff 1 2}\|f\|_{W^{1,2}_{\mu_\infty}}V_0(z)(1+|z|)^{\al_1\we \ff {2\al_1} {(1+\al_1)\al_2}}\\
&\leq  \ff {Ce^{Ct}} {\sq t}\|f\|_{W^{1,2}_{\mu_\infty}}V_0(z)(1+|z|)^{\al_1\we \ff {2\al_1} {(1+\al_1)\al_2}},
\end{align*}
where $C$ is a positive constant whose value may vary at each line.  Thus \eqref{nnPS-p} holds with $V_1(z)=V_0(z)(1+|z|)^{\al_1\we \ff {2\al_1} {(1+\al_1)\al_2}}$. Similarly, combining \eqref{Af-Af} with \eqref{nn2Pph}, we find that
\beg{align*}
|\nn^2 P_t^{\mu_\infty}\bar Af(z) |&\leq Ce^{Ct}V_0(z)(1+|z|)^{\al_1\we \ff {2\al_1} {(1+\al_1)\al_2}}t^{-1}\ph_0\left(t^{\ff 1 2\we \ff {(1+\al_1)\al_2-2\al_1} {2(1+\al_1)\al_2}}\right)\|\nn f\|_{L^2(\mu_\infty)}.
\end{align*}
Hence, \eqref{nn2PS-p} holds with $\th(t)=t^{-1}\ph_0\left(t^{\ff 1 2\we\ff {(1+\al_1)\al_2-2\al_1} {2(1+\al_1)\al_2}}\right)$, which is integrable since 
\beg{equation*}
\int_0^1t^{-1}\ph\left(t^{\ff 1 2\we\ff {(1+\al_1)\al_2-2\al_1} {2(1+\al_1)\al_2}}\right)\d t=2\left(\ff  {(1+\al_1)\al_2} {(1+\al_1)\al_2-2\al_1} \vee 1\right)\int_0^1 u^{-1}\ph\left(u\right)\d u<+\infty.
\end{equation*}
The H\"older inequality and  \eqref{xphV0-p} imply that $V_1,V_2\in L^2(\mu_\infty)$.  For $\th$, there is 
\beg{align*}
\Th(t)=\int_0^t \ff {\th(s)} {\sq{t-s}}\d s\leq \sq{\ff 2 t} \int_0^{t/2} \th(s)\d s+\ff {2\th(t)} t\int_{t/2}^t\ff {\d s} {\sq{t-s}}<+\infty,~t>0.
\end{align*}

\end{proof}

Up to now, we have proven that assumptions in Theorem \ref{pert} and Theorem \ref{thm-regu} hold for $P^{\mu_\infty}_t$ and $\bar A$. Finally, we prove assertions of (a) and (b) in Corollary \ref{cor-L+A}.
\beg{lem}
Under all the assumptions of Corollary \ref{cor-L+A}, assertions of (a) and (b) in Corollary \ref{cor-L+A} hold.
\end{lem}
\beg{proof}
We first prove assertions of (a).  For $f\in C^2_b$, it is clear that $L_{\mu_\infty} f\in L^2(\mu_\infty)$ since \eqref{xphV0-p}, {\bf (A1)} and {\bf (A2)}. It follows from the It\^o formula that
\[\ff {P_t^{\mu_\infty}f(x)-f(x)} t-L_{\mu_\infty}f(x)=\ff 1 t\int_0^t\left( \E L_{\mu_\infty}f(X^{\mu_\infty}_s(x))-L_{\mu_\infty}f(x)\right)\d s,~x\in\R^d.\]
Then
\[\left\|\ff {P_t^{\mu_\infty}f -f } t-L_{\mu_\infty}f \right\|_{L^2(\mu_\infty)}\leq \ff 1 t\int_0^t\|P_s^{\mu_\infty}L_{\mu_\infty}f-L_{\mu_\infty}f\|_{L^2(\mu_\infty)}\d s.\]
This, together with the dominated convergence and  that $P^{\mu_\infty}_t$ is a $C_0$-semigroup of contractions on $L^2(\mu_\infty)$, implies 
\[\lim_{t\ra 0^+}\left\|\ff {P_t^{\mu_\infty}f -f } t-L_{\mu_\infty}f \right\|_{L^2(\mu_\infty)}=0.\]
Thus, $C^2_b\subset\sD(L_{\mu_\infty})$.  Combining this with Lemma \ref{lem-PWnn} and Lemma \ref{lem-PAB2}, we can derive from Theorem \ref{pert} that 
\[\sD(L_{\mu_\infty}+\bar A)=\sD(L_{\mu_\infty}) \supset C_b^2.\]
Hence,
\beg{align*}
\ff {\d } {\d t}Q_tf&=(L_{\mu_\infty}+\bar A) Q_tf=\ff 1 2{\rm Tr}(\si \si^* \nn^2Q_tf) +b(\cdot,\mu_\infty)\cdot \nn Q_tf +\bar A Q_tf .
\end{align*}
Due to Lemma \ref{lem-PAB2} and Lemma \ref{nnPS-nn2PS}, we derive from  Theorem \ref{thm-regu} that $L_{\mu_\infty} Q_\cdot f(\cdot)\in C^{0}([0,+\infty)\times\R^d)$. $\bar A$ is continuous from $W^{1,2}_{\mu_\infty}$ to $\sG_{V_0,\ph_0}$ since \eqref{Af-Af}. This, together with that $Q_\cdot f$ is continuous on $W^{1,2}_{\mu_\infty}$ since  $f\in C^2_b\subset W^{1,2}_{\mu_\infty}$, we have that $\bar A Q_\cdot f(\cdot) \in C^{0}([0,+\infty)\times\R^d)$. Hence, $\ff {\d } {\d t}Q_tf\in C^{0}([0,+\infty)\times\R^d)$. Therefore, $Q_\cdot f(\cdot)\in C^{1,2}([0,+\infty)\times \R^d)$ and $Q_tf$ satisfies \eqref{equ-clas}.\\
Since  $\1\in C_b^2$ and
\[\ff {\d } {\d t}Q_t\1=Q_t\left(L_{\mu_\infty}+\bar A\right)\1=0,\] 
we have that 
\beg{equation}\label{Q1=1}
Q_t\1=\1.
\end{equation}

Next, we prove assertions of (2).  Taking into account that 
\[P_t^{\mu_\infty}\bar A f=P_t^{\mu_\infty}\bar A (f+c)~,c\in\R,~f\in W^{1,2}_{\mu_\infty},\]
\eqref{Qt-Pt-int} can be rewritten into the following form
\[Q_tf=P_t^{\mu_\infty}f+\int_0^t P_s^{\mu_\infty}\bar A Q_{t-s}(f-\mu_\infty(f))\d s,~f\in L^2(\mu_\infty),\]
and $\|f\|_{L^2(\mu_\infty)}$ in the right hand side of \eqref{nnG1} and  \eqref{nnG2} can be replaced by $\|f-\mu_\infty(f)\|_{L^2(\mu_\infty)}$. Combining these  with Lemma \ref{lem-grad-est}, Lemma \ref{nnPS-nn2PS},  Theorem \ref{thm-regu} and 
\[\nn Q_tf(x)=\nn Q_t(f+c)(x),\qquad \nn^2 Q_tf(x)=\nn^2 Q_t(f+c)(x),~x\in\R^d,~c\in\R,\] 
we have that $Q_tf\in C^2$, and \eqref{nnQ1-pp} and \eqref{nnQ2-pp} hold. 
For \eqref{0Q1-pp}, it follows from \eqref{Af-Af} and $\mu_\infty(\bar A f)=0$ that
\beg{align*}
|\bar A f(x)|&=|\bar Af(x)-\mu_\infty(\bar Af)|\leq \int_{\R^d}|\bar A f(x)-\bar A f(y)|\mu_\infty(\d y)\\
&\leq \|F\|_{L^2(\mu_\infty)}\|\nn f\|_{L^2(\mu_\infty)}\int_{\R^d}\ph_0(|x-y|)(V_0(x)+V_0(y))\mu_\infty(\d y).
\end{align*}
Thus,
\beg{align*}
|P_t^{\mu_\infty}\bar A f(x)|&\leq \E|\bar A f(X_t^{\mu_\infty}(x))|\\
&\leq C\|\nn f\|_{L^2(\mu_\infty)}\E\int_{\R^d}\ph_0(|X_t^{\mu_\infty}(x)-y|)(V_0(X_t^{\mu_\infty}(x))+V_0(y))\mu_\infty(\d y)\\
&\leq C\|\nn f\|_{L^2(\mu_\infty)} \left(\E\left(\mu_\infty(\ph_0^2(|X_t^{\mu_\infty}(x)-\cdot|)) \right)\right)^{\ff 1 2}\left(\E V_0(X_t^{\mu_\infty}(x))^{2}+\mu_\infty(V_0^2)\right)^{\ff 1 2}\\
&\leq Ce^{Ct}\|\nn f\|_{L^2(\mu_\infty)} V_0(x)\ph_0\left(\left(\E(\mu_\infty(|X_t^{\mu_\infty}(x)-\cdot|^2)\right)^{\ff 1 2}\right)\\
&\leq Ce^{Ct}\|\nn f\|_{L^2(\mu_\infty)} V_0(x)\ph_0\left(2^{\ff 1 2}\left(\E|X_t^{\mu_\infty}(x)|^2+\mu_\infty(|\cdot|^2)\right)^{\ff 1 2}\right)\\
&\leq Ce^{Ct}\|\nn f\|_{L^2(\mu_\infty)} V_0(x)\ph_0\left(Ce^{Ct}(1+|x|)\right)\\
&\leq Ce^{Ct}\|\nn f\|_{L^2(\mu_\infty)} (1+\ph_0\left(|x|\right))V_0(x),
\end{align*}
where $C$ is a constant whose value may vary at each line. Combining this with \eqref{Q1=1}, \eqref{Qt-Pt-int} and \eqref{Q-W12}, the inequality \eqref{0Q1-pp} follows.

\end{proof}

\section{Proof of Theorem \ref{thm0}}

Because \eqref{Inequ-nnb1}, \eqref{Inequ-nnb2} and {\bf (H2)}  imply that $(b(\cdot,\mu_\infty),\si(\cdot))$ satisfies {\bf (A1)}, {\bf (A2)} with $\al_1=0$ and $\al_2\in (1,2)$, and the assumption of (3) of Lemma \ref{lem-Malli}, the first assertion except \eqref{invQtmu} of Theorem \ref{thm0} can be derived from (1) of Corollary \ref{cor-L+A}.    For \eqref{invQtmu}, we find that
\beg{align*}
\mu_\infty(Q_tf)-\mu_\infty(f)&=\mu_\infty\left(\int_0^t(L_{\mu_\infty}+\bar A)Q_sf\d s\right),~f\in C_b^2.
\end{align*}
Due to Remark \ref{re-DF}, $\mu_\infty(D^F_{\mu_\infty}b(x,\cdot))=0$ for any $x\in\R^d$. Since \eqref{ph-p-q0} and \eqref{mu-inf-U0}, there is $\mu_\infty(V_0^p)<+\infty$ for any $p\geq 1$.  Then, \eqref{nnQ1-pp}  and the Fubini theorem imply $\mu_\infty(\bar A Q_{s}f)=0$, and 
\[\mu_\infty\left(\int_0^t\bar AQ_{s}f\d s\right)=\int_0^t \mu_\infty\left(\bar AQ_{s}f\right)\d s=0.\] 
It follows from $f\in C^2_b$ and Corollary \ref{cor-L+A} that $Q_sf\in C^2\cap L^2(\mu_\infty)$. Moreover, \eqref{Inequ-nnb2}, \eqref{nn2Pph} with $\ph(r)=r$ and $V_1\equiv 1$, \eqref{nnQ1-pp}, \eqref{Gf2W12} and Lemma \ref{nnPS-nn2PS}  with $\al_1=0$, there is 
\[|L_{\mu_\infty}Q_{s}f(x)|\leq C_sV_0(x)(1+|x|)^{ \be_1+2},~x\in\R^d,~s>0,\]
where
\beg{align*}
C_s=Ce^{Cs}\left(1+\ff 1 {\sq s}(\|f\|_\infty+\|\nn f\|_{\infty}) + \int_0^s \ff {\ph(\sq r)} r\d r\|f\|_{W^{1,2}_{\mu_\infty}}\right),~s\geq 0.
\end{align*}
Thus, by using {\bf (H0)}, \eqref{mu-inf-U0} and \eqref{ph-p-q0}, we find that $\|L_{\mu_\infty}Q_\cdot f\|_{L^2(\mu_\infty)}\in L^1_{loc}([0,+\infty)$, and
\[\mu_\infty(L_{\mu_\infty}Q_sf)=0,~s>0,\]
since  $\mu_\infty$ is the invariant probability measure of $P^{\mu_\infty}_t$.  Then the Fubini theorem implies that
\[\mu_\infty\left(\int_0^t  L_{\mu_\infty} Q_sf\d s\right)=\int_0^t \mu_\infty\left(L_{\mu_\infty} Q_sf\right)\d s=0.\]
Consequently, $\mu_\infty(Q_tf)=\mu_\infty(f)$ for $f\in C^2_b$. By the approximation argument, \eqref{invQtmu} holds.

In this section, we focus on proving the rest assertions of Theorem \ref{thm0}. Due to {\bf (H0)}-{\bf (H3)}, for any $p\geq 1$ and $X_0$ with $\sL_{X_0}=\mu_0\in\sP_{U_0^p}$, \eqref{main-equ} has a unique strong solution with $\sL_{X_\cdot}\in C([0,T],\sP_{U_0^p})$ for any $T>0$. Denote by $\mu_t=\sL_{X_t}$. The decoupled equation \eqref{freezing} has a unique strong solution for $X^{\mu}_0=x\in\R^d$ since {\bf (H1)}-{\bf (H3)} again. {\bf (H0)}-{\bf (H2)} imply the strong wellposedness of \eqref{equ-mu-inf}.  Moreover, \eqref{xphV0-p} follows from {\bf (H0)}, \eqref{mu-inf-U0} and \eqref{ph-p-q0}. Hence, all conclusions of Corollary \ref{cor-L+A} hold under {\bf (H0)}-{\bf (H3)} and \eqref{DF-DF}. Denote 
\[D^F_{\mu,\nu}u(x)=\int_0^1D^F_{r\mu+(1-r)\nu} u(x)\d r.\]
We first establish the following crucial  equality for $\mu_t$ and $\mu_\infty$. 
\beg{lem}\label{Dhr0}
Assume {\bf (H0)}-{\bf (H3)} and \eqref{DF-DF} hold. For any $f\in C^2_b$ and $\mu_0\in\sP_{U_0}$, there is 
\beg{equation}\label{DuHam}
\beg{split}
&(\mu_t-\mu_\infty)(f)=(\mu_0-\mu_\infty)(Q_tf)+ \int_0^t (\mu_s-\mu_\infty)((b(\cdot,\mu_s)-b(\cdot,\mu_\infty))\cdot\nn Q_{t-s}f(\cdot))\d s\\
&+\int_0^t \int_{\R^d} (\mu_s-\mu_\infty)\left( D^F_{\mu_s,\mu_\infty}b(x,\cdot)-D^F_{\mu_\infty}b(x,\cdot)\right)\cdot\nn Q_{t-s}f(x) \mu_\infty(\d x)\d s. 
\end{split}
\end{equation}
\end{lem}
\beg{proof}
Since $\mu_t(f)=\mu_0(P^{\mu}_tf)$, we have that
\beg{align*}
\mu_t(f)-\mu_0(Q_tf)=\mu_0(P^{\mu}_tf)-\mu_0(Q_tf).
\end{align*}
It follows from Corollary \ref{cor-L+A} and $f\in C_b^2$ that $Q_\cdot f(\cdot)\in C^{1,2}([0,+\infty)\times\R^d)$. Then, by using  the It\^o formula, there is  
\beg{align*}
\d Q_{t-s}f(X_s^\mu(x))&=-(L_{\mu_\infty}+\bar A)Q_{t-s}f(X_s^\mu(x))\d s+ L_{\mu_s}Q_{t-s}f(X_s^\mu(x))\d s\\
&\quad\,+\<\nn Q_{t-s}f(X_s^\mu(x)),\si(X_s^{\mu}(x))\d B_s\>\\
&= (b(X_s^\mu(x),\mu_s)-b(X_s^\mu(x),\mu_\infty))\cdot\nn Q_{t-s}f(X_s^\mu(x))\d s-\bar AQ_{t-s}f(X_s^\mu(x))\d s\\
&\quad\,+\<\nn Q_{t-s}f(X_s^\mu(x)),\si(X_s^{\mu}(x))\d B_s\>.
\end{align*}
By using {\bf (H2)}, \eqref{nnQ1-pp} with $\al_1=0$, \eqref{b1-mu-nu}, \eqref{DF-DF} and $f\in C^2_b$, we have that 
\beg{align*}
|(b(x,\mu_s)-b(x,\mu_\infty))\cdot\nn Q_{t-s}f(x)|&\leq Ce^{Ct}\left(\ff {1} {\sq{(t-s)\we 1}}+V_0(x)\right)\|f\|_\infty\|\mu_s-\mu_\infty\|_{V_0,\ph_0},\\
|\si(x)^*\nn Q_{t-s}f(x)|&\leq  Ce^{Ct}\left(\ff {1} {\sq{(t-s)\we 1}}+V_0(x)\right)\|f\|_\infty,\\
|\bar AQ_{t-s}f(x)|&\leq \mu_\infty(F|\nn Q_{t-s}f|)(1+\ph_0(x))V_0(x)\\
&\leq C\mu_\infty(F(1+V_0))e^{Ct} \ff {(1+\ph_0(x))V_0(x)} {\sq{(t-s)\we 1}}\|f\|_{L^2(\mu_\infty)}.
\end{align*}
Combining these with $\mu_\cdot\in C([0,T],\sP_{U_0})$, \eqref{ph-p-q0},  \eqref{norm-p0q0}, and the Fubini theorem, we have that
\begin{equation*}
\begin{split}
\mu_0(P_t^{\mu} f)-\mu_0(Q_tf)&=\int_{\R^d}\E \left(Q_{t-s}f(X_s^{\mu}(x))|_{s=t}- Q_{t-s}f(X_s^{\mu}(x))|_{s=0}\right)\mu_0(\d x)\\
&= \int_0^t \int_{\R^d}(b(x,\mu_s)-b(x,\mu_\infty))\cdot\nn Q_{t-s}f(x)\mu_s(\d x)\d s\\
&\quad\,-\int_0^t\int_{\R^d}(\bar AQ_{t-s}f)(x)\mu_s(\d x)\d s.
\end{split}
\end{equation*}
Combining this with \eqref{invQtmu}, we therefore have that
\beg{align*}
&(\mu_t-\mu_\infty)(f)-(\mu_0-\mu_\infty)(Q_tf)\\
&\quad\,=\mu_t(f)-\mu_0(Q_tf)+(\mu_\infty(Q_tf)-\mu_\infty(f))\\
&\quad\,=\int_0^t \mu_s((b(\cdot,\mu_s)-b(\cdot,\mu_\infty))\cdot\nn Q_{t-s}f(\cdot))\d s-\int_0^t (\mu_s-\mu_\infty)(\bar AQ_{t-s}f)\d s\\
&\quad\, =\int_0^t (\mu_s-\mu_\infty)((b(\cdot, \mu_s)-b(\cdot, \mu_\infty))\cdot\nn Q_{t-s}f(\cdot) )\d s\\
&\quad\,\quad\,+\int_0^t\mu_\infty((b(\cdot, \mu_s)-b(\cdot, \mu_\infty))\cdot\nn Q_{t-s}f(\cdot) )\d s\\
&\quad\,\quad\, -\int_0^t \int_{\R^d}\left(\mu_\infty\left( D^F_{\mu_\infty}b(\cdot,z)\cdot\nn Q_{t-s}f(\cdot)\right)\right)(\mu_s-\mu_\infty)(\d z)\d s\\
&\quad\, =\int_0^t (\mu_s-\mu_\infty)((b(\cdot,\mu_s)-b(\cdot,\mu_\infty))\cdot\nn Q_{t-s}f(\cdot))\d s\\
&\quad\,\quad\,+\int_0^t\mu_\infty((b(\cdot,\mu_s)-b(\cdot,\mu_\infty))\cdot\nn Q_{t-s}f(\cdot))\d s\\
&\quad\,\quad\, -\int_0^t \int_{\R^d}\left( (\mu_s-\mu_\infty)\left(  D^F_{\mu_\infty}b(x,\cdot)\right)\cdot\nn Q_{t-s}f(x)\right)\mu_\infty(\d x)\d s\\
&\quad\, =\int_0^t (\mu_s-\mu_\infty)((b(\cdot,\mu_s)-b(\cdot,\mu_\infty))\cdot\nn Q_{t-s}f(\cdot))\d s\\
&\quad\,\quad\,+\int_0^t\int_{\R^{d}}((\mu_{s}-\mu_{\infty})(D^{F}_{\mu_{s},\mu_{\infty}}b(x,\cdot))\cdot\nn Q_{t-s}f(\cdot))\mu_{\infty}(\d x)\d s\\
&\quad\,\quad\, -\int_0^t \int_{\R^d}\left( (\mu_s-\mu_\infty)\left(  D^F_{\mu_\infty}b(x,\cdot)\right)\cdot\nn Q_{t-s}f(x)\right)\mu_\infty(\d x)\d s\\
&\quad\, =\int_0^t (\mu_s-\mu_\infty)((b(\cdot,\mu_s)-b(\cdot,\mu_\infty))\cdot\nn Q_{t-s}f(\cdot))\d s\\
&\quad\,\quad\, +\int_0^t \int_{\R^d} (\mu_s-\mu_\infty)\left( D^F_{\mu_s,\mu_\infty}b(x,\cdot)- D^F_{\mu_\infty}b(x,\cdot)\right)\cdot\nn Q_{t-s}f(x) \mu_\infty(\d x)\d s.
\end{align*}

\end{proof}

\beg{lem}\label{exp-con-Q}
Assume {\bf (H0)}-{\bf (H3)} hold, and $D^F_{\mu_\infty} b(x,z)$ satisfies \eqref{DF-DF}. Suppose \eqref{con-W1} and \eqref{Q-exp-dec} holds. Then there is $C\geq 1$, which depends on $C_{Q},\la_Q,\la_P$, $C_W$, $\ph_0$, $V_0$ and $\mu_\infty$ such that
\beg{equation}\label{exp-Lip0}
\|Q_tf\|_{V_0,\ph_0}\leq CH(t)e^{-(\la_P\we\la_Q) t}  \|f\|_{V_0,\ph_0},~f\in \mathscr{G}_{V_0,\ph_0},
\end{equation}
where
\beg{align*}
H(t)=\left(\ff {1} {|\la_p-\la_Q|\we 1}\we (t\vee \sq t)\right)\1_{[\la_P\neq\la_Q]}+(t\vee \sq t)\1_{[\la_P=\la_Q]},~t\geq 0.
\end{align*}
\end{lem}
\beg{proof}
It follows from \eqref{Q-W12} that there is $C_1>0$ such that
\beg{align*}
\|\nn Q_tf\|_{L^2(\mu_\infty)}\leq \ff {C_1} {\sq t} \|f\|_{L^2(\mu_\infty)}\leq \ff {C_1e^{\la_Q}} {\sq t} e^{-\la_Q t} \|f\|_{L^2(\mu_\infty)},~0<t\leq 1,~f\in L^2(\mu_\infty).
\end{align*}
Combining this with the semigroup property  and \eqref{Q-exp-dec},  we have  that
\beg{align*}
\|\nn Q_tf\|_{L^2(\mu_\infty)}&=\|\nn Q_1 Q_{t-1}f\|_{L^2(\mu_\infty)}\leq C_1 \|Q_{t-1}f\|_{L^2(\mu_\infty)}\\
&\leq C_1C_Qe^{-\la_Q(t-1)}\|f\|_{L^2(\mu_\infty)}\\
&=C_1C_Qe^{\la_Q}e^{-\la_Qt}\|f\|_{L^2(\mu_\infty)},~\mu_\infty(f)=0,~t\geq 1.
\end{align*}
Taking into account \eqref{Q1=1}, there is $\tld C>0$, which depends on $C_Q,\la_Q$, such that
\[\|\nn Q_tf\|_{L^2(\mu_\infty)}=\|\nn Q_t(f-\mu_\infty(f))\|_{L^2(\mu_\infty)}\leq \ff {\tld C} {\sq{t\we 1}} e^{-\la_Qt}\|f-\mu_\infty(f)\|_{L^2(\mu_\infty)},~t>0.\]
Combining this with \eqref{con-W1} and \eqref{Af-Af}, we find that for any $f\in\mathscr{G}_{V_0,\ph_0}$ 
\beg{align*}
\| Q_t f\|_{V_0,\ph_0}&\leq \| P^{\mu_\infty}_t f\|_{V_0,\ph_0}+\int_0^t \|P^{\mu_\infty}_{t-s}\bar A Q_sf\|_{V_0,\ph_0}\d s\\
&\leq C_We^{-\la_P t}\| f\|_{V_0,\ph_0}+C_W\int_0^t e^{-\la_P (t-s)} \| \bar A Q_s f\|_{V_0,\ph_0}\d s\\
& \leq C_We^{-\la_P t}\| f\|_{V_0,\ph_0}+ C_W\|F\|_{L^2(\mu_\infty)}\int_0^t e^{-\la_P (t-s)} \|\nn Q_s f\|_{L^2(\mu_\infty)}\d s\\
&\leq C_We^{-\la_P t}\| f\|_{V_0,\ph_0}+ C_W\|F\|_{L^2(\mu_\infty)}\tld C\int_0^t \ff {e^{-\la_P (t-s)-\la_Qs} } {\sq{s\we 1}}\|f-\mu_\infty(f)\|_{L^2(\mu_\infty)}\d s
\end{align*}
Note that
\beg{align}\label{add-inteee}
\int_0^t \ff {e^{-\la_P (t-s)-\la_Qs} } {\sq{s\we 1}}\d s&=e^{-\la_Pt}\int_0^t \ff {e^{(\la_P-\la_Q)s}} {\sq{s\we 1}}\d s\nonumber\\
&=e^{-\la_Pt}\left(2e^{(\la_P-\la_Q)^+t}\sq t\1_{[t\leq 1]}+\ff {e^{(\la_P-\la_Q)t}- e^{(\la_P-\la_Q)} } {\la_P-\la_Q}\1_{[t\geq 1]}\right)\1_{[\la_p\neq\la_Q]}\nonumber\\
&\quad\, + e^{-\la_Pt}(2\sq{t\we1}+(t-1)^+)\1_{[\la_P=\la_Q]}\\
&\leq \left(2+\ff 1 {|\la_p-\la_Q|}\right)e^{-(\la_P\we\la_Q)t}\1_{[\la_p\neq\la_Q]}+ e^{-\la_Pt}\varphi_1(t)\1_{[\la_P=\la_Q]}\nonumber
\end{align}
and
\beg{align}\label{int-ewe}
\int_0^t \ff {e^{-\la_P (t-s)-\la_Qs} } {\sq{s\we 1}}\d s&\leq \int_0^t \ff {e^{-(\la_P\we\la_Q) (t-s)-(\la_Q\we\la_P)s} } {\sq{s\we 1}}\d s\nonumber\\
&\leq e^{-(\la_P\we\la_Q)t}(2\sq{t\we1}+(t-1)^+).
\end{align}
We have by using \eqref{add-inteee} that
\beg{align*}
\| Q_t f\|_{V_0,\ph_0}&\leq C_We^{-\la_P t}\| f\|_{V_0,\ph_0}\\
&\quad\, + Ce^{-(\la_{P}\we \la_Q) t}\left(\ff {\1_{[\la_P\neq\la_Q]}} {|\la_p-\la_Q|\we 1}+\varphi_1(t)\1_{[\la_P=\la_Q]}\right) \|f-\mu_\infty(f)\|_{L^2(\mu_\infty)},
\end{align*}
where $\tld C$ is a constant which depends on  $C_Q,\la_Q,\la_P,C_W,\|F\|_{L^2(\mu_\infty)}$. Note that for $f\in\mathscr{G}_{V_0,\ph_0}$, we have 
\beg{equation}\label{ph-ppp}
\beg{split}
|f(x)-f(y)|&\leq \ph_0(|x-y|)(V_0(x)+V_0(y))\\
&\leq C(1+|x|+|y|)(V_0(x)+V_0(y)).
\end{split}
\end{equation}
Thus
\beg{align}\label{f-L2-ph}
\|f-\mu_\infty(f)\|_{L^2(\mu_\infty)}^2& = \int_{\R^d}\left(f(x)-\int_{\R^d}f(y)\mu_\infty(\d y)\right)^2\mu_\infty(\d x)\nonumber\\
& = \int_{\R^d}\left(\int_{\R^d}\left(f(x)-f(y)\right)\mu_\infty(\d y)\right)^2\mu_\infty(\d x)\\
& \leq C\int_{\R^d}\left(\int_{\R^d} (1+|x|+|y|)(V_0(x)+V_0(y)) \mu_{\infty}(\d y)\right)^2\mu_\infty(\d x)\|f\|_{V_0,\ph_0}^2\nonumber\\
& \leq  C \left(\mu_\infty\left((1+|\cdot|)^{2}\right)\mu_\infty(V_0^2)+\mu_\infty((1+|\cdot|)^2V_0^2)\right)\|f\|_{V_0,\ph_0}^2.\nonumber
\end{align}
Hence,   there is $\tld C\geq 1$ depending on $C_Q,\la_Q,\la_P,V_0,\ph_0,F,\mu_\infty$ such that
\beg{align*}
\| Q_t f\|_{V_0,\ph_0}&\leq C_W\left(e^{-\la_Pt}+\tld C\left(\ff {\1_{[\la_P\neq\la_Q]}} {|\la_p-\la_Q|\we 1}+\varphi_1(t)\1_{[\la_P=\la_Q]}\right) e^{-(\la_P\we\la_Q)t}\right)\|f\|_{V_0,\ph_0}\\
&\leq C_W \left(\ff {\1_{[\la_P\neq\la_Q]}} {|\la_p-\la_Q|\we 1}+\tld C\varphi_1(t)\1_{[\la_P=\la_Q]}\right)e^{-(\la_P\we\la_Q)t}\|f\|_{V_0,\ph_0}\\
&\leq C_W\tld C \left(\ff {\1_{[\la_P\neq\la_Q]}} {|\la_p-\la_Q|\we 1 }+ ( t\vee \sq t )\1_{[\la_P=\la_Q]}\right)e^{-(\la_P\we\la_Q)t}\|f\|_{V_0,\ph_0},~f\in\mathscr{G}_{V_0,\ph_0}.
\end{align*}
Similarly, by using \eqref{int-ewe}, we have that 
\[\| Q_t f\|_{V_0,\ph_0}\leq 2C_W\tld C( t\vee \sq t )e^{-(\la_P\we\la_Q)t}\|f\|_{V_0,\ph_0},~f\in\mathscr{G}_{V_0,\ph_0}.\]
Therefore, there is $C\geq 1$ depending on $C_Q,\la_Q,\la_P$, $C_W$, $V_0$, $\ph_0$, $F$, $\mu_\infty$  such that  \eqref{exp-Lip0} holds.

\end{proof}

Next, we use Lemma \ref{Dhr0} and Lemma \ref{exp-con-Q} to estimate $\|\mu_t-\mu_\infty\|_{V_0,\ph_0}$. 
\beg{lem}
Suppose assumptions of Theorem \ref{thm0} hold. Let $H(t)$ be defined in Lemma \ref{exp-con-Q}. Then there exists  $C\geq 1$ depending on $C_Q,\la_Q$, $C_W$, $V_0$, $\ph_0$, $\mu_\infty$ such that  any $\mu_0\in\sP_{U_0}$ and $t\geq 0$, there is
\beg{equation}\label{DuHam-2}
\beg{split}
\|\mu_t-\mu_\infty\|_{V_0,\ph_0}&\leq  C \Big(\|\mu_0-\mu_\infty\|_{V_0,\ph_0}H(t)e^{-(\la_{P}\we \la_Q) t}\\
&\quad\, + \int_0^t e^{-(\la_{P}\we \la_Q)s} \ff {\ph_0(\sq{s\we 1})} {s\we 1}H(s)\|\mu_{t-s}-\mu_\infty\|_{V_0,\ph_0}^2 \d s\Big). 
\end{split}
\end{equation}
\end{lem}
\beg{proof}
We first estimate the second term of \eqref{DuHam}. Combining Lemma \ref{lem-grad-est} (setting $\al_1=0$, $\ph=\ph_0$ and $V_1=V_0$ there) with \eqref{PZf-f}, \eqref{f-L2-ph}, \eqref{Q1=1}, \eqref{nnQ1-pp} and \eqref{nnQ2-pp} with  $c=-f(\xi_t)$, where $\xi_t$ satisfies \eqref{eq-xi} and $\xi_0=x$, we find that there is $C>0$ whose value may vary at each line such that
\beg{align}\label{nnQx-pp}
|\nn Q_tf(x)|&\leq Ce^{Ct}V_0(x)\left(t^{-\ff 1 2}\ph_0\left(t^{\ff 1 2}\right)+1\right)\|f\|_{V_0,\ph_0}\\
|\nn^2 Q_tf(x)|&\leq Ce^{Ct}V_0(x)\left(t^{-1}\ph_0\left(t^{\ff 1 2}\right)+1\right)\|f\|_{V_0,\ph_0},~x\in\R^d,~f\in  \sG_{V_0,\ph_0},\label{nn2Qx-pp}
\end{align}
where in the second inequality, we have used
\beg{align*}
\Th(t)&=\int_0^t\ff {\ph_0(\sq s)} {s\sq{t-s}}\d s \leq \sq{\ff 2 t}\int_0^{t/2}\ff {\ph_0(\sq s)} {s}\d s+\ff {2\ph_0(\sq t)} {t}\int_{t/2}^t\ff {\d s} {\sq{t-s}}\\
&\leq \ff {2\sq 2} {\sq t}\int_0^{\sq{t/2}}\ff {\ph_0(u)} {u}\d s+\ff {2\sq 2 \ph_0(\sq t)} {\sq t}\\
&\leq C \left(1+\ff 1 {\sq t}\right)+\ff {2\sq 2 \ph_0(\sq t)} {\sq t}\\
&\leq C\left(1+\ff {\ph_0(\sq t)} {t}\right).
\end{align*}
Then, the semigroup property of $Q_t$, Lemma \ref{exp-con-Q}, \eqref{nnQx-pp} and \eqref{nn2Qx-pp}  imply that there is  $C\geq 1$ depending on $C_Q,\la_Q$, $C_W$, $p_0$, $\mu_\infty$ such that
\beg{align*}
|\nn Q_tf(x)|&\leq Ce^{-(\la_P\we\la_Q) t}\bar\ph_1(t)V_0(x)  \|f\|_{V_0,\ph_0},\\
|\nn^2 Q_tf(x)|&\leq  Ce^{-(\la_P\we\la_Q) t}\bar\ph_2(t)V_0(x)\|f\|_{V_0,\ph_0},~x\in\R^d,~f\in  \sG_{V_0,\ph_0},~t>0,
\end{align*}
where
\beg{align}\label{barp1p2}
\bar\ph_1(t)&=\ff {\ph_0(\sq{t\we 1})} {\sq{t\we 1}}H(t),\qquad \bar\ph_2(t)=\ff {\ph_0(\sq{t\we 1})} {t\we 1}H(t).
\end{align}
Denote 
\[\la=\la_P\we\la_Q,\qquad q_{t-s}(x)=\nn Q_{t-s}f(x).\] 
Then, together with \eqref{Cr<ph} and \eqref{vV0V0}, as proving \eqref{hhh}, we have that 
\beg{align*}
|q_{t-s}(x_1)-q_{t-s}(x_2)|&\leq  \Big\{\Big[Ce^{-\la (t-s)} \bar\ph_2(t-s)|x_1-x_2|\left(V_0(x_1)+V_0(x_2)\right)\1_{[|x_{1}-x_{2}|\leq 1]}\Big]\\
&\quad\,+\left[Ce^{-\la (t-s)}\bar\ph_1(t-s)\left(V_0(x_1)+V_0(x_2)\right)\1_{[|x_{1}-x_{2}|\geq 1]}\right]\Big\}\|f\|_{V_0,\ph_0}\\
&\leq  C e^{- \la (t-s)}\bar\ph_2(t-s) (|x_1-x_2|\we 1)\left(V_0(x_1)+V_0(x_2)\right)\|f\|_{V_0,\ph_0}\\
&\leq  C e^{- \la (t-s)}\bar\ph_2(t-s)\ph_0(|x_1-x_2|)\left(V_0(x_1)+V_0(x_2)\right)\|f\|_{V_0,\ph_0}.
\end{align*}
Combining this with \eqref{b1-mu-nu} and  \eqref{hhh0}, and denoting  $h(x)=b(x,\mu_s)-b(x,\mu_\infty)$,  we have that
\beg{align*}
&\left|h(x_1)\cdot q_{t-s}(x_1)-h(x_2)\cdot q_{t-s}(x_2)\right|\\
&\quad\,\leq  |(h(x_1)-h(x_2))\cdot q_{t-s}(x_1)|+|h(x_2)\cdot(q_{t-s}(x_1)-q_{t-s}(x_2))| \\
&\quad\,\leq  Ce^{-\la(t-s)}\bar\ph_1(t-s)\ph_0(|x_1-x_2|)V_0(x_1)\|\mu_s-\mu_\infty\|_{V_0,\ph_0}\|f\|_{V_0,\ph_0}\\
&\quad\,\quad\,+ Ce^{-\la (t-s)}  \bar\ph_2(t-s)\|\mu_s-\mu_\infty\|_{V_0,\ph_0} \ph_0(|x_1-x_2|)\left(V_0(x_1)+V_0(x_2)\right)\|f\|_{V_0,\ph_0}\\
&\quad\,\leq 2Ce^{-\la (t-s)}\bar\ph_2(t-s)\ph_0(|x_1-x_2|)\left(V_0(x_1)+V_0(x_2)\right)\|f\|_{V_0,\ph_0}\|\mu_s-\mu_\infty\|_{V_0,\ph_0}.
\end{align*}
Thus
\beg{align*}
(\mu_s-\mu_\infty)((b(\cdot,\mu_s)-b(\cdot,\mu_\infty))\cdot\nn Q_{t-s}f(\cdot))&\leq  Ce^{-\la (t-s)}\bar\ph_2(t-s)\|\mu_s-\mu_\infty\|_{V_0,\ph_0}^2\|f\|_{V_0,\ph_0},
\end{align*}
and as a consequence, 
\beg{equation}\label{5.1-2term}
\beg{split}
&\int_0^t(\mu_s-\mu_\infty)((b(\cdot,\mu_s)-b(\cdot,\mu_\infty))\cdot\nn Q_{t-s}f(\cdot))\d s\\
&\quad\,\leq  C\|f\|_{V_0,\ph_0}    \int_0^t e^{-\la (t-s)}\bar\ph_2(t-s)\|\mu_s-\mu_\infty\|_{V_0,\ph_0}^2\d s.
\end{split}
\end{equation}

Next, we estimate the third term in \eqref{DuHam}.  It follows from \eqref{DFb-pph} that
\beg{align*}
&\|D^F_{\mu_s,\mu_\infty}b(x,\cdot)-D^F_{\mu_\infty}b(x,\cdot)\|_{V_0,\ph_0}\\
&\quad\, \leq \int_0^1 \|D^F_{r\mu_s+(1-r)\mu_\infty}b(x,\cdot)-D^F_{\mu_\infty}b(x,\cdot)\|_{V_0,\ph_0}\d r\\
&\quad\, \leq C\int_0^1 \|r\mu_s+(1-r)\mu_\infty-\mu_\infty\|_{V_0,\ph_0}\d r\\
&\quad\,  =\ff C 2 \|\mu_s-\mu_\infty\|_{V_0,\ph_0}.
\end{align*}
This yields 
\beg{align*} 
&\left|(\mu_s-\mu_\infty)\left( D^F_{\mu_s,\mu_\infty}b(x,\cdot)-D^F_{\mu_\infty}b(x,\cdot)\right)\cdot\nn Q_{t-s}f(x)\right|\\
&\quad\, \leq C\|\mu_s-\mu_\infty\|_{V_0,\ph_0}^2|\nn Q_{t-s}f(x)|\\
&\quad\, \leq  Ce^{-\la (t-s)} \bar\ph_1(t-s) V_0(x)\|f\|_{V_0,\ph_0}\|\mu_s-\mu_\infty\|_{V_0,\ph_0}^2.
\end{align*}
Hence, 
\beg{equation}\label{5.1-3term}
\beg{split}
&\left|\int_0^t \int_{\R^d} (\mu_s-\mu_\infty)\left( D^F_{\mu_s,\mu_\infty}b(x,\cdot)-D^F_{\mu_\infty}b(x,\cdot)\right)\cdot\nn Q_{t-s}f(x) \mu_\infty(\d x)\d s\right|\\
&\quad\, \leq  C \mu_\infty\left(V_0\right) \|f\|_{V_0,\ph_0}\int_0^te^{-\la (t-s)} \bar\ph_1(t-s)\|\mu_s-\mu_\infty\|_{V_0,\ph_0}^2 \d s.
\end{split}
\end{equation}

Substituting \eqref{5.1-2term} and \eqref{5.1-3term} into \eqref{DuHam}, together with Lemma \ref{exp-con-Q} and $\bar\ph_1\leq \bar\ph_2$, there is some $C\geq 1$ depending on $C_Q,\la_Q$, $C_W$, $p_0$, $\mu_\infty$  such that for any $f\in C_b^2(\R^d)$
\beg{equation}\label{DuHam-1}
\beg{split}
(\mu_t-\mu_\infty)(f)&\leq C\|f\|_{V_0,\ph_0}\Big(\|\mu_0-\mu_\infty\|_{V_0,\ph_0}H(t)e^{-\la t}\\
&\quad\, + \int_0^t  e^{-\la (t-s)}\bar\ph_2(t-s)\|\mu_s-\mu_\infty\|_{V_0,\ph_0}^2 \d s\Big).
\end{split}
\end{equation} 

Finally, we extend \eqref{DuHam-1} to $f\in\sG_{V_0,\ph_0}$. Let $\ze \in C^2(\R^d)$ be a nonnegative function such that $\int_{\R^d}\ze(x)\d x=1$ and ${\rm supp}\ze\subset B(0,1)$, and let $\ze_n(x)=n^d\ze(nx)$. For each $f\in \sG_{V_0,\ph_0}$, $m,n\in\N$, we set 
\[f_m=f\vee (-m)\we  m,\qquad f_{m,n}=f_m*\ze_n.\]
Then $f_{m,n}\in C^2_b$ and
\beg{align*}
|f_{m,n}(x)-f_{m,n}(z)|&\leq \int_{\R^d} |f_m(x-y)-f_m(z-y)|\ze_n(y)\d y\\
&\leq \|f\|_{V_0,\ph_0}\ph_0(|x-z|) \int_{\R^d} ( V_0(x-y)+ V_0(z-y))n^d \ze(ny)\d y \\
&\leq \|f\|_{V_0,\ph_0}\ph_0(|x-z|) \int_{|u|\leq 1} \left( V_0\left(x-\ff u n\right)+V_0\left(z-\ff u n\right)\right) \ze(u)\d u\\
&\leq C_{V_0}\|f\|_{V_0,\ph_0}\ph_0(|x-z|)  \left( V_0(x)+V_0(z)\right),
\end{align*}
where we have used \eqref{vV0V0} in the last inequality.  Thus
\[\sup_{m,n\geq 1}\|f_{m,n}\|_{V_0,\ph_0}\leq C_{V_0}\|f\|_{V_0,\ph_0}.\]
Moreover, 
\beg{align*}
|f_{m,n}(x)-f_m(x)|&\leq \int_{\R^d} |f_m(x-y)-f_m(x)|\ze_n(y)\d y\\
&\leq \int_{\R^d} \ph_0(|y|)(  V_0(x-y)+ V_0(x))n^d \ze(ny)\d y\\
&\leq \int_{|u|\leq 1} \ph_0\left(\ff {|u|} n\right)\left( V_0\left(x-\ff u n\right)+V_0(x)\right) \ze(u)\d u,\\
|f_{m,n}(x)|&\leq \int_{\R^d}|f_m(x-y)|\ze_n(y)\d y\\
&\leq |f_m(0)|+\int_{\R^d}|f_m(x-y)-f_m(0)|\ze_n(y)\d y\\
&\leq |f(0)|+\int_{|u|\leq 1} \ph_0\left(|x-\ff u n|\right)\left( V_0\left(x-\ff u n\right)+V_0(0)\right) \ze(u)\d u\\
&\leq |f(0)|+C\ph_0(1+|x|)V_0(x)\\
&\leq |f(0)|+C(1+\ph_0(|x|))V_0(x),
\end{align*}
where  \eqref{vV0V0} has been used in the second-to-last  inequality and that $\ph_0^2(\sq{\cdot})$ is concave in the last inequality. Combining the both inequalities with {\bf (H0)}, \eqref{mu-inf-U0}, \eqref{ph-p-q0},  \eqref{ph-ppp}, $\mu_t\in\sP_{U_0}$ and $|f_m|\leq |f|$, we can derive from the dominated convergence theorem  that
\beg{align*}
\lim_{n\ra +\infty} (\mu_t+\mu_\infty)(|f_{m,n}-f_m|) =0,\\
\lim_{m\ra +\infty} (\mu_t+\mu_\infty)(|f_{m}-f|) =0.
\end{align*}
Hence, we can choose a sequence from $\{f_{m,n}\}_{m,n\geq 1}$, denoting by $\{\tld f_{m}\}_{m\geq 1}$,  such that $\tld f_m\in C^2_b(\R^d)$ and
\beg{align*}
\varlimsup_{m\ra +\infty}\|\tld f_m\|_{V_0,\ph_0}&\leq C_{V_0,\ph_0}\|f\|_{V_0,\ph_0},\\
\lim_{m\ra +\infty}(\mu_t-\mu_\infty)(\tld f_m)&=(\mu_t-\mu_\infty)(f).
\end{align*}
Hence,  \eqref{DuHam-2} follows from \eqref{DuHam-1} and the approximation argument.

\end{proof}

Finally, we prove the convergence of $\|\mu_t-\mu_\infty\|_{V_0,\ph_0}$. Denote $\la=\la_P\we\la_Q$. Let $\bar\ph_2$ be defined by \eqref{barp1p2} and  $C\geq 1$ be a constant such that \eqref{DuHam-2} holds, and let $\mu_0\in\sP_{U_0}$ satisfy \eqref{add-con-init}. Denote $\de=\|\mu_0-\mu_\infty\|_{V_0,\ph_0}(>0)$. 

If $\la_P\neq\la_Q$,  we let
\[\ta=\inf\left\{~t\geq 0~\Big|~\|\mu_t-\mu_\infty\|_{V_0,\ph_0}\geq   \ff {2C} {|\la_p-\la_Q|\we 1} e^{-\la t}\de  \right\}.\]
Then on $[0,\ta]$, there is 
\beg{align*}
&\int_0^t  e^{-\la (t-s)}\ff {\ph_0(\sq{(t-s)\we 1})} {(t-s)\we 1}\|\mu_s-\mu_\infty\|_{V_0,\ph_0}^2 \d s\\
&\quad\,\leq \ff {4 C^2} {(|\la_p-\la_Q|\we 1)^2}\int_0^t  e^{-\la(t-s)}\ff {\ph_0(\sq{(t-s)\we 1})} {(t-s)\we 1}e^{-2\la s}\de^2\d s\\
&\quad\,= \ff {4 C^2} {(|\la_p-\la_Q|\we 1)^2}\de^2e^{-2\la t} \int_0^t e^{\la s}\ff {\ph_0(\sq{s\we 1})} {s\we 1}\d s\\
&\quad\,\leq \ff {4 C^2} {(|\la_p-\la_Q|\we 1)^2}\de^2 e^{-2\la t}\left(\int_0^{t\we 1}e^{\la s}\ff {\ph_0(\sq{s\we 1})} {s\we 1}\d s+\int_1^t e^{\la s}\ff {\ph_0(\sq{s\we 1})} {s\we 1}\d s\1_{[t\geq 1]}\right)\\
&\quad\,\leq \ff {4 C^2} {(|\la_p-\la_Q|\we 1)^2}\de^2 e^{-2\la t}\left(e^{\la (t\we 1)}\int_0^{t\we 1}\ff {\ph_0(\sq{s})} {s}\d s+\ph_0(1)\left(\ff {e^{\la t}- e^{\la} } {\la }\right)\1_{[t\geq 1]}\right)\\
&\quad\,\leq \ff {4 C^2} {(|\la_p-\la_Q|\we 1)^2}\de^2 e^{-2\la t}e^{\la t}\left(\int_0^1\ff {\ph_0(\sq s) } s\d s+\ph_0(1)\left(\ff {1-e^{-\la (t-1)^+}} {\la }\right)\1_{[t\geq 1]}\right)\\
&\quad\,\leq \ff {16 C^2} {(|\la_p-\la_Q|\we 1)^2}\de^2 \left(\left(\int_0^1\ff{\ph_0(u)} u\d u\right)\vee \ff {\ph_0(1) } {2\la}\right) e^{-\la t}\\
&\quad\,=: \ff {16 C^2} {(|\la_p-\la_Q|\we 1)^2}\de^2\de_{\la} e^{-\la t},
\end{align*}
which, together with \eqref{DuHam-2}, yields that
\beg{align*}
&\|\mu_t-\mu_\infty\|_{V_0,\ph_0}\\
&\leq  \ff C {|\la_p-\la_Q|\we 1} \left(\|\mu_0-\mu_\infty\|_{V_0,\ph_0}e^{-\la t}  + \int_0^t  e^{-\la s}\ff {\ph_0(\sq{s\we 1})} {s\we 1}\|\mu_{t-s}-\mu_\infty\|_{V_0,\ph_0}^2 \d s\right)\\
&\leq  \ff {C} { |\la_p-\la_Q|\we 1}\left(1+\ff {16 C^2} {(|\la_p-\la_Q|\we 1)^2}\de\de_{\la}\right) \de e^{-\la t}. 
\end{align*}
For  
\[\de=\|\mu_0-\mu_\infty\|_{V_0,\ph_0}< \ff {(|\la_p-\la_Q|\we 1)^2} {16 C^2\de_{\la}},\]
there is
\beg{equation}\label{mu-ta<}
\|\mu_{\ta}-\mu_\infty\|_{V_0,\ph_0}< \ff {2 C\de } {|\la_p-\la_Q|\we 1} e^{-\la \ta}.
\end{equation}
Due to \eqref{norm-p0q0} and $t\mapsto\mu_t $ is continuous in $\sP_{U_0}$, we have that
\beg{align*}
\varlimsup_{t\ra s}\|\mu_t-\mu_s\|_{V_0,\ph_0}\leq \varlimsup_{t\ra s}\|\mu_t-\mu_s\|_{U_0}=0. 
\end{align*}
Hence, $s\mapsto \|\mu_s-\mu_\infty\|_{V_0,\ph_0}$ is continuous. Consequently, if $\ta<+\infty$, then 
\[\|\mu_{\ta}-\mu_\infty\|_{V_0,\ph_0}= \ff {2 C\de } {|\la_p-\la_Q|\we 1} e^{-\la \ta}.\]
This contradict \eqref{mu-ta<}. Hence, we have that $\ta=+\infty$.  Therefore, 
\[\|\mu_{t}-\mu_\infty\|_{V_0,\ph_0}< \ff {2 C} {|\la_p-\la_Q|\we 1} e^{-\la t}\|\mu_0-\mu_\infty\|_{V_0,\ph_0}.\]

By using 
\[t\vee \sq t\leq 1+t\leq \ff {e^{(\ep\la)\we 1-1}} {(\ep\la)\we 1} e^{\ep\la t}\leq \ff {1} {(\ep\la)\we 1} e^{\ep\la t},~t\geq 0,~\ep>0,\]
we have that, no matter whether $\la_Q=\la_P$ or not,  
\[H(t)e^{-\la t}\leq \ff {e^{-(1-\ep)\la t}} {(\ep\la)\we 1},~t\geq 0,~\ep\in (0,1).\]
Combining this with \eqref{DuHam-2}, there is 
\beg{equation}\label{DuHam-2add}
\beg{split}
&\|\mu_t-\mu_\infty\|_{V_0,\ph_0}\leq  C_{\ep} \Big(\|\mu_0-\mu_\infty\|_{V_0,\ph_0}e^{-(1-\ep)\la t}\\
&+ \int_0^t e^{-(1-\ep)\la(t-s)}\left(\ff {\ph_0(\sq{(t-s)\we 1}} {(t-s)\we 1}\right)\|\mu_s-\mu_\infty\|_{V_0,\ph_0}^2 \d s\Big),~t\geq 0,~\ep\in (0,1),
\end{split}
\end{equation}
where $C_{\ep}= \ff {C} {(\ep\la)\we 1}$.  Repeating the argument for $\la_P\neq\la_Q$, we have that
\[\|\mu_{t}-\mu_\infty\|_{V_0,\ph_0}< 2C_{\ep} e^{-\la(1-\ep) t}\|\mu_0-\mu_\infty\|_{V_0,\ph_0},\]
when $\mu_0$ satisfies
\[\|\mu_0-\mu_\infty\|_{V_0,\ph_0}<\ff 1 {16 C_\ep^2\de_{(1-\ep)\la}}.\]
Note that 
\[\ff 1 {16 C_\ep^2\de_{(1-\ep)\la}}=\ff {[(\ep\la)^2\we 1]} {16C^2}\left(\left(\int_0^1\ff{\ph_0(s)} s\d s\right)^{-1}\we\left(\ff { 2(1-\ep)\la} {\ph_0(1)}\right)\right).\]
Therefore, the proof is complete.

\section{Proofs of examples}

\beg{proof}[Proof of Example \ref{exa2}]

Consider the following equation 
\[\mu(\d x)=\ff { \exp\left\{-\ff 2 {\si^2}\left(\ff {x^4} 4-\ff {x^2} 2+\ff {\be} 2\int_{\R}(x-y)^2\mu(\d y)\right)\right\}\d x } {\int_{\R}\exp\left\{-\ff 2 {\si^2}\left(\ff {x^4} 4-\ff {x^2} 2+\ff {\be} 2\int_{\R}(x-y)^2\mu(\d y)\right)\right\}\d x }.\]
This equation can be reformulated as follows
\beg{equation}\label{imp-mu}
\beg{split}
\mu(\d x)&=\ff { \exp\left\{-\ff 2 {\si^2}\left(\ff {x^4} 4-\ff {x^2} 2+\ff {\be} 2(x^2-2mx)\right)\right\}\d x } {\int_{\R}\exp\left\{-\ff 2 {\si^2}\left(\ff {x^4} 4-\ff {x^2} 2+\ff {\be} 2(x^2-2mx)\right)\right\}\d x }\\
&=\ff { \exp\left\{-\ff 2 {\si^2}\left(\ff {x^4} 4-\ff {x^2} 2+\ff {\be} 2(x-m)^2\right)\right\}\d x } {\int_{\R}\exp\left\{-\ff 2 {\si^2}\left(\ff {x^4} 4-\ff {x^2} 2+\ff {\be} 2(x-m)^2\right)\right\}\d x },
\end{split}
\end{equation}
where $m=\int_{\R}y\mu(\d y)$. The existence of solutions for \eqref{imp-mu} can be established by solving the following equation
\beg{align}\label{ps=0}
\ps(m):=\ff {\int_{\R} (x-m) \exp\left\{-\ff 2 {\si^2}\left(\ff {x^4} 4-\ff {x^2} 2+\ff {\be} 2(x-m)^2\right)\right\}\d x } {\int_{\R}\exp\left\{-\ff 2 {\si^2}\left(\ff {x^4} 4-\ff {x^2} 2+\ff {\be} 2(x-m)^2\right)\right\}\d x }=0.
\end{align}
It is clear that $\ps(0)=0$. According to \cite[Theorem 3.3.1 and Theorem 3.3.2]{Daw} or \cite[Theorem 4.6]{HT10a}, there are $m_+>0$ and $m_-<0$ such that $\ps(m_{\pm})=0$. Due to \cite[Theorem 3.3.1 and Theorem 3.3.2]{Daw}, for $\si<\si_c$, $\ps'(0)>0$ and $\ps(\cdot)$ is concave on $[0,+\infty)$. Thus, there is $\ps'(m_{\pm})<0$. Note that $\ps'(m_{\pm})=-1+\ff {2\be} {\si^2}\mu_{\pm}((\cdot-m_{\pm})^2)$. Thus
\beg{equation}\label{mu+-2}
\ff {2\be} {\si^2}\mu_{ +}((\cdot-m_{+})^2)<1,\quad \ff {2\be} {\si^2}\mu_{-}((\cdot-m_{-})^2)<1,
\end{equation}
where $\mu_{\pm}$ are  solutions of \eqref{imp-mu} associated with $m_{\pm}$.

In this example, $b(x,\mu)=-x^3+x-\be(x-\int_{\R}y\mu(\d y))$. For simplicity, we consider $m_+$, and it is similar for $m_-$.  We first prove the convergence in the $\WW_1$ distance.  

Since $\mu_+$ satisfies \eqref{imp-mu} with $m=m_+$, {\bf (H0)} holds. It is clear that  {\bf (H1)}-{\bf (H4)} hold with $U_0(x)=(1+|x|^2)^{\ff 1 2}$, $K(\mu)\equiv (1-\be)^+\vee 6$, $\be_1=1$,  $V_0\equiv 1$, $\ph_0(r)=r$, $K_V=\be\vee \sq 2$, $\|\mu-\nu\|_{V_0,\ph_0}=\WW_1(\mu,\nu)$, 
\beg{align}\label{lip}
\mathscr{G}_{V_0,\ph_0}&=\mathscr{L}ip:=\left\{f~\Big|~\|f\|_{\mathscr{L}ip}:=\ff {|f(x)-f(y)|} {|x-y|}\leq 1\right\},\\
D^F_{\mu_{+}}b(x,y)&=\be(y-m_+).\label{DFb-exa1}
\end{align}
Set $\mu_\infty=\mu_+$  and $L_{\mu_\infty}= L_{\mu_+}$ defined as follows 
\[L_{\mu_{+}} f(x)=\ff {\si^2} 2\ff {\d^2 f} {\d x^2}(x)-(x^3-x)\ff {\d f} {\d x}(x)-\be(x-m_+)\ff {\d f} {\d x}(x),~f\in C^2_b.\]
Since Lemma \ref{den-C1-W12}, $L_{\mu_+}$ is essential self-adjoint in $L^2(\mu_+)$, and 
\beg{equation}\label{dirich}
\mu_{+}(gL_{\mu_+} f)=-\mu_{+}\left(\ff {\si^2} 2 g' f'\right),~f,g\in \sD(L_{\mu_+}).
\end{equation}
We also denote by $L_{\mu_+}$ the self-adjoint extension and $P^{\mu_+}_t$ the associated diffusion semigroup.  For the operator $\bar A$, there is 
\[\bar Af(x)=\be(x-m_+)\int_{\R} f'(x)\mu_+(\d x)=-\be(x-m_+)\int_{\R}f(x)\left(\ff {\d } {\d x}\log\ff {\d\mu_+}{\d x}\right)(x)\mu_+(\d x),\]
which can be extended to be a compact operator on $L^2(\mu_+)$ since $\left(\ff {\d } {\d x}\log\ff {\d\mu_+}{\d x}\right)\in L^2(\mu_+)$.

According to \cite[Theorem1.4]{LW}, there are $C>0$ and $\la_P>0$ such that 
\beg{equation*}
\WW_1(\left(P_t^{\mu_m}\right)^*\de_x,\left(P_t^{\mu_m}\right)^*\de_y)\leq Ce^{-\la_P t}|x-y|,~x,y\in\R.
\end{equation*}
Due to the Kantorovich-Rubinstein duality or Remark \ref{rem-WW-pp}, \eqref{con-W1} holds for $P^{\mu_+}_t$.

We prove that \eqref{Q-exp-dec} is satisfied. According to \cite[Proposition 3.21]{CGWW},   the supper Poincar\'e inequality holds for  the Dirichlet form associated with  $L_{\mu_+}$. Thus, \eqref{Poincare} holds.  It follows from \eqref{DFb-exa1} and the definition of $\bar A$  that we can set 
\[W(x,y)=\ff {2\be} {\si^2}(x-m_+)(y-m_+).\]
Then $W$ satisfies \eqref{WLW=0}, \eqref{LW-L2} and \eqref{DF=si2W}. Moreover, since \eqref{mu+-2},
\beg{align*}
\int_{\R^d\times\R^d}W(x,y)^2\mu_+(\d x)\mu_+(\d y)=\left(\ff {2\be} {\si^2}\mu_+((\cdot-m_+)^2)\right)^2<1.
\end{align*}
Hence, \eqref{Q-exp-dec} holds according to  Proposition \ref{lem-2.3}.

Therefore, Theorem \ref{thm0} and the Kantorovich-Rubinstein duality imply that the assertion of this example holds.

Finally, we set $U_0(x)=V_0(x)=(1+|x|^2)^{\ff {p_0} 2}$ and  prove the convergence in $\|\cdot\|_{p_0,\ph_0}$ for  $p_0\geq 1,\ph_0(r)=r\we 1$. In this case, it can be directly checked that {\bf (H3)} and {\bf (H4)} hold.  According to the proof above, we remain to prove that \eqref{con-W1}, then Theorem \ref{thm0} can be applied to derive the convergence under $\|\cdot\|_{p_0,\ph_0}$. To this end, we use \cite[Corollary 2.3]{EGZ} with the Lyapunov function $V(x)=(1+x^2)^{\ff {p_0} 2}$.  Then, there exist $\la'>0,\ep>0$ and 
\[\rh(x,y)=f(|x-y|)(1+\ep V(x)+\ep V(y)),~x,y\in\R^d,\]
where $f$ is some nondecreasing, bounded, and concave continuous function satisfying
\beg{equation}\label{add-rfr}
C_1r\leq f(r)\leq C_2r,~r\in [0,R]
\end{equation}
for some positive constants $C_1,C_2,R$ (see \cite[(5.4)]{EGZ}, such that
\[\WW_{\rh}(\sL_{X^{\mu_m}_t(x)},\sL_{X^{\mu_m}_t(y)})\leq e^{-\la' t}\rh(x,y),~x,y\in\R^d,\]
where
\[\WW_{\rh}(\mu,\nu)=\inf_{\pi\in\sC(\mu,\nu)}\int_{\R^d\times\R^d}\rh(x,y)\pi(\d x,\d y),~\mu,\nu\in\sP_{p_0}.\]
Notice that \eqref{add-rfr} and the definition of $\rh(x,y)$ imply that there exist positive constants $\tld C_1, \tld C_2$ such that 
\[\tld C_1\ph_0(|x-y|)((1+|x|)^{p_0}+(1+|y|)^{p_0})\leq \rh(x,y)\leq \tld C_2\ph_0(|x-y|)((1+|x|)^{p_0}+(1+|y|)^{p_0}).\]
Then, \eqref{WW-p-ph} holds for $\la_P=\la'$ and some $C_{W}\geq 1$. According to Remark \ref{rem-WW-pp}, \eqref{con-W1} holds.

\end{proof}

\bigskip

\beg{proof}[Proof of Example \ref{exa1}]
To obtain stationary distributions of \eqref{exa-Gaus}, we only need to solve the following equation
\beg{equation}\label{fix-exa1}
\mu(\d x)=\ff {\exp\left\{-\ff {x^2} 2+\be x \int_{\R}(\cos z)\mu(\d z)\right\}\d x} {\int_{\R}\exp\left\{-\ff {x^2} 2+\be x \int_{\R}(\cos z)\mu(\d z)\right\}\d x }.
\end{equation}
Let 
\[\mu_m(\d x)= \ff {e^{-\ff 1 2(x-\be m)^2}} {\sq{2\pi}}\d x.\]
Then $\mu_m$ is a solution of \eqref{fix-exa1} if and only if $m$ satisfies $\mu_m(\cos(\cdot))=m$.  Since $\mu_m$ is a Gaussian measure and the first equation in \eqref{eq-mm}, there is 
\[\mu_{m}(\cos(\cdot))=e^{-\ff 1 2}\cos(m\be)=m.\]  Thus for $m$ satisfying the first equation in \eqref{eq-mm}, $\mu_m$ is a stationary distribution of \eqref{exa-Gaus}.

In this example, $b(x,\mu)=-x+\be \mu(\cos(\cdot))$. We first prove the convergence in the $\WW_1$ distance.   It is clear that {\bf (H1)}-{\bf (H4)} hold with $U_0(x)=(1+|x|^2)^{\ff 1 2}$, $V_0\equiv1$, $K(\mu)\equiv 0$, $\be_1=0$, $\ph_0(r)=r$, $K_V=\be\vee 1$, $\|\mu-\nu\|_{V_0,\ph_0}=\WW_1(\mu,\nu)$, and $\mathscr{G}_{V_0,\ph_0}=\mathscr{L}ip$ (see \eqref{lip}). It is clear that {\bf (H0)} holds for $\mu_\infty:=\mu_m$. Taking into account $\mu_m(\cos(\cdot))=m$, we have that 
\beg{align*}
L_{\mu_\infty}f&\equiv L_{m}f:=\ff {\d^2 f } {\d x^2}-x\ff {\d f} {\d x}+\be m\ff {\d f} {\d x},~f\in C_b^2,\\
D^F_{\mu_m}b(x,z)&=\be(\cos z-\mu_{m}(\cos(\cdot)))=\be(\cos z-m)\\
\bar Af(x)&=\be (\cos x-m)\mu_{m}(f'),~f\in W^{1,2}_{\mu_m}.
\end{align*}
It follows from the integration by part formula that 
\[\bar Af(x)=\be(\cos x-m)\int_{\R}(x-\be m)f(x)\mu_m(\d x),~f\in W^{1,2}_{\mu_m}.\]
Thus $\bar A$ can be extended to a compact operator on $L^2(\mu_m)$. 

Since Lemma \ref{den-C1-W12}, $L_m$ is an essential self-adjoint operator in $L^2(\mu_m)$, and we also denote by $L_m$ the self-adjoint extension. Moreover, \eqref{Poincare}  holds for $L_m$. We denote by $P_t^{\mu_m}$ the diffusion semigroup generated by $L_m$. Then $P_t^{\mu_m}$ is symmetric w.r.t. $\mu_m$ and  
\beg{equation}\label{exa1-W1}
\WW_1(\left(P_t^{\mu_m}\right)^*\de_x,\left(P_t^{\mu_m}\right)^*\de_y)\leq e^{-t}|x-y|,~x,y\in\R.
\end{equation}
Then,  following from the Kantorovich-Rubinstein duality or Remark \ref{rem-WW-pp}, there is  
\[\|P_t^{\mu_m}f\|_{\mathscr{L}ip}\leq e^{-t}\|f\|_{\mathscr{L}ip},~t\geq 0,~f\in \mathscr{L}ip.\] 

Next, we show the semigroup $Q_t$ generated by $L_{m}+\bar A$ satisfies \eqref{Q-exp-dec}. According to Remark \ref{rem2.8} and Proposition \ref{prop-22}, it sufficients to prove that $I+\bar A(L_{m}-z)^{-1}$ has a bounded  inverse on $\cH=\{f\in L^2_{\mathbb{C}}(\mu_m)~|~\mu_m(f)=0\}$  for any $z$ with ${\rm Re}z\geq 0$.  For simplicity, we denote $e_1(x)=x-\be m$ and $e_\infty(x)=\cos x-m$. Then $e_1$ is the eigenfunction of $L_m$ associated with the eigenvalue $-1$. For any $f\in\cH$ and $z$ with ${\rm Re}z\geq 0$, we have that $(L_{m}-z)^{-1}f\in\sD(L_{m})\subset W^{1,2}_{\mu_m}$ and 
\beg{align*}
\bar A(L_{m}-z)^{-1}f(x) & =\be e_\infty(x)\mu_m(((L_{m}-z)^{-1}f)') = \be e_\infty(x)\mu_m(e_1(L_{m}-z)^{-1}f)\\
& = \be e_\infty(x)\left\< (L_{m}-z)^{-1}f,e_1\right\>_{L^2_{\mathbb{C}}}= \be e_\infty(x)\left\<f,(L_{m}-\bar z)^{-1} e_1\right\>_{L^2_{\mathbb{C}}}\\
&= -(1+z)^{-1}\be e_\infty(x)\left\<f,e_1\right\>_{L^2_{\mathbb{C}}}.
\end{align*}
Then, $I+\bar A(L_{m}-z)^{-1}$ has a bounded  inverse if and only if 
\beg{equation}\label{eq-z-ei-e1}
1-(1+z)^{-1}\be\<e_\infty,e_1\>_{L^2_{\mathbb{C}}}\neq 0.
\end{equation}
It follows from $\cos(\be m)=m\sq e$ (the first equality in \eqref{eq-mm}) that
\beg{align*}
\left\<e_\infty,e_1\right\>_{L^2_{\mathbb{C}}}&=\mu_m(e_1e_\infty)=\int_{\R} x\cos x\mu_m(\d x)-\be m^2\\
&= e^{-\ff 1 2}\left(\be m\cos(\be m)-\sin(\be m)\right)-\be m^2\\
&= \be m^2-e^{-\ff 1 2}\sin(\be m)-\be m^2\\
&=-e^{-\ff 1 2} \sin(\be m).
\end{align*}
Thus \eqref{eq-z-ei-e1} holds for any $z$ with ${\rm Re}z\geq 0$ if and only if 
\[1+z+\ff {\be} { \sq e}\sin(\be m)\neq 0,~{\rm Re}z\geq 0,\]
which follows from the second inequality of \eqref{eq-mm}.

Therefore, Theorem \ref{thm0} yields that \eqref{th-loc-cov} is satisfied for $\mu_m$.  The exponential convergence  in $\WW_1$ follows from the Kantorovich-Rubinstein duality.

The finial assertion can be proved  by using \cite[Corollary 2.3]{EGZ} with the Lyapunov function $V(x)=(1+x^2)^{\ff {p_0} 2}$ and Remark \ref{rem-WW-pp} as in Example \ref{exa2}, and we omit the proof.

\end{proof}

\bigskip

\beg{proof}[Proof of Example \ref{exa3}]
Let $(m_1,m_2)$ be a solution of \eqref{XY-m1m2}, and
\beg{align*}
\mu_m(\d x,\d y)&=\ff 1 {2\pi}\exp\left\{-\ff 1 2\left((x-\be m_2)^2+(y-\be m_1)^2\right)\right\}\d x\d y\\
&=\ff {\exp\{-\ff 1 2(x^2+y^2)+\be m_2 x+\be m_1 y\}} {\int_{\R^2}\exp\{-\ff 1 2(x^2+y^2)+\be m_2 x+\be m_1 y\}\d x\d y }\d x\d y.
\end{align*}
Then
\beg{align*}
\int_{\R^2}\cos(x)\mu_m(\d x,\d y)&=\ff 1 {\sq e}\cos(\be m_2)=m_1,\\
\int_{\R^2}\cos(y)\mu_m(\d x,\d y)&=\ff 1 {\sq e}\cos(\be m_1)=m_2.
\end{align*}
Thus, $\mu_m$ is a stationary distributions of \eqref{exa-XY}.

In this example, $\mu_\infty=\mu_m$, 
\beg{align*}
L_{\mu_\infty}=L_{m}&:=-(x-\be m_2)\ff {\pp } {\pp x}-(y-\be m_1)\ff {\pp } {\pp y}+\ff {\pp^2 } {\pp x^2}+\ff {\pp^2 } {\pp y^2},\\
D_{\mu_m}^Fb(u,v,x,y)&=\be (\cos x-m_2,\cos y-m_1),~x,y,u,v\in\R,\\
\bar A f(x,y)&=\be\int_{\R^2}D_{\mu_m}^Fb(u,v,x,y)\cdot \nn f(u,v) \mu_m(\d u,\d v)\\
&= \be (\cos y-m_2)\mu_m\left(\ff {\pp f} {\pp x}\right)+\be (\cos x-m_1)\mu_m\left(\ff {\pp f} {\pp y}\right),~f\in C_b^1.
\end{align*}
It follows from Lemma \ref{den-C1-W12} that $L_m$ is an essential self-adjoint operator in $L^2(\mu_m)$. According to the integration by part formula, 
\beg{align*}
\mu_m\left(\ff {\pp f} {\pp x}\right)=\int_{\R^2}f(u,v)(u-\be m_2)\mu_m(\d u,\d v),\\
\mu_m\left(\ff {\pp f} {\pp y}\right)=\int_{\R^2}f(u,v)(v-\be m_1)\mu_m(\d u,\d v).
\end{align*}
Thus, $\bar A$ is a compact operator on $L^2(\mu_m)$. 

The rest proof is similar to Example \ref{exa1}, and we  only  prove that  $I+\bar A(L_{m}-z)^{-1}$ has a bounded inverse on $\cH$ for any $z$ with ${\rm Re}z\geq 0$. Let 
\beg{equation*}
\beg{array}{ll}
e_{11}(x,y)=x-\be m_2, &e_{12}(x,y)=y-\be m_1,\\
e_{\infty 1}(x,y)=\cos x- m_1, &e_{\infty 2}(x,y)=\cos y-m_2,~x,y\in\R.
\end{array}
\end{equation*}
Then $e_{11}$ and $e_{12}$ are two linearly independent eigenfunctions of  $L_{m}$ associated with the eigenvalue $-1$. For any $f\in\cH$ and $z$ with ${\rm Re}z\geq 0$,  there are $(L_{m}-z)^{-1}f\in\sD(L_{m})\subset W^{1,2}_{\mu_m}$, and
\beg{align*}
\bar A(L_{m}-z)^{-1}f&=\be e_{\infty 2}\mu_m(e_{11}(L_{m}-z)^{-1}f)+\be e_{\infty 1}\mu_m(e_{12}(L_{m}-z)^{-1}f)\\
&=\be e_{\infty 2}\< (L_{\mu_m}-z)^{-1}f,e_{11}\>_{L^2_{\mathbb{C}}}+\be e_{\infty 1}\<(L_{m}-z)^{-1}f,e_{12}\>_{L^2_{\mathbb{C}}}\\
&=\be e_{\infty 2}\<f, (L_{m}-\bar z)^{-1}e_{11}\>_{L^2_{\mathbb{C}}}+\be e_{\infty 1}\<f,(L_{m}-\bar z)^{-1} e_{12}\>_{L^2_{\mathbb{C}}}\\
&=-\ff {\be} {1+z}\left(e_{\infty 2}\<f,e_{11}\>_{L^2_{\mathbb{C}}}+e_{\infty 1}\<f, e_{12}\>_{L^2_{\mathbb{C}}}\right).
\end{align*}
Then $I+\bar A(L_{m}-z)^{-1}$ has a bounded inverse on $\cH$ if and only if the following matrix is invertible 
\[
\left[
\beg{array}{cccc} 
1  & 0 & 0 & 0\\
0 & 1 & 0 &0\\
0 &-\ff {\be} {1+z}\|e_{12}\|_{L^2_{\mathbb{C}}}^2 & 1-\ff {\be} {1+z}\<e_{\infty 1},e_{12}\>_{L^2_{\mathbb{C}}} & -\ff {\be} {1+z}\<e_{\infty 2},e_{12}\>_{L^2_{\mathbb{C}}}\\
-\ff {\be} {1+z}\|e_{11}\|_{L^2_{\mathbb{C}}}^2 & 0 &  -\ff {\be} {1+z}\<e_{\infty 1},e_{11}\>_{L^2_{\mathbb{C}}} & 1-\ff {\be} {1+z}\<e_{\infty 2},e_{11}\>_{L^2_{\mathbb{C}}}
\end{array}
\right], 
\]
that is also equivalent to 
\beg{equation}\label{det-exa3}
\left|
\beg{array}{cc} 
 1-\ff {\be} {1+z}\mu_{m}(e_{\infty 1}e_{12}) & -\ff {\be} {1+z}\mu_m(e_{\infty 2}e_{12})\\
 -\ff {\be} {1+z}\mu_m(e_{\infty 1}e_{11}) & 1-\ff {\be} {1+z}\mu_m(e_{\infty 2}e_{11})
\end{array}
\right|\neq 0.
\end{equation}
Note that $\mu_m(e_{\infty 2}e_{11})=\mu_m(e_{\infty 1}e_{12})=0$, 
\beg{align*}
\mu_m(e_{\infty 1}e_{11})&=\int_{\R} x\cos x\mu_m(\d x,\d y)-\be m_1m_2\\
&= e^{-\ff 1 2}\left(\be m_2\cos(\be m_2)-\sin(\be m_2)\right)-\be e^{-\ff 1 2} m_2\cos(\be m_2)\\
&=-e^{-\ff 1 2} \sin(\be m_2),
\end{align*}
where we have used \eqref{XY-m1m2} in the second equality, and $\mu_m(e_{\infty 2}e_{12})=-e^{-\ff 1 2} \sin(\be m_1)$.  Then \eqref{det-exa3} holds for any $z$ with ${\rm Re}z\geq 0$ if and only if 
\beg{align*}
(1+z)^2-\ff {\be^2} e\sin(\be m_1)\sin(\be m_2)\neq 0,~{\rm Re}z\geq 0,
\end{align*}
which follows from \eqref{bm1m2}.

\end{proof}

\beg{proof}[Proof of Example \ref{exa-2.5}]
Let $u_0(x)=\ff {x} {(1+x^2)^{\ff 1 3}}$, $m\in\R$ and let 
\beg{align*}
\mu_m(\d x)&=\ff {\exp\left\{-\ff 1 {2\si^2} \left(u_0^2(x)-1\right)^2  - \ff {\be } {\si^2 }\left(u_0^2(x)-2u_0(x)m \right)+\log \ff {1+\ff 1 3 x^2} {(1+x^2)^{\ff 4 3}} \right\}\d x } {\int_{\R} \exp\left\{-\ff 1 {2\si^2} \left(u_0^2(x)-1\right)^2 - \ff {\be } {\si^2 }\left(u_0^2(x)-2u_0(x)m \right)+\log \ff {1+\ff 1 3 x^2} {(1+x^2)^{\ff 4 3}} \right\}\d x }.
\end{align*}
Let $\si_c$ be the same as in Example \ref{exa2}. Then for $\si\in (0,\si_c)$, there are three solutions of \eqref{ps=0}. Then for $m$ being a solutions of \eqref{ps=0}, there is 
\beg{align*}
\int_{\R}u_0(x)\mu_m(\d x) &= \ff {\int_{\R} u_0(x)\exp\left\{-\ff 1 {2\si^2} \left(u_0^2(x)-1\right)^2  - \ff {\be } {\si^2 }\left(u_0^2(x)-2u_0(x)m \right)\right\}\d u_0(x) } {\int_{\R} \exp\left\{-\ff 1 {2\si^2} \left(u_0^2(x)-1\right)^2 - \ff {\be } {\si^2 }\left(u_0^2(x)-2u_0(x)m \right)\right\}\d u_0(x)}\\
&=\ff {\int_{\R} x\exp\left\{-\ff 1 {2\si^2} \left(x^2-1\right)^2  - \ff {\be } {\si^2 }\left(x^2-2x m \right)\right\}\d x } {\int_{\R} \exp\left\{-\ff 1 {2\si^2} \left(x^2-1\right)^2 - \ff {\be } {\si^2 }\left(x^2-2x m \right)\right\}\d x}\\
&=m.
\end{align*}
Thus, $\mu_m$ satisfies  the following equation
\beg{align*}
\mu(\d x)&=\ff {\exp\left\{-\ff {\left(u_0^2-1\right)^2} {2\si^2} + \ff {\be } {\si^2 }\left(u_0^2(x)-2u_0(x)\int_{\R}u_0(y)\mu(\d y) \right)+\log \ff {1+\ff 1 3 x^2} {(1+x^2)^{\ff 4 3}} \right\}\d x } {\int_{\R} \exp\left\{-\ff {\left(u_0^2-1\right)^2} {2\si^2} + \ff {\be } {\si^2 }\left(u_0^2(x)-2u_0(x)\int_{\R}u_0(y)\mu(\d y) \right)+\log \ff {1+\ff 1 3 x^2} {(1+x^2)^{\ff 4 3}} \right\}\d x },
\end{align*}
and $\mu_m$ is a stationary distribution of \eqref{equ-2.5}.

Let  $m_+$ be the positive solution of \eqref{ps=0}. Then $\mu_\infty=\mu_{m_+}=\mu_+$ in this case, and
\beg{align*}
L_{\mu_+}f&=\ff {\si^2} {2} \ff {\d^2 f} {\d x^2}(x)+\ff {1+\ff 1 3x^2} {(1+x^2)^{\ff 4 3}}\left(-\ff {x^3 } {(1+x^2) } +\ff {(1-\be)x } {(1+x^2)^{\ff 1 3}}+\be m_+\right)\ff {\d f} {\d x}(x)\\
&\quad\, -\ff {\si^2 x(1+\ff 1 9 x^2)} {(1+\ff 1 3 x^2)(1+x^2)}\ff {\d f} {\d x}(x),~f\in C^2_b,\\
D^F_{\mu}b(x,y)&=\be\ff {1+\ff 1 3x^2} {(1+x^2)^{\ff 4 3}}\left(u_0(y)-\mu(u_0)\right)=\be u_0'(x)\left(u_0(y)-\mu(u_0)\right),\\
\bar Af(x)&=\be\left(u_0(x)-m_+\right)\int_{\R}u_0'(y) f'(y)\mu_+(\d y),~f\in C^1_b,~x,y\in\R.
\end{align*}
We choose $U_0=V_0$.  It is clear that {\bf (H2)}, \eqref{mu-inf-U0}, \eqref{Inequ-nnb1}, \eqref{Inequ-nnb2}, \eqref{ph-p-q0} and \eqref{b1-mu-nu} hold. \eqref{nnb1-Lypu1} holds for $U_0$ and $V_0$. \eqref{vV0V0} holds since 
\beg{align*}
(1+(x+v)^2)^{\ff 2 3}=\left((1+(x+v)^2)^{\ff 1 2}\right)^{\ff 2 3}\leq \left((1+x^2)^{\ff 1 2}+|v|\right)^{\ff 2 3}\leq (1+x^2)^{\ff 1 3}+|v|^{\ff 2 3}.
\end{align*}
There are positive constants $C_0,r_0$ which is independent of $\mu$ such that
\beg{align}\label{bmu-U0}
b(x,\mu)\ff {\d }{\d x} e^{(1+x^2)^{\ff 1 3}}&\leq \left(- \ff {(1+\ff 1 3 x^2)x^4} {(1+x^2)^3}+ \ff {x^2(1+\ff 1 3 x^2)} {(1+x^2)^{\ff 7 3}}+ \ff {2\si^2x^2(1+\ff 1 9 x^2)} {(1+\ff 1 3 x^2)(1+x^2)^{\ff 5 3}}\right) e^{(1+x^2)^{\ff 1 3}}\nonumber\\
&\quad\, + \ff {2\be(1+\ff 1 3x^2)|x|} {3(1+x^2)^2}e^{(1+x^2)^{\ff 1 3}}\mu(|u_0|)\nonumber\\
&\leq \left(- \ff 1 4\1_{[|x|> r_0]}+ C_0\1_{[|x|\leq r_0]}\right) e^{(1+x^2)^{\ff 1 3}} + \ff {2\be e^{(1+x^2)^{\ff 1 3}}} {3(1+x^2)^{\ff 1 2}}\mu(|u_0|).
\end{align}
By using the Jensen inequality and the following Young inequality:
\[st\leq (1+s)\log(1+s)-s+e^t-t-1,~s,t\geq 0,\]
there are $C_1\geq C_0$ and $r_1>r_0$ which are independent of $\mu$ such that 
\beg{align*}
 \ff {2\be e^{(1+x^2)^{\ff 1 3}}} {3(1+x^2)^{\ff 1 2}}\mu(|u_0|)&\leq  \left(1+\ff {2\be e^{(1+x^2)^{\ff 1 3}}} {3(1+x^2)^{\ff 1 2}}\right)\log\left(1+\ff {2\be e^{(1+x^2)^{\ff 1 3}}} {3(1+x^2)^{\ff 1 2}}\right)+e^{\mu(|u_0|)}\\
&\leq \left(1+\ff {2\be e^{(1+x^2)^{\ff 1 3}}} {3(1+x^2)^{\ff 1 2}}\right)\left((1+x^2)^{\ff 1 3}+\log(1+\ff {2\be} 3)\right)+ \mu\left(e^{|u_0|}\right)\\
 &\leq \left(\ff 1 8\1_{|x|>r_1}+C_1\1_{[|x|\leq r_1]}\right)e^{(1+x^2)^{\ff 1 3}}+ \mu(U_0).
\end{align*}
Putting this into \eqref{bmu-U0}, we have that
\beg{align}\label{Lya-e13}
b(x,\mu)\ff {\d }{\d x} e^{(1+x^2)^{\ff 1 3}}&\leq \left(- \ff 1 4\1_{[|x|> r_1]}+ C_1\1_{[|x|\leq r_1]}\right) e^{(1+x^2)^{\ff 1 3}}+\mu(U_0)\nonumber\\
&\leq -\ff 1 4e^{(1+x^2)^{\ff 1 3}}+\left(C_1\vee \ff 1 4\right)e^{(1+r_1^2)^{\ff 1 3}}+\mu(U_0).
\end{align}
Thus, \eqref{nnb1-Lypu0} and \eqref{LypuV0} hold.   It is clear that \eqref{DF-DF} holds with $F(x)=Cu_0'(x)$ for some  constant $C$ depending on $\be,m_+$. Due to that $u_0'$, $u_0''$ are bounded and
\beg{align*}
|u_0(x)-u_0(y)|&\leq |x-y|\we (|u_0(x)|+|u_0(y)|)\leq |x-y|\we ((1+x^2)^{\ff 1 6}+(1+y^2)^{\ff 1 6})\\
&\leq (|x-y|\we 1)((1+x^2)^{\ff 1 6}+(1+y^2)^{\ff 1 6})\leq \ph_0(|x-y|)(V_0(x)+V_0(y)),
\end{align*}
 \eqref{hhh0} follows from Remark \ref{rem-hhh0},  and \eqref{DFb-pph} follows from Remark \ref{rem-DFb-pph}. Hence, {\bf (H0)}-{\bf (H4)} hold.

Due to \eqref{Lya-e13}, $V_0=U_0$, Remark \ref{rem-WW-pp} and \cite[Example 2.1]{Wan23}, we see that \eqref{con-W1} holds. Then we focus on proving \eqref{Q-exp-dec} in the remainder of the  proof.  Since
\beg{align*}
\ff {\d } {\d m}\mu_m(u_0)=\ff {2\be} {\si^2}\mu_m(u_0^2)-\ff {2\be} {\si^2}\mu_m(u_0)^2=\ff {2\be} {\si^2}\mu_m((u_0-\mu_m(u_0))^2),
\end{align*}
we can derive from \eqref{mu+-2} that 
\beg{align}\label{mu+u0-m<1}
\ff {\d } {\d m}\Big|_{m=m_+}\mu_m(u_0)=\ff {2\be} {\si^2}\mu_+((u_0-\mu_+(u_0))^2)=\ff {2\be} {\si^2}\mu_+((u_0-m_+)^2)<1.
\end{align}
According to the definition of $D^F_{\mu}b(x,y)$ and $\bar A$, we can set
\beg{align*}
W(x,y)=\ff {2\be} {\si^2} (u_0(x)-m_+)(u_0(y)-m_+).
\end{align*}
Then $W$ satisfies \eqref{WLW=0}, \eqref{LW-L2} and \eqref{DF=si2W}.  Moreover, it follows from \eqref{mu+u0-m<1} that 
\beg{align*}
\int_{\R^d\times\R^d}W(x,y)^2\mu_{+}(\d x)\mu_{+}(\d y)=\left(\ff {2\be} {\si^2} \mu_{+}((u_0-m_+)^2)\right)^2<1. 
\end{align*}
By using the integration by part formula, we can extend $\bar A$ to a compact operator on $L^2(\mu_+)$.  $L_{\mu_+}$ is essential self-adjoint and satisfies \eqref{dirich}.  According to \cite[(3.21.4) of Proposition 3.21.]{CGWW}, the super Poincar\'e inequality holds for $\mu_+$, which yields \eqref{Poincare}. Hence, \eqref{Q-exp-dec} follows from Proposition \ref{lem-2.3}.

\end{proof}

\section{Appendix}
Let $\ga_1$ be a nonnegative and locally integrable function on $(0,+\infty)$ such that 
\beg{align*}
\Ga_1(t):=\int_0^t \ff {\ga_1(s)\d s} {\sq{(t-s)\we 1}}<+\infty,~t>0,
\end{align*}
and let $\ga_2$ be a positive and nondecreasing function on $[0,+\infty)$.  
\beg{lem}\label{lem-Gron}
Assume that a nonnegative function $\ps$ satisfies
\beg{equation}\label{in-Gron}
\ps(t)\leq \ga_1(t)+\ga_2(t)\int_0^t \ff {\ps(s)\d s} {\sq{(t-s)\we 1}},~t>0.
\end{equation}
Let $B\left(\cdot,\cdot\right)$ be the beta function,
\beg{align}\label{app-phi0}
\varphi_0(t)&=B\left(\ff 1 2,\ff 1 2\right)+2\sq{t-1}\1_{[1\leq t\leq 2]}+t\1_{[t>2]}\\
\Ga_2(t)&=\int_0^t\left(\ga_1(s)+ \ga_2(s)\Ga_1(s)\right)\d s,~t>0.\nonumber
\end{align}
Then
\beg{align}\label{in-App1}
&\ps(t)\leq \left(\ga_1(t)+\ga_2(t)\Ga_1(t)\right)+\ga_2(t)^2\varphi_0(t) \Ga_2(t)\exp\left\{\int_0^t \ga_2(s)^2\varphi_0(s) \d s\right\}.
\end{align}

\end{lem}
\beg{proof}
For any $r>t$, multiplying both side of \eqref{in-Gron} by $\ff 1 {\sq{(r-t)\we 1}}$ and integrating on $[0,r]$, we have that
\beg{equation}\label{int-rtps-0}
\beg{split}
\int_0^r\ff {\ps(t)\d t} {\sq{(r-t)\we 1}}&\leq \int_0^r \ff {\ga_1(t)\d t} {\sq{(r-t)\we 1}}+\int_0^r\ff {\ga_2(t)} {\sq{(r-t)\we 1}}\int_0^t\ff {\ps(s)\d s} {\sq{(t-s)\we 1}}\d t\\
&\leq \Ga_1(r)+\ga_2(r)\int_0^r\left(\int_s^r\ff {\d t} {\sq{(r-t)\we 1}\cdot\sq{(t-s)\we 1}}\right)\ps(s)\d s\\
&=  \Ga_1(r)+\ga_2(r)\int_0^r\varphi_0(r-s)\ps(s)\d s\\
&\leq \Ga_1(r)+\ga_2(r)\varphi_0(r)\int_0^r\ps(s)\d s.
\end{split}
\end{equation}
This, together with \eqref{in-Gron}, yields that 
\beg{equation}\label{int-rtps-1}
\beg{split}
\ps(t)&\leq \ga_1(t)+\ga_2(t)\left(\Ga_1(t)+\ga_2(t)\varphi_0(t)\int_0^t\ps(s)\d s\right)\\
&= \left(\ga_1(t)+\ga_2(t)\Ga_1(t)\right)+\ga_2(t)^2\varphi_0(t)\int_0^t\ps(s)\d s.
\end{split}
\end{equation}
This implies that 
\beg{align*}
\int_0^r\ps(t)\d t&\leq  \int_0^r\left(\ga_1(t)+\ga_2(t)\Ga_1(t)\right)\d t+\int_0^r\left(\ga_2(t)^2\varphi_0(t)\int_0^t\ps(s)\d s\right)\d t\\
&\leq \Ga_2(r)+\ga_2(r)^2\varphi_0(r)\int_0^r\left(\int_0^t\ps(s)\d s\right)\d t,~r>0.
\end{align*}
Then the Gronwall inequality implies 
\[\int_0^t\ps(s)\d s\leq  \Ga_2(t)\exp\left\{\int_0^t \ga_2(s)^2\varphi_0(s) \d s\right\}.\]  
Combining this with \eqref{int-rtps-0} and \eqref{int-rtps-1}, we obtain 
\[\int_0^t\ff {\ps(s)\d s} {\sq{(t-s)\we 1}} \leq \Ga_1(t)+\ga_2(t)\varphi_0 (t)\Ga_2(t)\exp\left\{\int_0^t \ga_2(s)^2\varphi_0(s) \d s\right\},\]
and \eqref{in-App1}.

\end{proof}

\beg{proof}[{Proof of Lemma \ref{lem-Malli}}]
(1)  It follows from \eqref{Inequ-nnZ1} that
\beg{equation}\label{ZZ}
\beg{split}
\<Z(y)-Z(x),y-x\>&=\int_0^1\<\nn_{y-x} Z(x+\th(y-x)),y-x\>\d\th\\
&\leq \left(K_1-\int_0^1K_2(x+\th(y-x))\d\th\right) |y-x|^2\\
&= \left(K_1-\ff 1 {|y-x|}\int_0^{|y-x|}K_2\left(x+\th\ff {y-x} {|y-x|}\right)\d\th\right) |y-x|^2\\
&\leq \left(K_1- K_2^*(|y-x|)\right)|y-x|^2.
\end{split}
\end{equation}
This, together with that $\si$ is Lipschitz,  implies by \cite[Theorem 3.1.1]{LiuRBook} that \eqref{equ-Y} has a unique strong solution.  Proofs of rest inequalities are routine, we omit them or one can consult proofs for estimates of $Y^\ep_t$ below.

(2) It follows from \cite[Theorem 3.1.1]{LiuRBook} and \eqref{Inequ-nnZ1} that the equation 
\bequ\label{equ-et}
\d \et_{t}=\nn_{\et_{t}}Z(Y_{t})\d t+\nn_{\et_{t}}\si(Y_{t})\d B_t,~\et_{0}=v\in\R^d, t\geq 0.
\enqu
has a unique non-explosive strong solution. Moreover, the B-D-G inequality, the It\^o formula and \eqref{Inequ-nnZ1}  imply that for all $p\geq 1$, there is $C_p>0$ depending on $K_1$ and $\|\nn\si\|_\infty$ such that 
\bequ\label{Inequ-et-sup}
\beg{split}
\E\left(\sup_{s\in[0,t]}|\et_{s}|^{2p}+\int_0^tK_2(Y_s)|\et_s|^{2p}\d r\right)\leq C_pe^{C_pt}|v|^{2p}.
\end{split}
\enqu
The proof of that $Y_{t}$ is derivable w.r.t. the initial value along $v\in \R^d$, and $\nn_vY_t$ satisfies \eqref{equ-et}   is similar to the following discussion on Malliavin differentiable of $Y_t$, and we omit it. Moreover, by replacing $K_1$ by $K_1+1$ and $K_2$ by $K_2+1$, the inequality \eqref{Inequ-et-sup} implies \eqref{Eest-nnY-2p}.

Let $Y^\ep$ be the solution of the following equation:
\beg{equation}\label{Yep0}
Y^\ep_t=Y_0+\int_0^t Z(Y^\ep_s)\d s+\int_0^t\si(Y^\ep_s)\d B_s+\ep \int_0^t\si(Y^\ep_s)h'_s\d s.
\end{equation}
We first prove the well-posedness of \eqref{Yep0} and  estimate  $Y^\ep_t$. 
Since $\nn\si$ is bounded,  there is 
\[|\si(x)h'_{s}-\si(y)h'_s|\leq  \|\nn\si\|_\infty |w(t)|\cdot|\nn_v Y_t|\cdot|x-y|.\]
Then, taking into account that $\si$ is Lipschitz and \eqref{ZZ}, \cite[Theorem 3.1.1]{LiuRBook} implies that \eqref{Yep0} has a unique strong solution. \\
Next, we estimate $Y^\ep_t$.  By \eqref{lim-al1al2}, there are positive constants $r_0, C_1,C_2$ such that
\beg{equation}\label{rr_0}
\beg{split}
(1+r)^{(1+\al_1)\al_2}&\leq C_1K^*_2(r)r^2\1_{[r\geq r_0]}+(1+r)^{(1+\al_1)\al_2}\1_{[r< r_0]}\\
&\leq C_1K^*_2(r)r^2+\left(C_1K^*_2(r)r^2+(1+r)^{(1+\al_1)\al_2}\right)\1_{[r\leq r_0]}\\
&\leq C_1K^*_2(r)r^2+C_2.
\end{split}
\end{equation}
For  $\ep\in [0,1]$ and $\tld\ep\in (0,1]$, by using the It\^o formula, \eqref{lim-K2}, \eqref{ZZ} and \eqref{rr_0}, the H\"older inequality and the Jensen inequality imply that 
\beg{equation*}
\beg{split}
\d |Y^\ep_t|^2-2\<Y_t^\ep,\si(Y_t^\ep)\d B_t\>&=2\<Z(Y^\ep_t),Y^\ep_t\>\d t+2\ep\<\si(Y^\ep)h'_t,Y^\ep_t\>\d t+\|\si(Y^\ep_t)\|_{HS}^2\d t\\
&=  2\<Z(Y_t^\ep)-Z(0),Y^\ep_t\>\d t+2\<Z(0),Y_t^\ep\>\d t\\
&\quad\, +2\ep|\si(Y^\ep_t)|\cdot |Y_t^\ep| \cdot |h'_t|\d t+K_{\si}^2(1+|Y^\ep_t|)^{2\al_1}\d t\\
&\leq  2(K_1- K_2^*(|Y_t^\ep|) )|Y_t^\ep|^2\d t+2|Z(0)|\cdot|Y_t^\ep|\d t\\
&\quad\, +2K_\si(1+|Y^\ep_t|)^{\al_1}|Y_t^\ep|\cdot |w(t)\nn_vY^\ep_t|\d t+K_{\si}^2(1+|Y^\ep_t|)^{1+\al_1}\d t\\
&\leq  2(K_1- K_2^*(|Y_t^\ep|) )|Y_t^\ep|^2\d t+\left(\tld\ep^{-1}|Z(0)|^2+\tld\ep |Y_t^\ep|^2\right)\d t\\
&\quad\, +K_\si(1+|Y^\ep_t|)^{1+\al_1}\left(2 |w(t)\nn_vY^\ep_t|\d t+K_{\si}\right)\d t\\
&\leq  (2K_1+\tld\ep- 2K_2^*(|Y_t^\ep|) )|Y_t^\ep|^2\d t+ \tld\ep^{-1}|Z(0)|^2\d t\\
&\quad\, +K_\si\al_2^{-1}\tld\ep (1+|Y^\ep_t|)^{(1+\al_1)\al_2}\d t\\
&\quad\,+\ff {K_\si (\al_2-1)} {\al_2\tld\ep^{\ff 1 {\al_2-1}}} \left(2 |w(t)\nn_vY^\ep_t| +K_{\si}\right)^{\ff {\al_2} {\al_2-1}}\d t\\
&\leq (2K_1+\tld\ep- 2K_2^*(|Y_t^\ep|) )|Y_t^\ep|^2\d t+ \tld\ep^{-1}|Z(0)|^2\d t\\
&\quad\, +K_\si\al_2^{-1}\tld\ep (C_1K_2^*(|Y^\ep_t|)|Y_t^{\ep}|^2+C_2)\d t\\
&\quad\,+\ff {2^{\ff 1 {\al_2-1}} K_\si (\al_2-1)} {\al_2\tld\ep^{\ff 1 {\al_2-1}}}\left( \left(2 |w(t)\nn_vY^\ep_t| \right)^{\ff {\al_2} {\al_2-1}}+K_{\si}^{\ff {\al_2} {\al_2-1}}\right)\d t\\
&\leq \left(C_1(\tld\ep)-C_2(\tld\ep)K^*_2(|Y_t^\ep|)\right)|Y_t^\ep|^2\d t\\
&\quad\, +C_3(\tld\ep)|w(t)\nn_v Y_t|^{\ff {\al_2} {\al_2-1}}\d t+C_4(\tld\ep)\d t,
\end{split}
\end{equation*}
where \eqref{rr_0} has been used in the second-to-last inequality, and 
\beg{align*}
C_1(\tld\ep)&=2K_1+\tld\ep,\qquad C_2(\tld\ep)=2-C_{\si}\al_2^{-1} C_1\tld\ep,\qquad C_3(\tld\ep)=\ff {2^{\ff {\al_2+1} {\al_2-1}} K_\si (\al_2-1)} {\al_2\tld\ep^{\ff 1 {\al_2-1}}},\\
C_4(\tld\ep)&=\tld\ep^{-1}|Z(0)|^2+K_\si C_2\al_2^{-1}\tld\ep+\ff {2^{\ff 1 {\al_2-1}} K_\si^{\ff {2\al_2-1} {\al_2-1}} (\al_2-1)} {\al_2\tld\ep^{\ff 1 {\al_2-1}}}.
\end{align*} 
Then for any $p\geq 1$,  there are $C_{1,p},C_{2,p},C_{3,p}$ such that 
\beg{align}\label{Yep2p}
\d |Y_t^\ep|^{2p}&\leq  p\left(C_1(\tld\ep)-C_2(\tld\ep)K^*_2(|Y_t^\ep|)\right)|Y_t^\ep|^{2p}\d t +pC_3(\tld\ep)|Y_t^\ep|^{2(p-1)}\d t\nonumber\\
&\quad\,+p C_4(\tld\ep)|Y_t^\ep|^{2(p-1)} |w(t)\nn_v Y_t|^{\ff {\al_2} {\al_2-1}}\d t\nonumber\\
&\quad\, +2p|Y_t^\ep|^{2p-2}\<Y_t^\ep,\si(Y_t^\ep)\d B_t\>+2p(p-1)  |Y^{\ep}_t|^{2(p-2)}|\si^*(Y^{\ep}_t)Y^{\ep}_t|^2\d t\nonumber\\
&\leq   p\left(C_1(\tld\ep)-C_2(\tld\ep)K^*_2(|Y_t^\ep|)\right)|Y_t^\ep|^{2p}\d t +\tld\ep^{-(p-1)}C_3(\tld\ep)^{p}\d t+\tld\ep(p-1) |Y_t^\ep|^{2p}\d t\nonumber\\
&\quad\,+\tld\ep (p-1)|Y_t^\ep|^{2p}\d t+\tld\ep^{-(p-1)}C_4(\tld\ep)^{p}|w(t)\nn_v Y_t|^{\ff {\al_2 p} {\al_2-1}}\d t\nonumber\\
&\quad\, +2p|Y_t^\ep|^{2p-2}\<Y_t^\ep,\si(Y_t^\ep)\d B_t\>+2p(p-1)K_\si^2  |Y^{\ep}_t|^{2(p-1)}(1+|Y^{\ep}_t|)^{2\al_1}\d t\nonumber\\
&\leq  \left(pC_1(\tld\ep)+2\tld\ep(p-1)-pC_2(\tld\ep)K^*_2(|Y_t^\ep|)\right)|Y_t^\ep|^{2p}\d t +\tld\ep^{-(p-1)}C_3(\tld\ep)^{p}\d t\nonumber\\
&\quad\,+\tld\ep^{-(p-1)}C_4(\tld\ep)^{p}|w(t)\nn_v Y_t|^{\ff {\al_2 p} {\al_2-1}}\d t+2p|Y_t^\ep|^{2p-2}\<Y_t^\ep,\si(Y_t^\ep)\d B_t\>\nonumber\\
&\quad\, +2p(p-1)K_\si^2  |Y^{\ep}_t|^{2(p-1)}\left(\tld\ep\al_2^{-1}(1+|Y^{\ep}_t|)^{(1+\al_1)\al_2}+\ff {\al_2-1} {\al_2\tld\ep^{\ff 1 {\al_2-1}}} \right)\d t\\
&\leq  \left(pC_1(\tld\ep)+2\tld\ep(p-1)-pC_2(\tld\ep)K^*_2(|Y_t^\ep|)\right)|Y_t^\ep|^{2p}\d t +\tld\ep^{-(p-1)}C_3(\tld\ep)^{p}\d t\nonumber\\
&\quad\,+\tld\ep^{-(p-1)}C_4(\tld\ep)^{p}|w(t)\nn_v Y_t|^{\ff {\al_2 p} {\al_2-1}}\d t+2p|Y_t^\ep|^{2p-2}\<Y_t^\ep,\si(Y_t^\ep)\d B_t\>\nonumber\\
&\quad\, +2p(p-1)K_\si^2  |Y^{\ep}_t|^{2(p-1)}\left(\ff {\tld\ep C_1} {\al_2}K^*_2(|Y^\ep_t|)|Y_t^{\ep}|^2+\ff {\tld\ep C_2} {\al_2}+\ff {\al_2-1} {\al_2\tld\ep^{\ff 1 {\al_2-1}}} \right)\d t\nonumber\\
&\leq  \left(pC_1(\tld\ep)+2\tld\ep(p-1)-pC_2(\tld\ep)K^*_2(|Y_t^\ep|)\right)|Y_t^\ep|^{2p}\d t +\tld\ep^{-(p-1)}C_3(\tld\ep)^{p}\d t\nonumber\\
&\quad\,+\tld\ep^{-(p-1)}C_4(\tld\ep)^{p}|w(t)\nn_v Y_t|^{\ff {\al_2 p} {\al_2-1}}\d t+2p|Y_t^\ep|^{2p-2}\<Y_t^\ep,\si(Y_t^\ep)\d B_t\>\nonumber\\
&\quad\, +\tld\ep\ff {2p(p-1)K_\si^2 C_1} {\al_2}K^*_2(|Y^\ep_t|)|Y_t^{\ep}|^{2p} \d t\nonumber\\ 
&\quad\, +\ff {2(p-1)^2K_\si^2\tld\ep } {\al_2} |Y^{\ep}_t|^{2p}+\ff {2(p-1)K_\si^2} {\al_2\tld\ep^{p-1}}  \left( \tld\ep C_2 +\ff {\al_2-1} { \tld\ep^{\ff 1 {\al_2-1}}} \right)^{p}\d t\nonumber\\
&= \left(C_5(\tld\ep)-C_6(\tld\ep)K^*_2(|Y_t^\ep|)\right)|Y_t^\ep|^{2p}\d t +C_7(\tld\ep)|w(t)\nn_v Y_t|^{\ff {\al_2 p} {\al_2-1}}\d t\nonumber\\
&\quad\,+C_8(\tld\ep)\d t+2p|Y_t^\ep|^{2p-2}\<Y_t^\ep,\si(Y_t^\ep)\d B_t\>,\nonumber
\end{align}
where
\beg{align*}
C_5(\tld\ep)&=pC_1(\tld\ep)+2\tld\ep(p-1)+\ff {2(p-1)^2K_\si^2\tld\ep } {\al_2},\qquad C_6(\tld\ep)=pC_2(\tld\ep)-\tld\ep\ff {2p(p-1)K_\si^2 C_1} {\al_2},\\
C_7(\tld\ep)&=\tld\ep^{-(p-1)}C_4(\tld\ep)^{p},\qquad C_8(\tld\ep)=\tld\ep^{-(p-1)}C_3(\tld\ep)^{p}+\ff {2(p-1)K_\si^2} {\al_2\tld\ep^{p-1}}  \left( \tld\ep C_2 +\ff {\al_2-1} { \tld\ep^{\ff 1 {\al_2-1}}} \right)^{p}.
\end{align*}
Applying the B-D-G inequality and the H\"older inequality, we get that 
\beg{align*}
&\E\left[\sup_{s\in [0,t]}\left|\int_0^s |Y_r^\ep|^{2p-2}\<Y_r^\ep,\si(Y_r^\ep)\d B_r\>\right|\Big|\sF_0\right]\leq C \E\left[\left(\int_0^t |Y_r^{\ep}|^{4p-2}(1+|Y^{\ep}_r|)^2\d r\right)^{\ff 1 2}\Big|\sF_0\right]\\
&\leq C\E\left[\left(\int_0^t \left(|Y_r^{\ep}|^{4p}+1\right)\d r\right)^{\ff 1 2}\Big|\sF_0\right]\leq C\left(\sq{t}+ \E\left[\sup_{r\in [0,t]}|Y_r^{\ep}|^p\left(\int_0^t |Y_r^{\ep}|^{2p}\d r\right)^{\ff 1 2}\Big|\sF_0\right]\right)\\
&\leq C\left(\sq t+\E\left[ \int_0^t|Y_r^{\ep}|^{2p}\d r\Big|\sF_0\right]\right)+\ff 1 {2p}\E\left[ \sup_{r\in [0,t]}|Y_r^{\ep}|^{2p}\Big|\sF_0\right].
\end{align*}
Combining this with \eqref{Yep2p} and \eqref{Inequ-et-sup}, the Gronwall inequality and the stopping time argument imply that there exists $C_p>0$ such that  for any $t>0$
\beg{equation}\label{supYep2p}
\beg{split}
&\sup_{\ep\in [0,1]}\E\left[\sup_{s\in [0,t]}|Y_t^\ep|^{2p}+ \int_0^tK^*_2(|Y_s^\ep|)|Y_s^\ep|^{2p}\d s\Big|\sF_0\right]\\
&\quad\,\leq C_p\left(\sq t(\sq t+1)+|Y_0|^{2p}+\int_0^t\E\left[ |w(s)\nn_v Y_s|^{\ff {\al_2 p} {\al_2-1}}\Big|\sF_0\right]\d s \right)e^{C_{p} t}.
\end{split}
\end{equation}
By replacing $K_1$ by $K_1+1$ and $K_2$ by $K_2+1$, and taking $w(\cdot)\equiv 0$, we obtain \eqref{Eest-Y-2p}.

Next, we estimate $Y_t^\ep-Y_t$. For $h_t$, because of \eqref{Eest-nnY-2p} and that $w$ is a bounded process, for any $s>0$ and $p>0$,  we have that
\beg{equation}\label{Inhsup}
\beg{split}
\E\sup_{t\in [0,s]}|h'_{t}|^p&\leq  \sup_{t\in[0,s],\om\in \Om}|w(t,\om)|\E\sup_{t\in [0,s]}|\nn_v Y_t|^p\\
&\leq  C_pe^{C_p t}(\sup_{t\in[0,s],\om\in\Om}|w(t,\om)|)|v|^p
\end{split}
\end{equation}
for some $C_p>0$. \\
Due to \eqref{ZZ}, (A2) and the It\^o formula, for any $p\geq 1$, there is   
\beg{align*}
\d |Y_t^\ep-Y_t|^{2p}&\leq   p\left(2K_1+ \ff {2p-1} p+(2p-1)\|\nn\si\|_\infty^2- 2K^*_2(|Y_t^\ep-Y_t|)\right)|Y_t^\ep-Y_t|^{2p}\d t\\
&\quad\,+\ep^{2p}|\si(Y_t^\ep)h'_t|^{2p}\d t+2p|Y_t^\ep-Y_t|^{2p-2}\<Y_t^\ep-Y_t,(\si(Y^\ep_t)-\si(Y_t))\d B_t\>.
\end{align*}
Combining this with the B-D-G inequality, the Gronwall inequality and the H\"older inequality,  there is  $C_{p}>0$ such that 
\beg{align*}
&\E\sup_{s\in [0,t]}|Y_s^\ep-Y_s|^{2p}+\E\int_0^tK^*_2(|Y_s^\ep-Y_s|)|Y_s^\ep-Y_s|^{2p}\d s\\
&\quad\,\leq \ep^{2p}C_pe^{C_p t}\E \int_0^t|\si(Y_s^\ep)h'_s|^{2p}\d s.
\end{align*}
This, together with \eqref{supYep2p}, \eqref{Inhsup} and \eqref{Inequ-et-sup}, yields that
\beg{equation}\label{YYep}
\sup_{\ep\in (0,1],|v|\leq 1}\E\sup_{s\in [0,t]}\left|\ff  {Y_s^\ep-Y_s} {\ep}\right|^{2p}<+\infty,~t\geq 0.
\end{equation}

The conditions (A1) and (A2), together with \cite[Theorem 3.1.1]{LiuRBook}, yield that the following equation  has a unique strong solution
\[\d D_h Y_t= \nn_{D_h Y_t}Z(Y_t)\d t+\nn_{D_h Y_t}\si(Y_t)\d B_t+\si(Y_t)h'_t\d t,~D_h Y_0=0.\]
By using the It\^o formula and the B-D-G inequality directly, it can be proved that $D_h Y_t$ satisfies \eqref{Eest-DhY-2p}.
Let 
\[V_t^\ep=\ff {Y_t^\ep-Y_t} {\ep}-D_hY_t,\qquad  U_t^\ep=\ff {Y_t^\ep-Y_t} {\ep}.\]
Then
\beg{align*}
\d V_t^\ep=&\int_0^1\left(\nn_{U_t^\ep}Z(Y_t+\th(Y_t^\ep-Y_t))-\nn_{U_t^\ep}Z(Y_t)\right)\d\th\d t\\
&+\nn_{V_t^\ep}Z(Y_t)\d t+\nn_{V_t^\ep}\si(Y_t)\d B_t+\left(\si(Y_t^\ep)-\si(Y_t)\right)h'_t\d t\\
&+\int_0^1\left(\nn_{U_t^\ep}\si(Y_t+\th(Y_t^\ep-Y_t))-\nn_{U_t^\ep}\si(Y_t)\right)\d \th\d B_t.
\end{align*}
By the It\^o formula, {\bf (A1)} and {\bf (A2)}, and the H\"older inequality, for any $p\geq 1$, there is a positive constant $C$ which depends on $p,K_1,\|\nn\si\|_\infty$ such that
\beg{align*}
\d |V_t^\ep|^{2p}\leq & 2p\left(C- K_2(Y^\ep_t)\right)|V_t^\ep|^{2p}\d t+|\si(Y_t^\ep)-\si(Y_t)|^{2p}|h_t'|^{2p}\d t\\
& +\int_0^1 |\nn\si(Y_t+\th\ep U_t^\ep)-\nn \si(Y_t) |^{2p}\d \th|U_t^\ep|^{2p}\d t\\
& +\int_0^1|\nn Z(Y_t+\ep\th U_t^\ep)-\nn Z(Y_t)|^{2p}\d \th |U_t^\ep|^{2p}\d t\\
& + 2p |V_t^\ep|^{2p-2}\<V_t^\ep,\int_0^1\left(\nn_{U_t^\ep}\si(Y_t+\ep\th U_t^\ep)-\nn_{U_t^\ep}\si(Y_t)\right)\d \th\d B_t\>\\
& +2p |V_t^\ep|^{2p-2}\<V_t^\ep,\nn_{V_t^\ep}\si(Y_t)\d B_t\>.
\end{align*}
It follows from the B-D-G inequality and the Gronwall  inequality that there is a positive constant $C$ which depends on $p,K_1,\|\nn\si\|_\infty$ such that for any $t\geq 0$
\beg{align*}
&\E\sup_{s\in [0,t]}|V_s^\ep|^{2p}+ \E\int_0^t K_2(Y_s^\ep)|V_s^\ep|^{2p}\d s\\
&\quad\, \leq Ce^{Ct}\left( \E\int_0^t\int_0^1 |\nn\si(Y_s+\th\ep U_s^\ep)-\nn \si(Y_s) |^{2p}\d \th|U_s^\ep|^{2p}\d s\right.\\
&\quad\,\quad\, +\E\int_0^t\int_0^1|\nn Z(Y_s+\ep\th U_s^\ep)-\nn Z(Y_s)|^{2p}\d \th |U_s^\ep|^{2p}\d s\\
&\qquad\,\, \left. + \ep^{2p}\|\nn\si\|_\infty^{2p}\E\int_0^t |U_s^{\ep}|^{2p}|h_s'|^{2p}\d s\right)\\
&=: I_1+I_2+I_3.
\end{align*}
Since $\ep\th U_s^\ep\ra 0$ in probability, which is implied by \eqref{YYep}, and $\nn\si(\cdot)$ is continuous, the dominated convergence theorem yields that $\lim_{\ep\ra 0^+}I_1=0$.\\
Since $\nn Z$ has polynomial growth, there exist $q>0$ and $C>0$ such that 
\[|\nn Z(y)|\leq C(1+|y|)^q.\]
Then
\beg{align*}
|\nn Z(Y_s+\ep\th U_s^\ep)-\nn Z(Y_s)|^{2p}|U_s^\ep|^{2p}\leq 2C(|Y_s|+|U^\ep_s|)^{2pq}|U_s^\ep|^{2p},
\end{align*}
which is integrable on $\Om\times [0,t]\times [0,1]$ under $\d \P\otimes\d s\otimes\d\th$. This, together with that $\ep\th U_s^\ep\ra 0$ in probability, implies by using the dominated convergence theorem that $\lim_{\ep\ra 0^+}I_2=0$.\\
By using \eqref{Inhsup} and \eqref{YYep}, it is clear that $\lim_{\ep\ra 0^+}I_3=0$.\\
Hence
\[\lim_{\ep\ra 0^+}\E\sup_{s\in [0,t]}|V_s^\ep|^{2p}=0.\]
Therefore, $Y_t$ is Malliavin differentiable along $h$.

(3) Similarly, following from {\bf {(A1)}}, {\bf {(A2)}} except \eqref{nondege0}, \eqref{Inequ-et-sup} and \eqref{supYep2p},  we can prove that $\nn_v Y_t^y$ is differentiable w.r.t. $y$, i.e. $Y_t$ is twice differentiable w.r.t. the initial value. Moreover, $\nn_u\nn_v Y_t$ satisfies 
\beg{align*}
\d \nn_u\nn_v Y_t= & \nn_{\nn_u Y_t}\nn_{\nn_v Y_t}Z(Y_t)\d t+ \nn_{\nn_u Y_t}\nn_{\nn_v Y_t} \si(Y_t)\d B_t\\
& +\nn_{\nn_u\nn_v Y_t}Z(Y_t)\d t+\nn_{\nn_u\nn_v Y_t}\si(Y_t)\d B_t,~\nn_u\nn_v Y_0=0,
\end{align*}
and the It\^o formula, {\bf {(A1)}} and \eqref{Inequ-nnZ2} imply
\beg{equation}\label{nn2Y2}
\beg{split}
&\d |\nn_u\nn_v Y_t|^2-2\<\nn_{u}\nn_{v}Y_t,\left(\nn_{\nn_u Y_t}\nn_{\nn_vY_t}\si(Y_t)+\nn_{\nn_u\nn_vY_t}\si(Y_t)\right)\d B_t\>\\
&\quad\, = 2\<\nn_{\nn_u Y_t}\nn_{\nn_vY_t}Z(Y_t),\nn_u\nn_vY_t\>\d t+2\<\nn_u\nn_vY_t,\nn_{\nn_u\nn_vY_t}Z(Y_t)\>\d t\\
&\quad\,\quad\, +\|\nn_{\nn_u Y_t}\nn_{\nn_vY_t}\si(Y_t)+\nn_{\nn_u\nn_vY_t}\si(Y_t)\|_{HS}^2\d t\\
&\quad\,\leq 2|\nn^2Z(Y_t)|\cdot |\nn_{u}\nn_{v}Y_t|\cdot |\nn_u Y_t|\cdot|\nn_v Y_t|\d t+2(K_1-K_2(Y_t))|\nn_{u}\nn_{v}Y_t|^2\d t\\
&\quad\, \quad\, +2\left(\|\nn^2\si\|_\infty^2 |\nn_u Y_t|^2|\nn_v Y_t|^2+\|\nn\si\|_\infty^2|\nn_u\nn_vY_t|^2\right)\d t.
\end{split}
\end{equation}
Then for any $p\geq 1$, there are positive constants $C_1,\cdots,C_5$ such that
\beg{align*}
\d |\nn_u\nn_v Y_t^y|^{2p}&\leq p|\nn_u\nn_v Y_t^y|^{2(p-1)}\left(C_1|\nn_u\nn_v Y_t^y|^{2}+C_2(|\nn^2Z(Y_t^y)|^2+1)|\nn_u Y^y_t|^2\cdot|\nn_v Y^y_t|^2\right)\d t\\
&\quad\,+2p|\nn_y\nn_v Y^y_t|^{2(p-1)}\<\nn_u\nn_v Y_t^y,\left(\nn_{\nn_u Y^y_t}\nn_{\nn_v Y^y_t}\si(Y^y_t)+\nn_{\nn_u\nn_v Y_t^y}\si(Y^y_t)\right)\d B_t\>\\
&\quad\,+C_3|\nn_u\nn_v Y_t^y|^{2(p-1)}\left(|\nn_u Y^y_t|^2 |\nn_v Y_t^y|^2+|\nn_u\nn_v Y_t^y|^2\right)\d t\\
&\leq C_4|\nn_u\nn_v Y_t^y|^{2p}\d t+C_5\left(1+|\nn^2Z(Y_t^y)|\right)^{2p}|\nn_u Y^y_t|^{2p} |\nn_v Y_t^y|^{2p}\d t\\
&\quad\, +2p|\nn_y\nn_v Y^y_t|^{2(p-1)}\<\nn_u\nn_v Y_t^y,\left(\nn_{\nn_u Y^y_t}\nn_{\nn_v Y^y_t}\si(Y^y_t)+\nn_{\nn_u\nn_v Y_t^y}\si(Y^y_t)\right)\d B_t\>.
\end{align*}
Then, the B-D-G inequality, the Gronwall inequality and the H\"older inequality imply that there is $C>0$ such that
\beg{align*}
&\E\sup_{t\in [0,s]} |\nn_u\nn_v Y_t^y|^{2p}\leq Ce^{Cs}\E\int_0^s\left(1+|\nn^2Z(Y_t^y)|\right)^{2p}|\nn_u Y^y_t|^{2p} |\nn_v Y_t^y|^{2p}\d t\\
&\quad\,\leq Ce^{Cs}\left(\E\int_0^s\left(1+|\nn^2Z(Y_t^y)|\right)^{4p}\d t\right)^{\ff 1 2}\left(\E\int_0^s|\nn_u Y^y_t|^{8p}\d t\right)^{\ff 1 4}\left(\E\int_0^s |\nn_v Y_t^y|^{8p}\d t\right)^{\ff 1 4}.
\end{align*}
Combining this with \eqref{Inequ-nnZ2}  and \eqref{Eest-Y-2p}, we find that there is a positive constant $C$ which depends on  $p,K_1,\be_1,r_0,K_3,\|\nn\si\|_\infty,\|\nn^2\si\|_\infty$ such that
\bequ\label{nn2Y-be1-uv}
\E\sup_{t\in [0,s]}|\nn_u\nn_v Y_t^y|^{2p}\leq C_pe^{C_ps}(1+|y|)^{2p\be_1}|u|^{2p}|v|^{2p},~s\geq 0.
\enqu
This yields \eqref{Inequ-nn2Y}.

\end{proof}

\noindent\textbf{Acknowledgements}

\medskip

The  author was grateful to the financial support by program for innovation research in Central University of Finance and Economics, and was supported by the National Natural Science Foundation of China (Grant No. 12371153).


\end{document}